\documentclass[12pt]{article}
\usepackage{amsmath,amsthm,amssymb,graphicx}
\usepackage{lineno}
\usepackage{overpic}
\usepackage{fullpage}
\newtheorem{thm}{Theorem}[section]

\theoremstyle{definition}
\newtheorem{defn}[thm]{Definition}

\usepackage{graphicx}
\usepackage{booktabs}
\usepackage{authblk}

\begin{document}

\title{Performance of the Uniform Closure Method for open knotting as a Bayes-type classifier}

\author[1]{Emily Tibor}
\author[2]{Elizabeth M. Annoni}
\author[2]{Erin Brine-Doyle}
\author[2]{Nicole Kumerow}
\author[2]{Madeline Shogren}
\author[3]{Jason Cantarella}
\author[4]{Clayton Shonkwiler}
\author[2]{Eric J. Rawdon}

\affil[1]{Department of Mathematics, University of Minnesota, Minneapolis, MN 55455, USA}
\affil[2]{University of St. Thomas, St. Paul, MN 55105, USA}
\affil[3]{Department of Mathematics, University of Georgia, Athens, GA 30602, USA}
\affil[4]{Department of Mathematics, Colorado State University, Fort Collins, CO 80523, USA}

\date{\today}
\maketitle

\begin{abstract}
  The discovery of knotting in proteins and other macromolecular
  chains has motivated researchers to more carefully consider how to
  identify and classify knots in open arcs. Most definitions 
  classify knotting in open arcs by constructing an ensemble of 
  closures and measuring the probability of different knot types
  among these closures. In this paper, we think of assigning knot
  types to open curves as a classification problem and compare the 
  performance of the Bayes MAP classifier to the standard Uniform
  Closure Method. Surprisingly, we find that both methods are essentially
  equivalent as classifiers, having comparable accuracy and positive
  predictive value across a wide range of input arc lengths and
  knot types.
 \end{abstract}



\section{Introduction}

Many physical systems contain long, open chain-like objects with two
free ends, e.g.~DNA, RNA, and other proteins.  Anyone who has ever
packed away a string of lights, a garden hose, or headphone cables
knows that these sorts of objects tend to be entangled.  However, 
it is not clear how to measure this entanglement mathematically.

Several techniques have been proposed~\cite{openknotting,knotencyclopedia}.
In most (the exception being
knotoids \cite{turaev2012}), entanglement of open chains is measured by
generating an ensemble of closures and measuring their knot types using traditional
tools from mathematical knot theory.

There are a number of proposals for how to generate the closures (for a review, see
\cite{openknotting,knotencyclopedia}) which vary in computation speed. 
In the simplest of cases (such as when the endpoints are close to each other
or lie on the surface of the convex hull), all of the
techniques generally agree.  However, in more ambiguous situations,
such as when the endpoints are well inside the convex hull, the
different techniques can disagree.

These ambiguous situations cannot be neglected. In proteins, for example, it has been proposed that the borders of
the ``knotted cores'' (i.e.~the shortest subchains realizing a
particular knot type) are rich with intra-chain contacts, provide
structural stability, and may be correlated with their active sites
\cite{functionproteins}.  The determination of which subchains
determine the knotted cores requires that one computes the knotting in
all subchains
\cite{stevedore,knotprot2,knotprot,mansfield1,mansfield2,ourpnas,virnauweb},
which necessarily include the difficult-to-classify subchains.  
These papers show that the ``core'' knotting and linking has been
preserved within classes of proteins which have the same function in
different organisms.  The proteins within these classes are separated
by hundreds of millions of years of evolution, suggesting that their
entanglement is critical to their function
\cite{joannapnaslinks,ourpnas}.

While trying to compare different approaches to building ensembles of
closures for open arcs, we were struck by the similarity 
to a problem in machine learning: design a classifier which predicts
the knot type of a closed polygon $P$ given a subarc $A$ of that
polygon (for some existing machine learning approaches to knot classification, see~\cite{Gukov:2020vr,Hughes:2020hh,Jejjala:2019gm,Levitt:2019ut,Vandans:2020cp}). In this paper, we will focus on equilateral closed polygons
chosen from the uniform\footnote{This probability distribution is not
  meant to be a good approximation of the distribution of shapes of
  biological macromolecules (for instance, it ignores steric effects
  and bending stiffness), but it is a mathematically natural object of
  study.} probability measure~\cite{jasoneqknots,Cantarella:2016iy} on such
polygons. For polygons $P$ of various knot types randomly sampled from
this measure, we will compare the performance of the classical Bayes
maximum a posteriori probability (MAP) classifier to the performance of the
Uniform Closure Method, which is one of the most-used definitions for
measuring entanglement in open
chains~\cite{openknotting,knotprot2,knotencyclopedia,andrzejchromatin,fingerprints,subknots,ourpnas}.

\begin{figure}
  \centering
  {\ }\hfill
  \begin{overpic}[angle=0,width=0.32\textwidth]{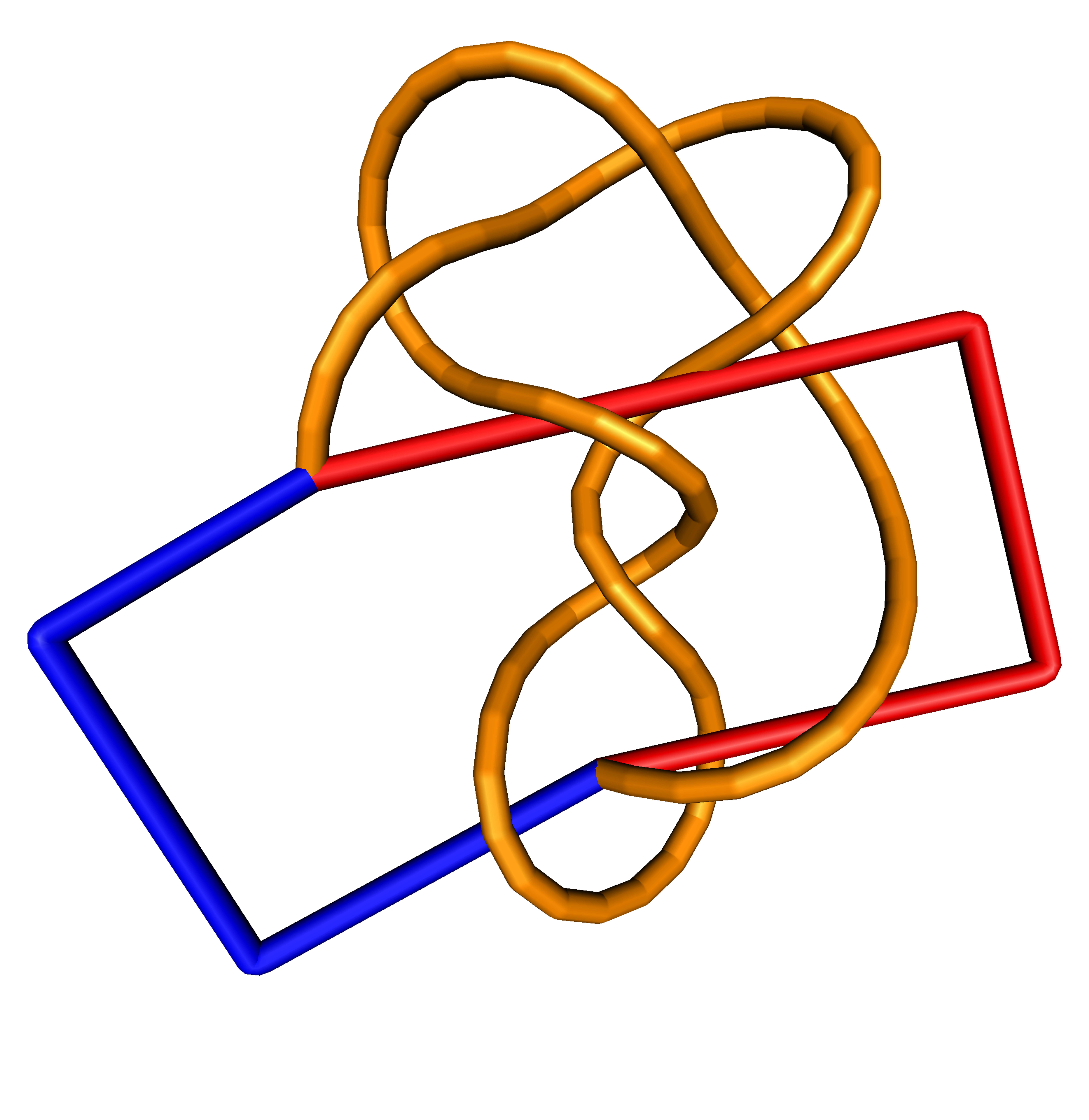}
    \end{overpic}
  \hfill
  \begin{overpic}[angle=0,width=0.32\textwidth]{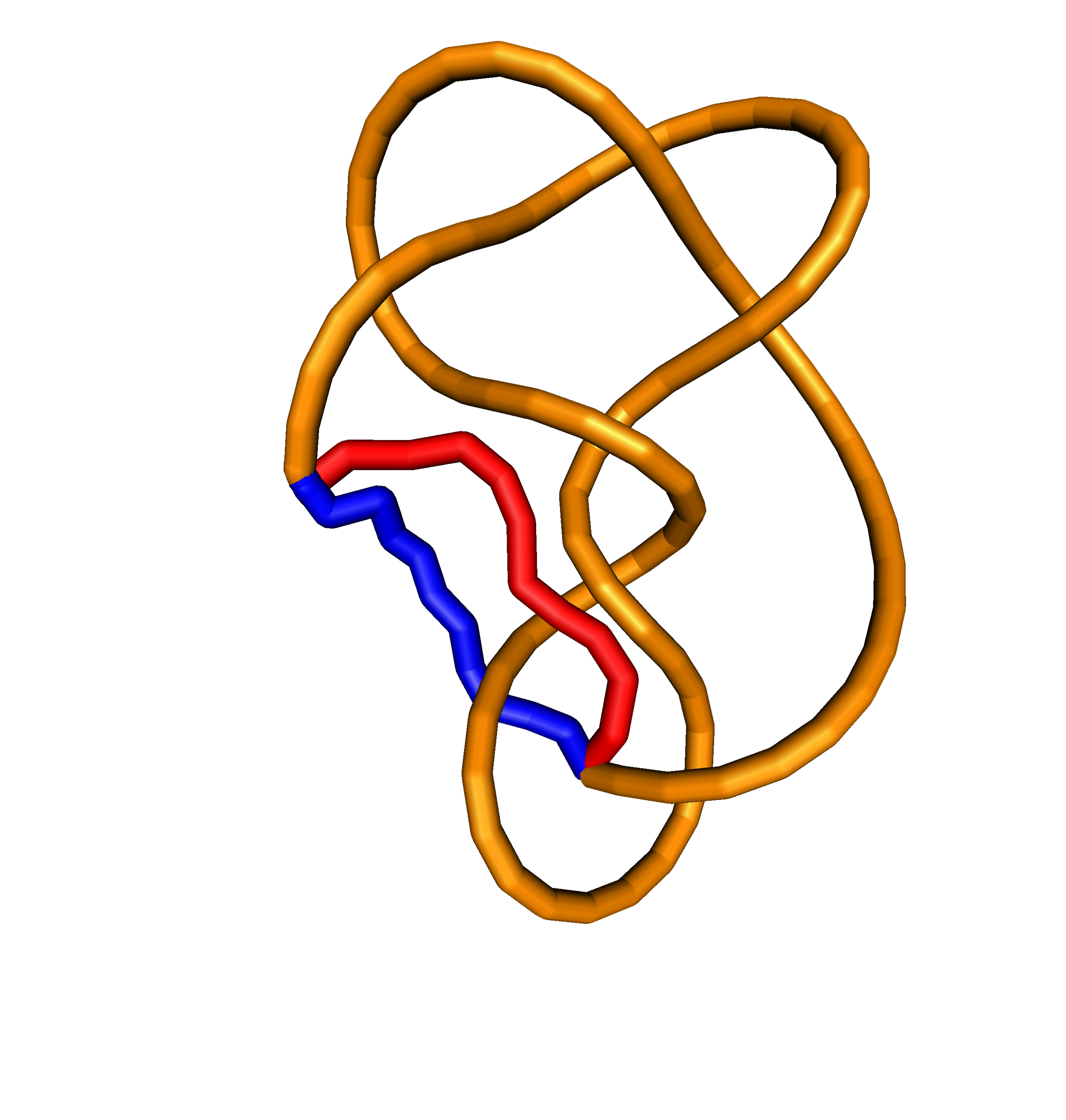}
    \end{overpic}
  \hfill{\ }
  \caption{Examples of closures used by $\operatorname{PU}$ (left) and $\operatorname{PR}$
    (right) associated with a subchain of the 98-edge KnotPlot $+7_5$
    knot \cite{knotplot} with 11 edges removed.  Note that the
    configurations analyzed in this paper are random, and thus are not as
    elegant as the configuration shown here.}
  \label{stochasticclosure}
\end{figure}

The Uniform Closure Method is defined as follows.  For a given open
chain and a given direction (seen as point on the unit sphere $S^2$),
one can create a closed knot by extending rays in that direction from
the two endpoints and closing the knot at infinity.  Then the Uniform
Closure Method ``knot type'' of the open chain is the probability
distribution of closed knot types obtained over all points on $S^2$.
A close relative of this method was initially introduced in
\cite{mdsmacro}, with the small difference that the authors closed the
knot by extending line segments from the endpoints to random points on
a large sphere encompassing the chain.  Note that in practice, one
needs only to have the rays extend beyond the convex hull of the
chain, and then the configuration can be closed via a straight line
segment, as seen on the left of Figure~\ref{stochasticclosure}.

Our main experimental finding is that the Uniform Closure Method is
almost exactly as good at classifying subarcs as the Bayes MAP
classifier.

\begin{defn}
Given a $k$-edge arc $A$ and any $n > k$, the \emph{Bayes classifier with maximum a posteriori probability (MAP) decision rule} for the knot type of $A$ is
\begin{equation*}
\operatorname{BC}(A,n) := 
\text{ 
\begin{minipage}{4in}
the most common knot type among closed \\ equilateral $n$-gons with subarc $A$.
\end{minipage}
}
\end{equation*}
If several knot types are equally common, the classifier selects among them at random. If there are no closed $n$-gons with subarc $A$ the classifier is undefined.
\end{defn}

It is a standard theorem that the Bayes MAP classifier minimizes the probability of misclassification~\cite{Devroye:1996wq}, so $\operatorname{BC}$ is the maximum accuracy classifier for the knot type of~$P$. If we can sample from closed equilateral $n$-gons with subarc $A$ by choosing random $(n-k)$-edge arcs to close $A$ (see Figure~\ref{stochasticclosure}B), we can approximate $\operatorname{BC}$ as follows:
\begin{defn}
Given a $k$-edge arc $A$ and any $n > k$, the \emph{random closure classifier} of the knot type of $A$ is
\begin{equation*}
\operatorname{PR}(A,n) := \text{
\begin{minipage}{4in}
the most common knot type among 
100 random closed \\ equilateral $n$-gons with subarc $A$.
\end{minipage}
}
\end{equation*}
\end{defn}
As one would expect, we could build improved versions of the classifier $\operatorname{PR}$ by increasing the number of samples or (potentially) by choosing closures from a low discrepancy sequence in the space of possible closures rather than at random. However, in practice, computer time is limited and the metric structure of the space of closures is not well-understood. So we investigated the performance of $\operatorname{PR}$ as a reasonable proxy for the performance of the maximum accuracy classifier $\operatorname{BC}$.

The posterior distribution on knot types defined by the Uniform Closure Method has its own MAP classifier which is, like $\operatorname{BC}$, impractical to compute. It is typically approximated as follows:

\begin{defn}
Given a $k$-edge arc $A$, the \emph{uniform closure classifier} of the knot type of $A$ is 
\begin{equation*}
\operatorname{PU}(A) := \text{
\begin{minipage}{4in}
the most common knot type among 
100 closures of $A$ by \\ parallel rays from the endpoints of $A$ (see Figure~\ref{stochasticclosure}).
\end{minipage}
}
\end{equation*}
The directions of the rays are a fixed set chosen by the applet~\cite{geodesic100} to cover the sphere evenly.\footnote{We do this because this matches the implementation of approximations to the Uniform Closure Method used elsewhere in the literature. Choosing evenly-spaced points provides a better approximation to the true posterior distribution than choosing random points.} We weight the knot type associated with each direction by the relative size of the spherical Voronoi cell of the corresponding ray, and define the ``most common'' knot type to be the one with the largest total weight.
\end{defn}

$\operatorname{PU}$ is clearly important for several reasons. The entanglement measured by $\operatorname{PU}$ has shown itself to be biologically relevant (again, see \cite{openknotting,knotencyclopedia,fingerprints,knotprot,subknots}). $\operatorname{PU}$ is considerably faster to compute than $\operatorname{PR}$ because it deals with $k+2$ edge polygons instead of $n$ edge polygons. And perhaps most usefully, given an open arc $A$, $\operatorname{PU}$ does not require one to fix a length $n-k$ for the closing arc. 

The posterior distribution on knot types given by the Uniform Closure Method appears very different from the ``true'' posterior distribution of random arc closures. 
Because the Bayes MAP classifier has the maximum possible accuracy, $\operatorname{PU}$ cannot be more accurate than $\operatorname{BC}$, and it is very unlikely to be more accurate than $\operatorname{PR}$, which approximates $\operatorname{BC}$. Indeed, we might expect $\operatorname{PU}$ to be \emph{much} less accurate than $\operatorname{PR}$ because the posterior distributions are so different.

We now argue informally that, under plausible assumptions, the distribution of knot types in the posterior distributions used by $\operatorname{PR}$ and $\operatorname{PU}$ may not be so different. If one believes that knots in random polygons tend to be localized, then the $(n-k)$-edge section of a random polygon $P$ joining the ends of a $k$-edge subarc $A$ containing a local knot is likely to stay well away from the small knot in $A$. In this case, if the $(n-k)$-edge arc is short enough that it does not introduce additional knotting, the knot type of $P$ is likely to be the same as that of a knot formed by joining the ends of $A$ ``to infinity'' with parallel rays. If this argument is right, at least for some $n$, $\operatorname{PU}$ and $\operatorname{PR}$ may perform similarly.

We can go a little further.\footnote{Again, this is an intuitive argument and is not meant as a proof!}
In probability problems with a large number of degrees of freedom, concentration of measure often leads the distribution of some relatively coarse labeling function on the space to be highly concentrated on its expected value or mode. Therefore, we might guess that almost all closures of $A$ (whether by random arcs or rays) have the same knot type. If so, then the following classifiers might be almost as accurate as $\operatorname{PR}$ and $\operatorname{PU}$ at much lower computational cost.

\begin{defn}
Given a $k$-edge arc $A$ and any $n > k$, the single \emph{random closure classifier} of the knot type of $A$ is
\begin{equation*}
\operatorname{SR}(A,n) := \text{
\begin{minipage}{4in}
the knot type of one 
random closed equilateral $n$-gon \\ with subarc $A$.
\end{minipage}
}
\end{equation*}
\end{defn}

\begin{defn}
Given a $k$-edge arc $A$, the \emph{single uniform closure classifier} of the knot type of $A$ is 
\begin{equation*}
\operatorname{SU}(A) := \text{
\begin{minipage}{4in}
the knot type of the 
closure of $A$ by parallel rays from the endpoints of $A$ whose direction is chosen randomly from the 100 fixed directions used in $\operatorname{PU}$.
\end{minipage}
}
\end{equation*}
\end{defn}

We now have two hypotheses about accuracy to test: when $n$ is not too much larger than~$k$, $\operatorname{PR}$ and $\operatorname{PU}$ should be approximately equally accurate (and both close to the maximum possible accuracy), and for reasonably large $n$ and $k$, even $\operatorname{SR}$ and $\operatorname{SU}$ should be close to the maximum possible accuracy. We will see data below which supports both of these hypotheses. We will also see support for a much stronger result: the accuracies of all four classifiers are very close even when we restrict our attention to the tiny subspace of quite rare and complicated random knots. 

However, when studying the performance of classifiers, accuracy is not the only important quantity. So we also compare the positive predictive value\footnote{The chance that a prediction of knot type $K$ is accurate.} (PPV) of a variety of predictions made by each of our four classifiers from unknots to 9-crossing knots. Here we see that the (easier to compute) $\operatorname{PU}$ is just as good as $\operatorname{PR}$, but that there are large differences between the PPV of these two and the PPV of $\operatorname{SU}$ and $\operatorname{SR}$. Interestingly, though $\operatorname{SR}$ may seem ``more natural'' than $\operatorname{SU}$,\footnote{At least to some of the authors of this paper!} the observed PPV for $\operatorname{SU}$ is substantially higher for every predicted knot type than the PPV of $\operatorname{SR}$.

Our experiments are computationally expensive, so we could not choose a wide range of $n$ and $k$.    Instead, we focused on random closed 100-edge equilateral polygons, analyzing the performance of our four classifiers at recovering the knot type of a particular polygon given the partial information in a subarc of length $k$. We computed accuracy and PPV for all $k$ from $1$ to $100$, finding as expected that all of our classifiers performed better when given more information.

\section{Experimental procedure}

\subsection{Generation of random closures for $\operatorname{PR}$ and $\operatorname{SR}$}

In the random closure method, we must sample $100$ random $n$-gons containing the initial $k$-edge arc $A$. It is a classical observation~\cite{Rayleigh:1919do} that the uniform probability distribution\footnote{That is, the Hausdorff measure on the subspace of $(S^2)^n$ consisting of edge directions which sum to $\vec{0}$.} on equilateral $n$-edge polygons can be written as the probability distribution of $(n-k)$-edge open arcs and $k$-edge open arcs,\footnote{The distribution of a $m$-edge arc is the uniform distribution of edge directions in $(S^2)^{m}$.} conditioned on the hypothesis that the end-to-end distances of the two arcs agree. 

Since the end-to-end distance $\ell$ of $A$ is fixed, we can therefore construct a sample from the conditional distribution of closed equilateral $n$-gons containing $A$ by constructing a random $(n-k)$-edge arc $B$ with the end-to-end distance $\ell$ and joining its ends to $A$. Constructing such a $B$ is equivalent to sampling a random \emph{closed} $(n-k+1)$-gon with  $n-k$ edges of length $1$ and one edge of length $\ell$. To do so, we used {\tt plcurve}~\cite{plcurve}, which implements an algorithm of Cantarella and Shonkwiler~\cite{Cantarella:2016iy}. We note that $A$ and $B$ can be rotated independently around the vector joining their endpoints, and that we must choose the rotation angle uniformly to get the correct posterior distribution of $n$-gons. We see an example on the right of Figure~\ref{stochasticclosure}.  

\subsection{Determination of knot types}

Throughout this paper, the knot types of the polygons
were determined using a combination of Ewing and Millett's HOMFLYPT
\cite{HOMFLY,PT}
polynomial software \cite{millettewing} and Hoste and Thistlethwaite's
\texttt{knotfind} program (which is a part of the 
\texttt{Knotscape} software \cite{knotscape}).  Thistlethwaite's
\texttt{unraveller} program \cite{unraveller} was used to simplify
projections, easing the burden in computing the knot types.  Using
these programs, we were able to identify knot types of closed polygons
as long as the knot type has a crossing number of 16 or less
(implicitly assuming that the crossing number of a composite knot type
is the sum of the crossing numbers of the factor knots, which is
conjectured but has not been proven).  

For all but a handful of knot
types (those for which the HOMFLYPT polynomials match for chiral
pairs), we can distinguish between the chiralities of the knot types
as well.  For chiral knot types, we designated one of the
pair as the $+$ version and one as the $-$ version, based on the
writhe of the standard diagram for the knot, or the spatial writhe of
the ropelength-minimization \cite{expmath} if the writhe of the
standard diagram for the knot is zero.  

\subsection{Sampling procedure for test data}

We first generated 5,000,000 random 100-edge equilateral
polygons using the \texttt{plCurve} library from Cantarella's lab~\cite{plcurve}. We computed the knot type of each polygon as above, obtaining the distribution of knot types in Table~\ref{knotprobtable}.

\begin{table}
  \centering
  \hfill{\ }
  \begin{tabular}[t]{cc} \toprule
    Knot & \% \\
    \midrule
    $0_1$        & 70.6613 \\
    $+3_1$       &  8.4252 \\
    $4_1$        &  3.2827 \\
    $+5_1$       &  0.4994 \\
    $+5_2$       &  0.8546 \\
    $+6_1$       &  0.2066 \\
    $+6_2$       &  0.2358 \\
    $6_3$        &  0.2694 \\
    $+3_1\#+3_1$ &  0.3865 \\
	\bottomrule
  \end{tabular}
  \hfill
  \begin{tabular}[t]{cc} \toprule
    Knot & \% \\
    \midrule
    $+3_1\#-3_1$ & 0.7738 \\ 
    $+7_1$       & 0.0234 \\ 
    $+7_2$       & 0.0499 \\ 
    $+8_2$       & 0.0113 \\ 
    $+8_8$       & 0.0201 \\ 
    $+8_{16}$    & 0.0039 \\ 
    $+8_{19}$    & 0.0263 \\ 
    $+8_{20}$    & 0.0518 \\ 
    $+8_{21}$    & 0.0293 \\ 
	\bottomrule
  \end{tabular}
  \hfill
  \begin{tabular}[t]{cc} \toprule
    Knot & \% \\
    \midrule
    $+9_2$    & 0.0024 \\ 
    $+9_{15}$ & 0.0035 \\ 
    $+9_{21}$ & 0.0028 \\ 
    $+9_{36}$ & 0.0023 \\ 
    $+9_{43}$ & 0.0107 \\ 
    $+9_{44}$ & 0.0162 \\ 
    $+9_{46}$ & 0.0054 \\ 
	\bottomrule
  \end{tabular}
  \hfill{\ }
  \caption{Percentage of the 5,000,000 randomly generated equilateral
    100-gons which form the knot types in our sample.}
  \label{knotprobtable}
\end{table}

Since the frequency of different knot types varies by several orders of magnitude, it is interesting to compare the performance of the classifiers both on the overall data set and on particular knot types (even very rare ones). The reason for this is that a classifier which was particularly good at recognizing unknots might have a very good overall accuracy score (because unknots are very common), but perform very poorly at recognizing the complicated knots that one is presumably most interested in.

So we took 100 samples from each of the knot types $0_1$, $+3_1$, $4_1$, $+5_1$, $+5_2$, $+6_1$, $+6_2$,
$6_3$, $+3_1\#+3_1$, $+3_1\#-3_1$, $+7_1$, $+7_2$, $+8_2$, $+8_8$,
$+8_{16}$, $+8_{19}$, $+8_{20}$, $+8_{21}$, $+9_2$, $+9_{15}$,
$+9_{21}$, $+9_{36}$, $+9_{43}$, $+9_{44}$, and $+9_{46}$.  
These knot types were chosen to provide a variety of different
crossing numbers, and a combination of alternating and
non-alternating, composite and prime, and chiral and amphichiral knot
types.  Note that we need to include chirality in the computations to
avoid false-positives.  For the chiral knot types, we always analyzed
the configurations with the positive version of the knot type.  The
results for the other chirality were similar (up to computational
error).  The amphichiral prime knots analyzed are $0_1$, $4_1$, and $6_3$.
The knots $+3_1\#+3_1$ and $+3_1\#-3_1$ are the composite granny and square
knots, respectively.  The knots $+8_{19}$, $+8_{20}$, $+8_{21}$,
$+9_{43}$, $+9_{44}$, and $+9_{46}$ are non-alternating.

\subsection{Experiments measuring accuracy of classifiers}

\begin{figure}
  \centering
  {\ }\hfill
  \begin{overpic}[angle=0,width=0.46\textwidth]{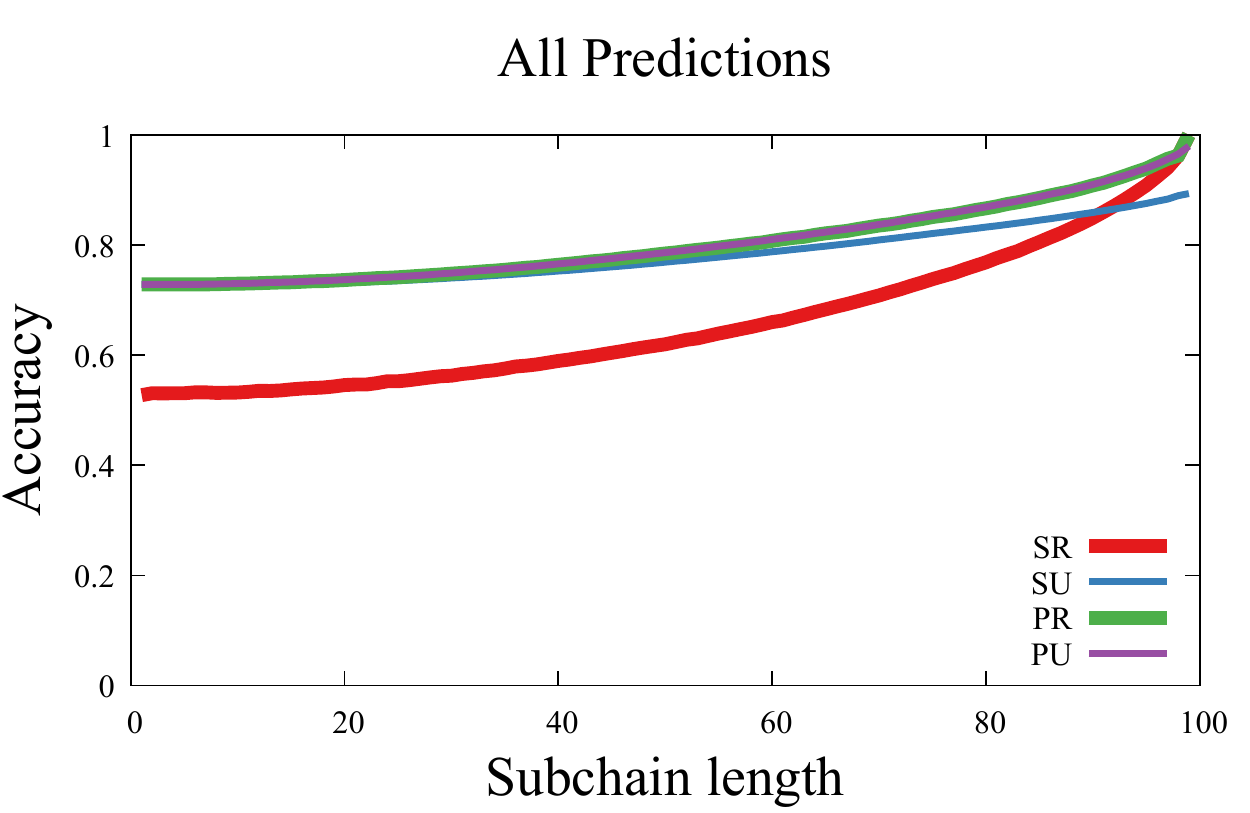}
    \end{overpic}
  \hfill
  \begin{overpic}[angle=0,width=0.46\textwidth]{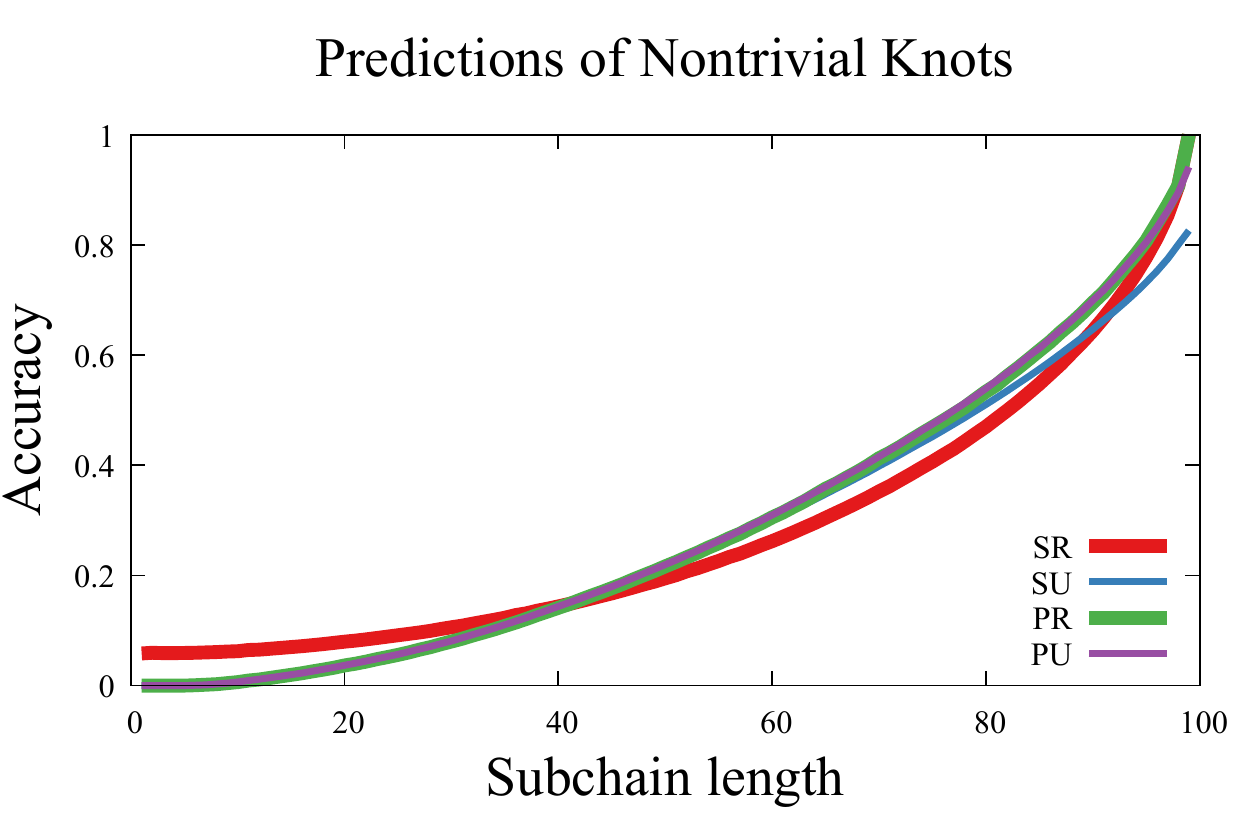}
    \end{overpic}
  \hfill{\ }
  \caption{The accuracy of each classifier is shown as a function of subchain length $k$ for a random 100-gon of any knot type including unknots (left plot), and a random 100-gon known to be nontrivially knotted (right plot).}
  \label{fig:combinedaccuracy}
\end{figure}

For each 100-edge polygon, there were 9900 
different open, connected subchains of various lengths. We ran all four classifiers on each subchain,\footnote{Because we analyze all subchains of a given polygon, there is
some worry about self-correlations within our data set.  We did an
experiment where we analyzed only a single subchain per closed polygon
and the results were similar to what we report below, albeit with
more noise in the graphs.} recording their predictions as accurate when they matched the knot type of the overall closed polygon. We also computed the PPV. We recorded the accuracy for each classifier by subchain length and knot type. We recovered an overall accuracy score for all samples and all knots by reweighting each polygon by the probability of its knot type. Our overall accuracy scores are shown in Figure~\ref{fig:combinedaccuracy}. Note that each of the data
points is based on 1,000,000 closures (100 samples per knot
type, 100 subchains of the given length per sample, 100 closures per
subchain).

We can check that the data makes sense by considering the small $k$ limit. Every ray closure of a small arc is unknotted, so $\operatorname{PU}$ and $\operatorname{SU}$ always classify short arcs as unknots. Since about $70\%$ of the random 100-gons are unknots, these classifiers are accurate about $70\%$ of the time, which is the (best possible) Bayes accuracy. The $\operatorname{PR}$ classifier constructs 100 completely random equilateral 100-gons and returns the most common knot type observed. This most common knot type is $0_1$ virtually all of the time,
so $\operatorname{PR}$ is also very close to Bayes accuracy for small $k$. On the other hand, $\operatorname{SR}$ is only accurate for small $k$ when two independent random 100-gons happen to have the same knot type. The probability that both samples agree \emph{and are unknots} is roughly $0.7^2 = 0.49$. Other possible coincidences of knot type are unlikely, but possible; the total probability of agreement is about $52\%$, consistent with the results we see for the red $\operatorname{SR}$ curve on the left hand side of the left plot in Figure~\ref{fig:combinedaccuracy}.

A striking feature of Figure~\ref{fig:combinedaccuracy} is the essentially equal accuracy of $\operatorname{PU}$ and $\operatorname{PR}$ for all $k$, whether over the entire sample (left plot) or just the knotted curves (right plot). As $k$ approaches $n$, we see the reassuring result that all classifiers approach $100\%$ accuracy; that is, they essentially\footnote{We note that when $k=n-1$, $\operatorname{SR}$ and $\operatorname{PR}$ are always correct. In theory, $\operatorname{SU}$ or (less likely) $\operatorname{PU}$ \emph{could} be wrong if the ray closure directions pass back through the body of the knot in just the wrong way. But this does not seem to occur often enough to be visible in our experimental data.} always yield the correct knot type when given the entire knot as input!

It might be somewhat disconcerting that the red $\operatorname{SR}$ classifier outperforms our approximation $\operatorname{PR}$ of the Bayes MAP classifier for small $k$ in the right hand plot. Is it not the case that the MAP classifier always has the maximum possible accuracy? The solution to this apparent paradox is that the inputs on the right-hand plot are chosen only from knotted polygons, so $\operatorname{PR}$ is approximating the wrong Bayes MAP classifier. The correct Bayes MAP classifier would predict the most likely \emph{nontrivial} knot ($\pm 3_1$) and be accurate $\sim 8.42\%$ of the time for small $k$; better than both $\operatorname{SR}$ and $\operatorname{PR}$.

\begin{figure}[ht]
  \centering
   {\ }\hfill
  \includegraphics[width=0.45\textwidth]{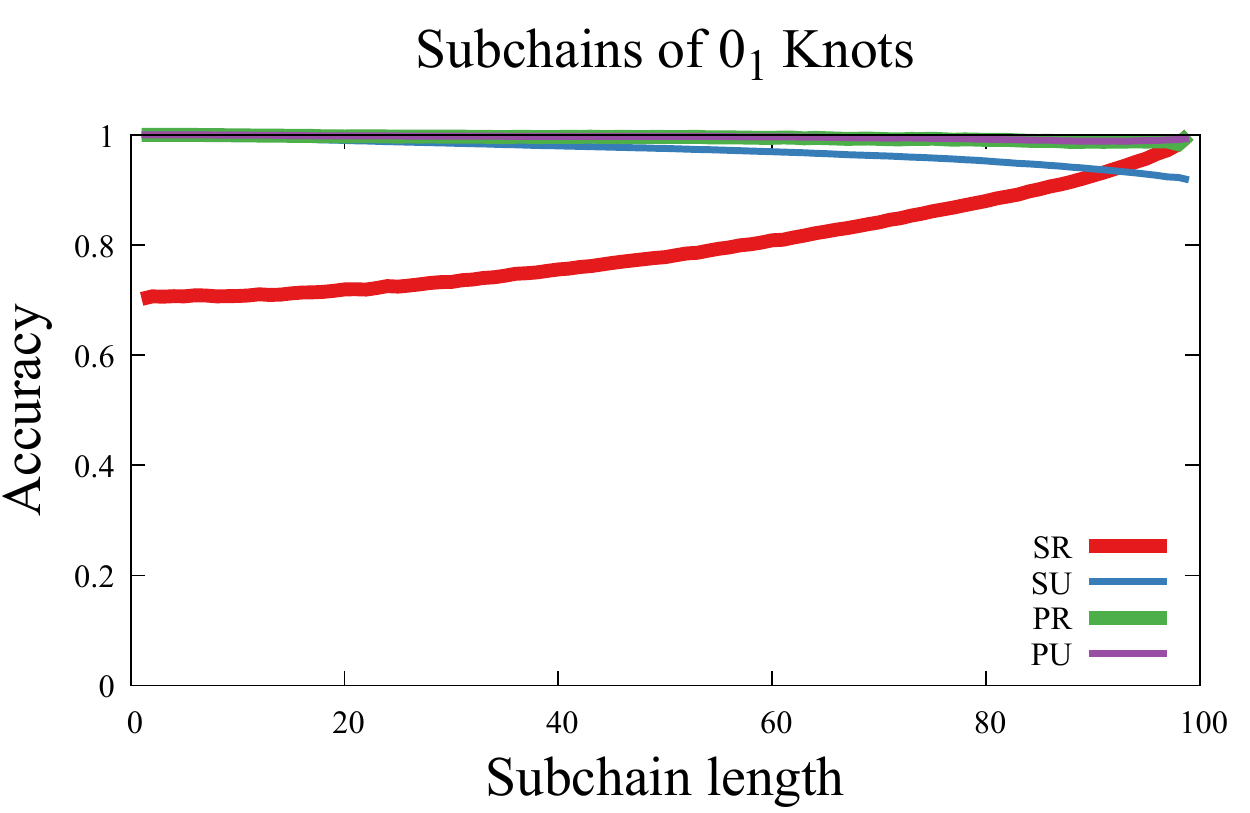}
  \hfill
  \includegraphics[width=0.45\textwidth]{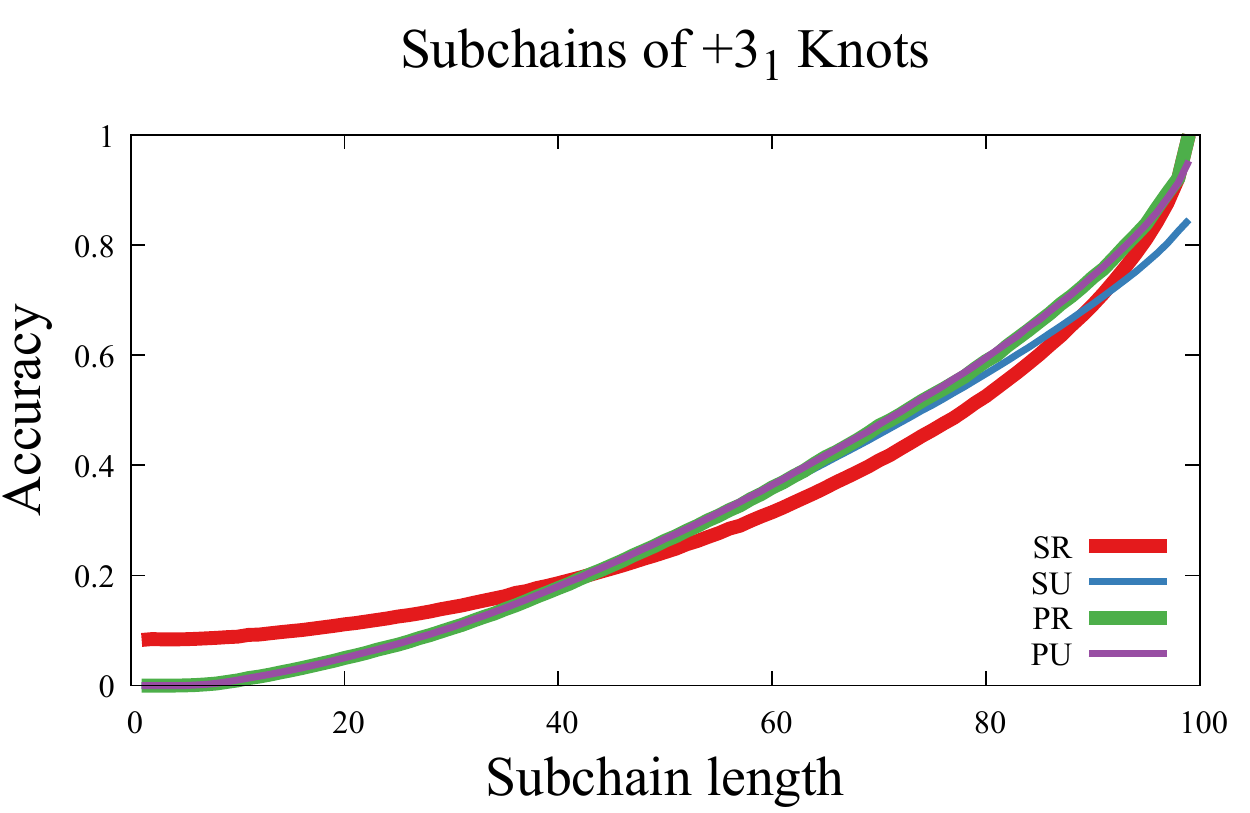}
  \hfill{\ }
  
  \centering
  {\ }\hfill
  \includegraphics[width=0.45\textwidth]{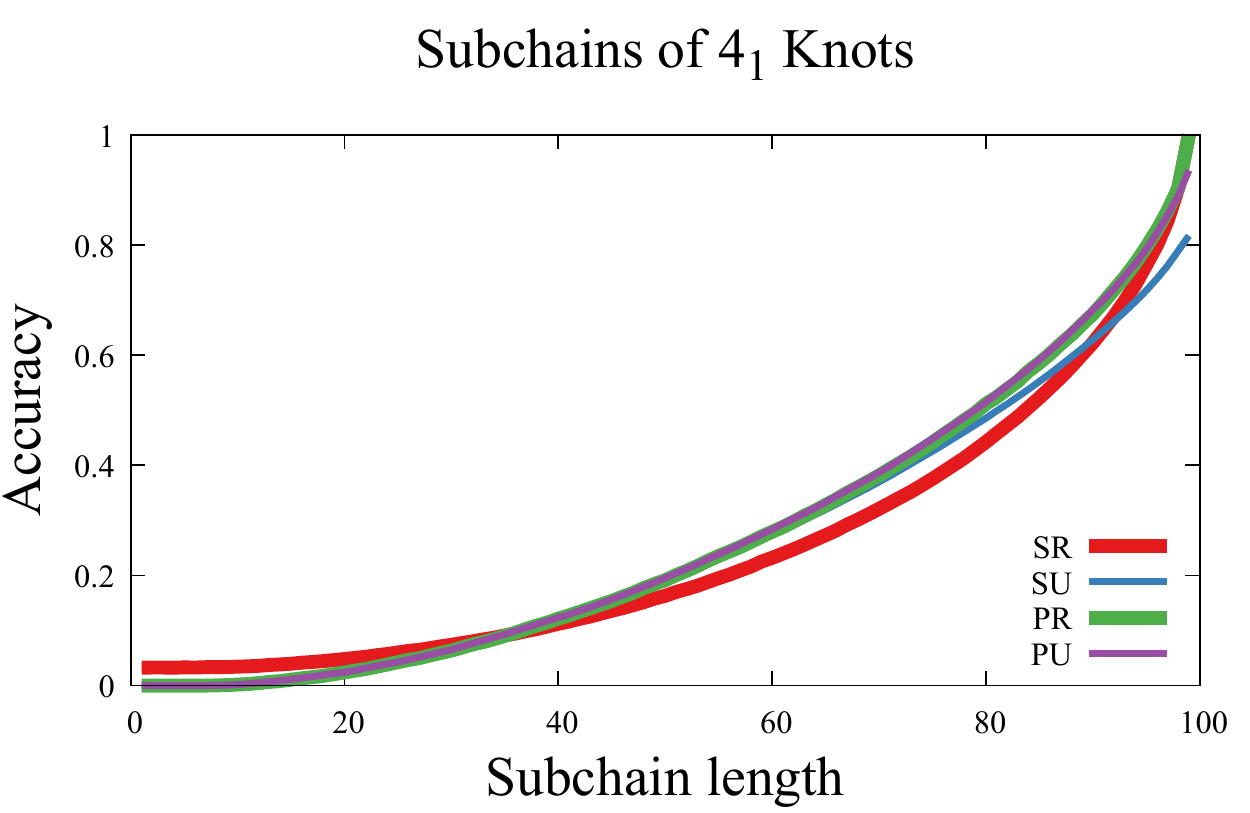}
  \hfill
  \includegraphics[width=0.45\textwidth]{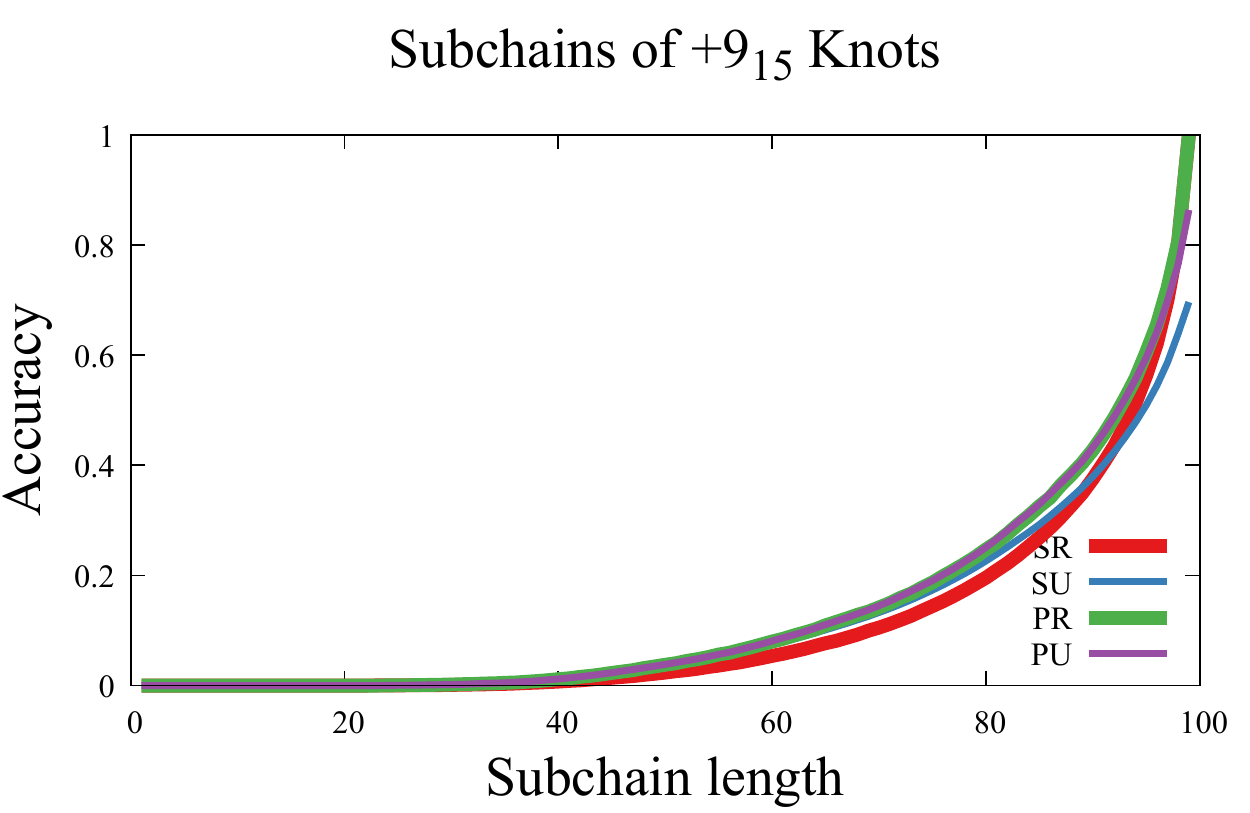}
  \hfill{\ }
  \caption{With the exception of the unknot, our four classifiers have comparable
  accuracy across all $k$ when the input set was restricted to curves of a given knot type. We observed this effect in all the knot types we tested.}
  \label{accuracycomparable}
\end{figure}

We found that the accuracy of all four classifiers remained comparable if we restricted our 
attention to curves of any single (nontrival) knot type. Figure~\ref{accuracycomparable} shows some example accuracy plots for the $+3_1$, $4_1$ and $+9_{15}$ knots. The only significant difference in performance was observed for the $\operatorname{SR}$ classifier on unknots, which started at the expected $\sim 70\%$ accuracy and gradually improved to $100 \%$ as $k$ increased.

We can see from Figure~\ref{accuracycomparable} that our classifiers require a larger fraction of the curve to recognize knot type when the knot type of the host is more complicated. Table \ref{passes50graph} shows the number of segments required to reach $50\%$ accuracy for each of our classifiers for the non-trivial knot types studied here. We note that the $\operatorname{PU}$ and $\operatorname{PR}$ classifiers are consistently very slightly faster at recognizing knots (though only by a few segments). One way to think about the number of segments required to recognize a knot is that it is a qualitative measure of something like an average size of the knot in a random configuration. This size generally increases with crossing number. However, the 8-crossing non-alternating knots ($+8_{19}$, $+8_{20}$, $+8_{21}$) are a few segments smaller than the 8-crossing alternating knots ($+8_{2}$, $+8_{8}$, $+8_{16}$), behaving more like 7-crossing alternating knots. The same is true for the 9-crossing non-alternating knots ($+9_{43}$, $+9_{44}$, $+9_{46}$) which compare to 8-crossing alternating knots. This phenomenon is shared with other ``geometric'' measures of knot complexity: generally non-alternating knot types have higher
probabilities in random knot studies \cite{modelsofrandomknots} and lower knot energy values \cite{expmath,torusenergy,spatial} than
alternating knot types with the same crossing number.

\begin{table}
   \centering
  {\ }\hfill
  \begin{tabular}{ccccc} \toprule
    Knot & SU & PU & SR & PR \\
    \midrule
    $+3_1$       & 75 & 73 & 79 & 73 \\
    $4_1$        & 82 & 79 & 84 & 80 \\
    $+5_1$       & 86 & 83 & 87 & 84 \\
    $+5_2$       & 87 & 84 & 88 & 85 \\
    $+6_1$       & 90 & 87 & 90 & 88 \\
    $+6_2$       & 90 & 87 & 90 & 87 \\
    $6_3$        & 90 & 87 & 90 & 88 \\
    $+3_1\#{+3_1}$ & 90 & 88 & 90 & 89 \\
    $+3_1\#{-3_1}$ & 90 & 88 & 90 & 88 \\
    $+7_1$       & 92 & 89 & 91 & 90 \\
    $+7_2$       & 92 & 89 & 92 & 89 \\
    $+8_2$       & 94 & 92 & 93 & 92 \\
    \bottomrule
  \end{tabular}
  \hfill
  \begin{tabular}{ccccc}
   \toprule
   Knot & SU & PU & SR & PR \\ 
    \midrule
     $+8_8$       & 94 & 92 & 93 & 92 \\
    $+8_{16}$    & 94 & 92 & 93 & 92 \\
     $+8_{19}$    & 91 & 88 & 91 & 88 \\
     $+8_{20}$    & 91 & 89 & 91 & 89 \\
     $+8_{21}$    & 92 & 89 & 92 & 89 \\
    $+9_2$       & 95 & 93 & 94 & 93 \\
     $+9_{15}$    & 95 & 93 & 94 & 93 \\
     $+9_{21}$    & 95 & 93 & 94 & 93 \\
     $+9_{36}$    & 95 & 93 & 94 & 93 \\
     $+9_{43}$    & 93 & 91 & 93 & 91 \\
     $+9_{44}$    & 93 & 91 & 93 & 91 \\
     $+9_{46}$    & 93 & 90 & 92 & 90 \\
     \bottomrule
  \end{tabular}
  \hfill{\ }
  \caption{Minimum number of edges at which $\operatorname{SU}$, $\operatorname{PU}$, $\operatorname{SR}$, and $\operatorname{PR}$ reach $50\%$ classification accuracy.}
  \label{passes50graph}
\end{table}

\subsection{Experiments measuring positive predictive value}

In the last section, we analyzed the accuracy of each of our classifiers and found them to be comparable in almost all cases. However, accuracy is not the only important measure of a classifier. In particular, one might be interested in the positive predictive value (PPV), which, as discussed in the Introduction, measures the probability that a knot type detected in a subchain actually matches the knot type of the host polygon. Like accuracy values, higher PPV values represent better performance. Unlike the accuracy, which is always bounded by the Bayes accuracy, the PPV can be arbitrarily close to 1.

We start by considering the mean PPV for each method as a function of $k$. Unlike accuracy, this mean is not weighted by the frequency with which the classifier makes a given prediction: it is simply the total PPV of all outputs of the classifier, divided by the number of possible outputs. It is important to note that this number of possible outputs changes with $k$, at least for the $\operatorname{SU}$ and $\operatorname{PU}$ classifiers; for small numbers of segments, ray closure simply cannot generate any very complicated knots. $\operatorname{SR}$, by contrast, can make any prediction, even for $k=1$, as can $\operatorname{PR}$.\footnote{Although it is overwhelmingly probable that $\operatorname{PR}$ will predict the unknot in this case, as we discussed above.}

\begin{figure}
  \centering
  {\ }\hfill
  \includegraphics[width=0.45\textwidth]{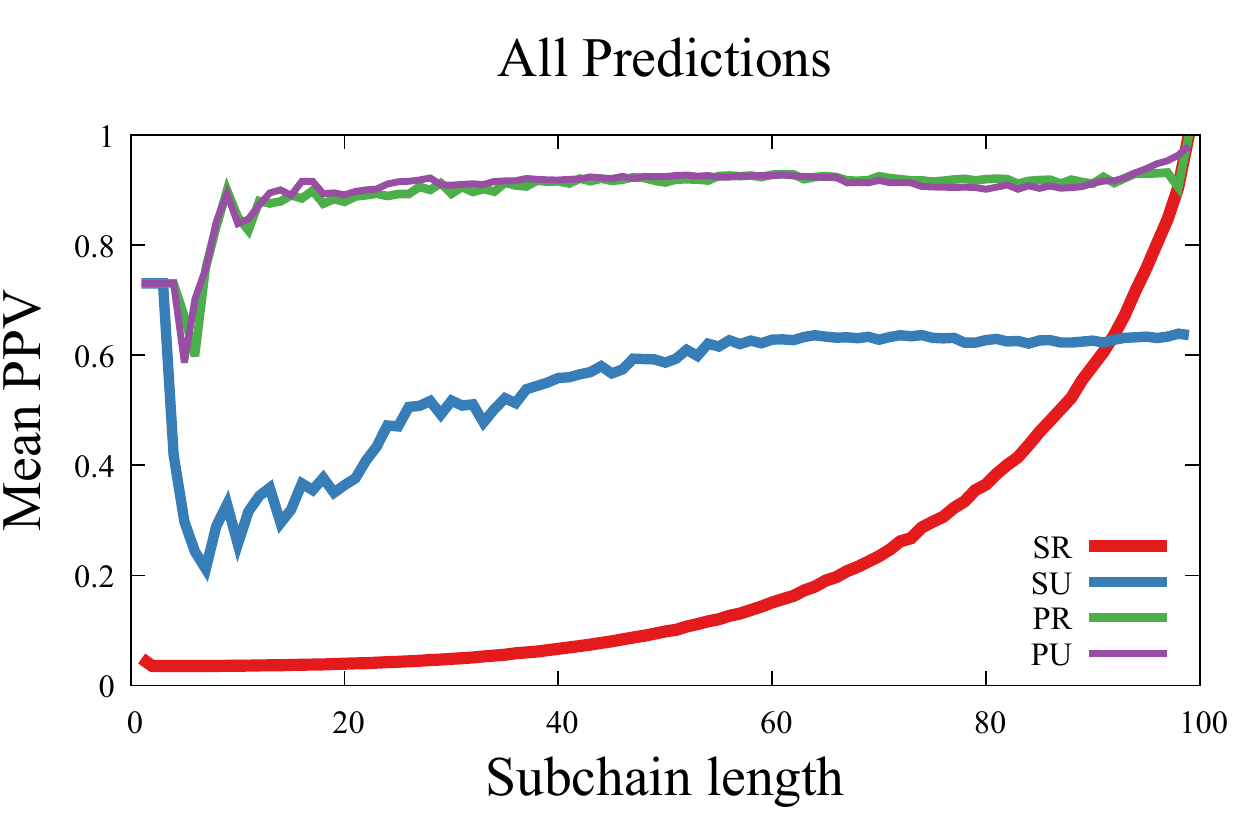}
  \hfill
  \includegraphics[width=0.45\textwidth]{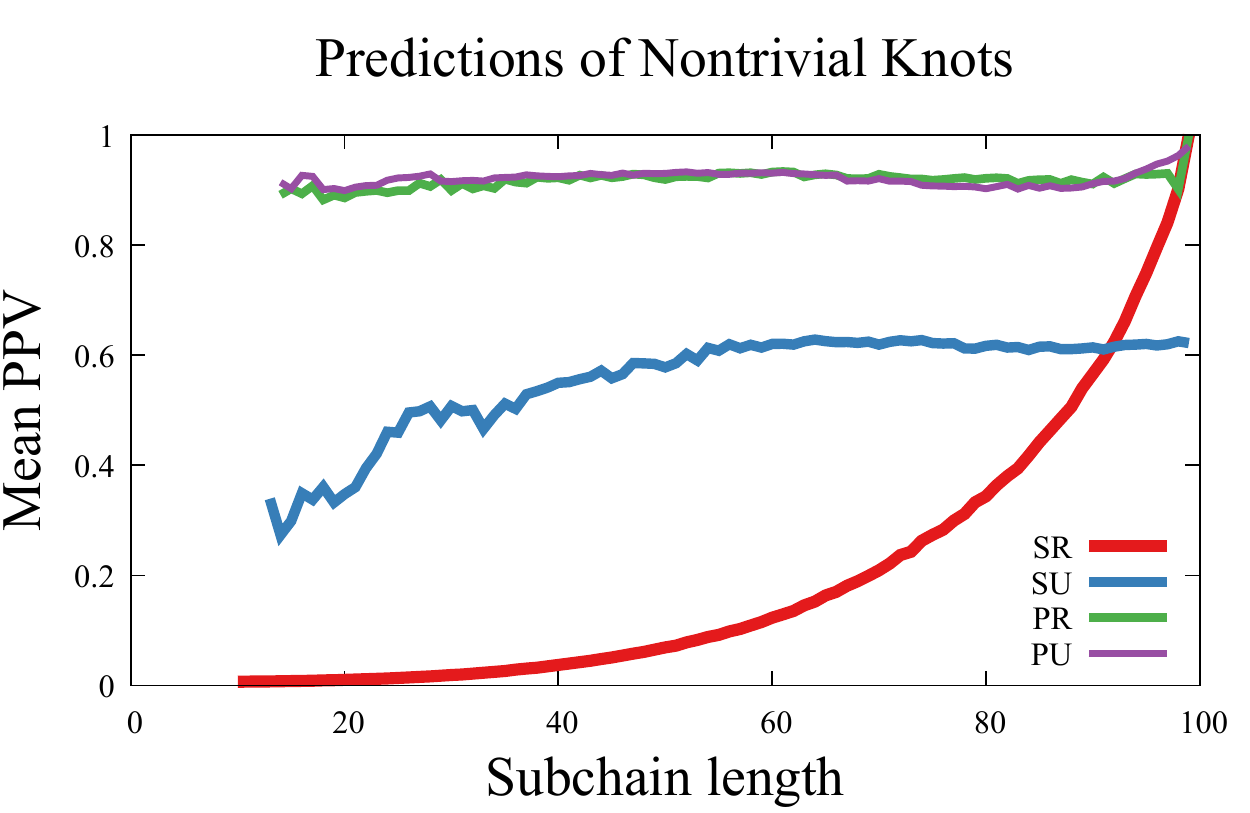}
  \hfill{\ }
  \caption{The mean PPV of the predictions produced by each classifier on our entire data set (100 knots of each type), for all predictions (left) and for predictions of nontrivial knots (right). Since the $\operatorname{SU}$ and $\operatorname{PU}$ classifiers cannot predict a nontrivial knot for very small $k$, we start the right hand plot at $k=10$. 
  }
  \label{meanppvgraphs}
\end{figure}
Figure~\ref{meanppvgraphs} shows the mean PPV of all predictions produced by each classifier for each value of $k$. This reveals two surprises. First, $\operatorname{PU}$ and $\operatorname{PR}$ are again essentially tied! Second, we see that $\operatorname{SU}$ is
significantly better than $\operatorname{SR}$ for almost all $k$.  

\begin{figure}
  \centering
   {\ }\hfill
  \includegraphics[width=0.45\textwidth]{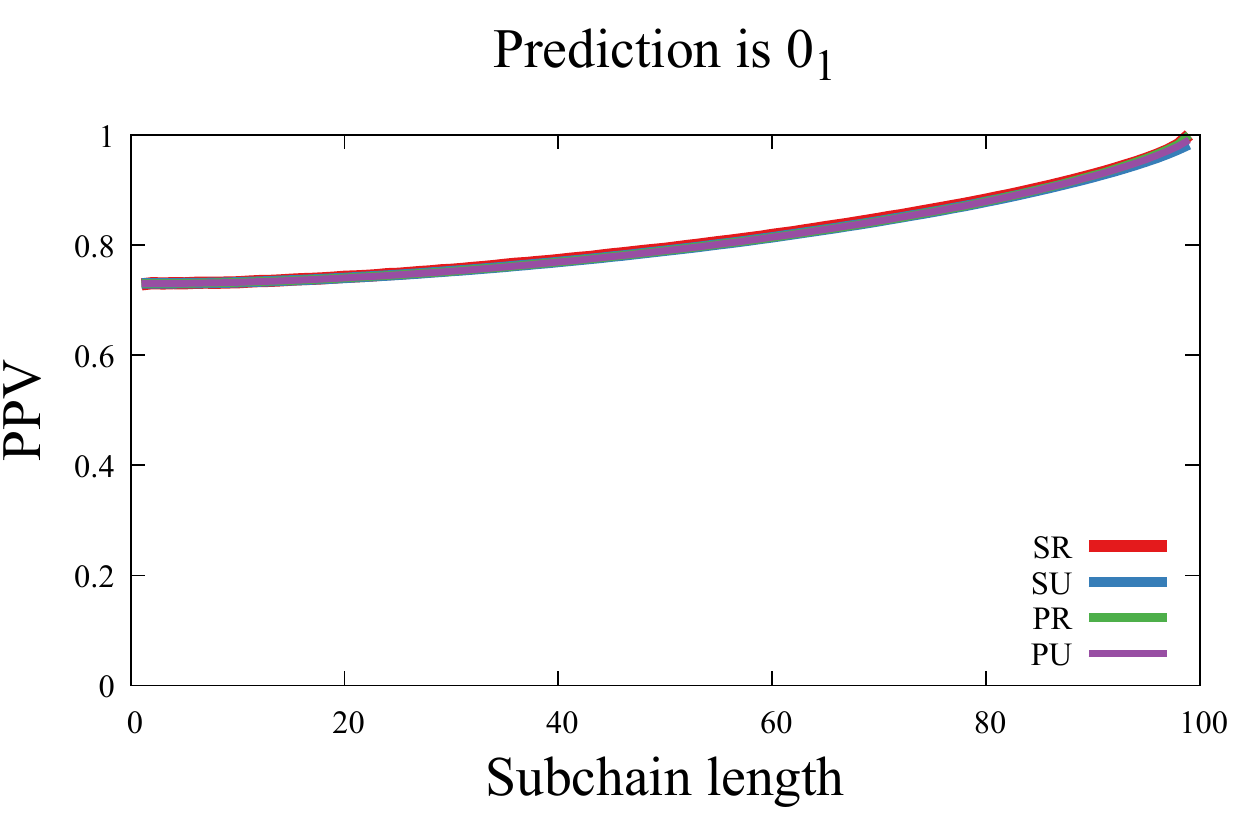}
  \hfill
  \includegraphics[width=0.45\textwidth]{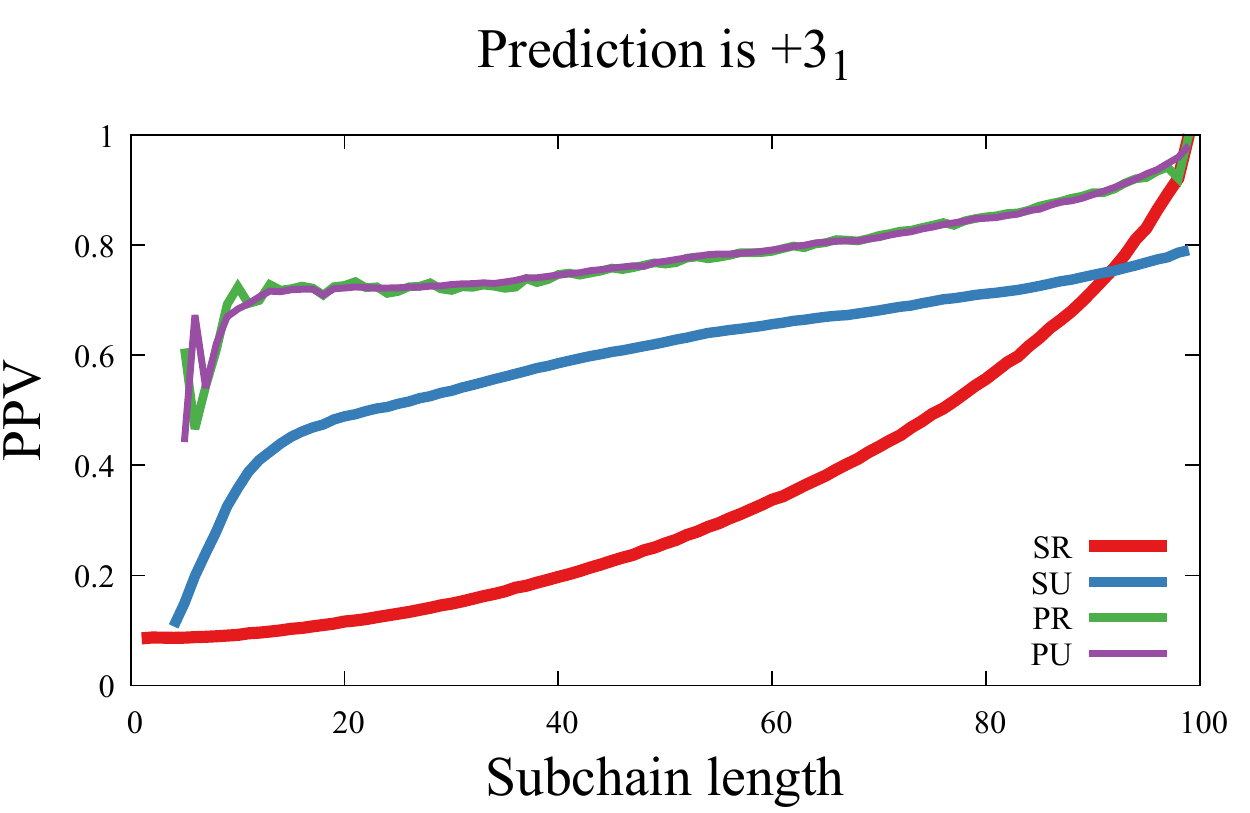}
  \hfill{\ }
  
  \centering
  {\ }\hfill
  \includegraphics[width=0.45\textwidth]{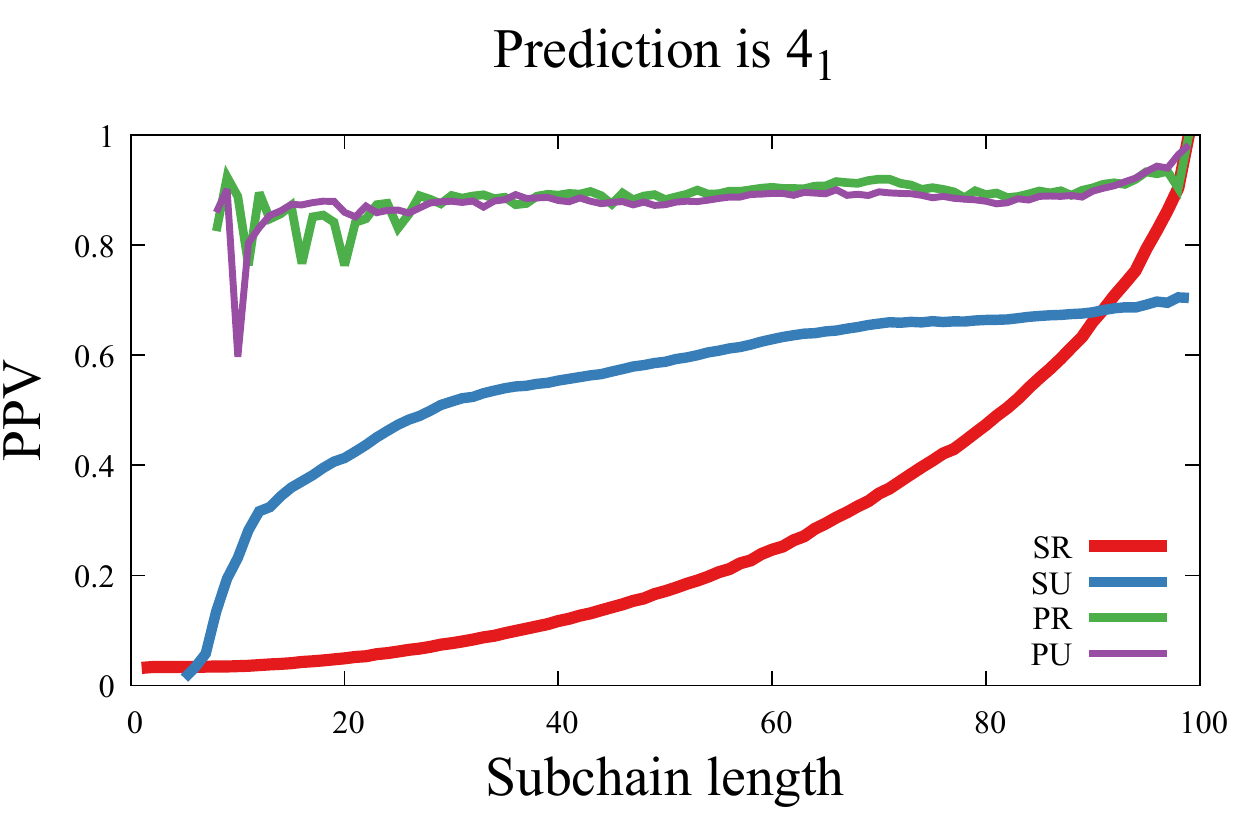}
  \hfill
  \includegraphics[width=0.45\textwidth]{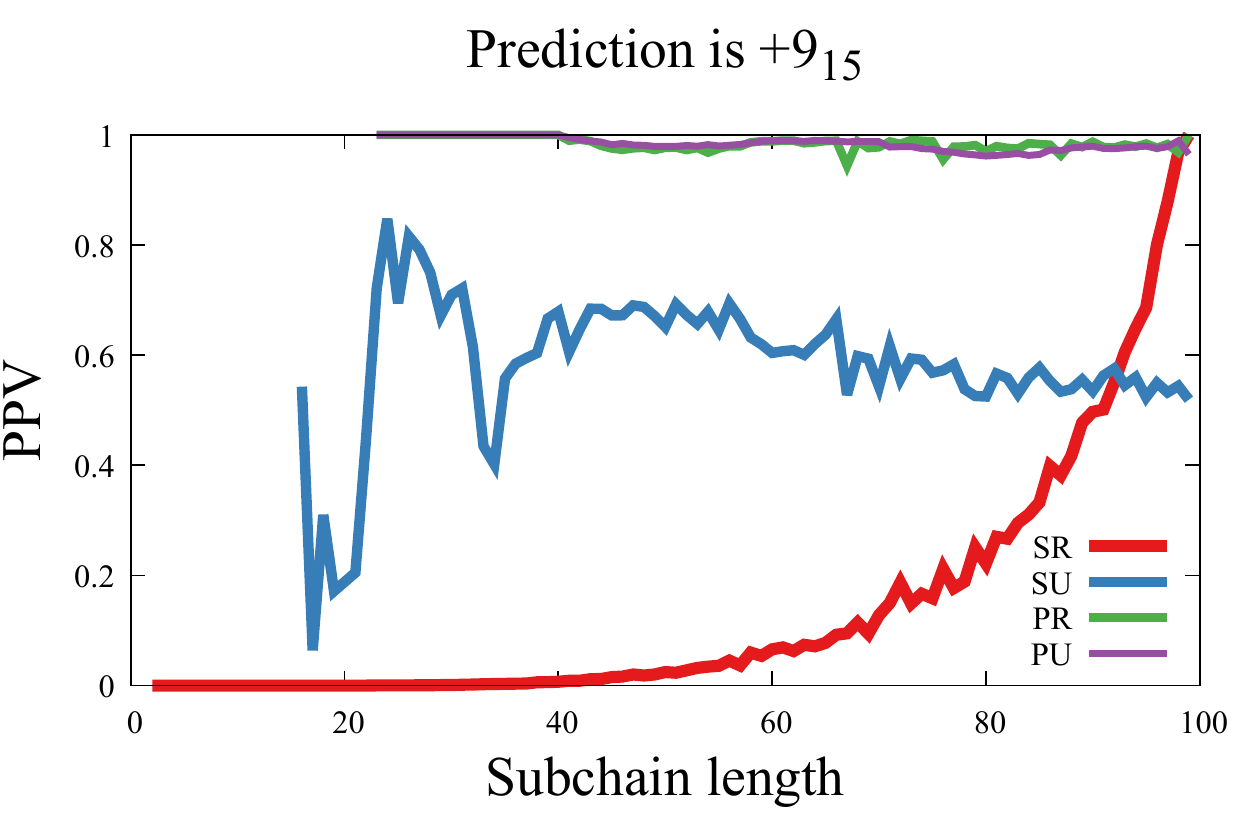}
  \hfill{\ }
  \caption{The PPV of the $\operatorname{PU}$ and $\operatorname{PR}$ classifiers was essentially the same, and increased with the complexity of the predicted knot type, eventually getting close to 1. The data is noisier for predictions of very complicated knots such as $+9_{15}$ because this prediction was made very rarely.}
  \label{ppvcomparable}
\end{figure}
We see a clear trend in the data that predictions of more complicated knotting are more likely to be correct. This is especially true for the MAP classifiers $\operatorname{PU}$ and $\operatorname{PR}$, which seem very unlikely to see a complicated knot a preponderance of the time unless the knot is really present in the host polygon.

In fact, the low PPV for predictions of trefoil knots tell us that many other host knots have arcs which are classified as trefoils, even by the high standards of the MAP classifiers $\operatorname{PU}$ and $\operatorname{PR}$. Some of this effect is doubtless due to the presence of composite knot types in our data such as ${+3_1}\#{+3_1}$. However, we suspect that other knots in the same family, such as $+5_1$, also contain subarcs classified as the trefoil knot~\cite{fertility,Chapman:2018ua,andrzejclassification,fingerprints,Millett:2016jz,Millett:2017iy,subknots}.

\section{Discussion}

In this paper, we found a surprising near-equivalence between the 
performance of the (approximate) Bayes MAP classifier $\operatorname{PR}$
and the uniform closure MAP classifier $\operatorname{PU}$ despite
what seemed to be large differences in the posterior distributions 
sampled by each method. The methods were comparable both in terms of 
accuracy and PPV, so experimenters should feel confident using them 
interchangeably, at least for $n \simeq 100$. 

In the process, we noticed that the $\operatorname{SU}$ 
classifier (which is 100$\times$ faster than the standard $\operatorname{PU}$ method)
has comparable accuracy. Therefore, $\operatorname{SU}$ may be an 
acceptable substitute for $\operatorname{PU}$ or $\operatorname{PR}$ when 
speed is more important than PPV.

Our experiments leave open the question of whether these phenomena persist
when $n$ is much larger. As discussed above, it seems very plausible that for 
very large $n$, $\operatorname{PR}$ would predict additional knotting that 
$\operatorname{PU}$ cannot. Whether this is desirable seems to be up to the 
experimenter. 

It would be interesting to compare $\operatorname{PU}$ to a different
measurement of ``surplus knotting''. Suppose $\mathcal{P}(n,K)$
is the probability of finding knot type $K$ in the prime decomposition of the
knot type of a random polygon of $n$ edges and $\mathcal{P}(n,K,A)$ was the 
corresponding probability for a random polygon of $n$ edges containing $A$ 
as a subarc. Intuitively, a measure of surplus knotting would compare $\mathcal{P}(n,K)$ to 
$\mathcal{P}(n,K,A)$ as $n \rightarrow \infty$. It seems reasonable 
to conjecture that this might be a better approximation of what $\operatorname{PU}$ 
is detecting.

\section{Acknowledgments}

Rawdon was supported by NSF DMS \#1115722, \#1418869, and \#1720342.
Annoni was supported by NSF DMS \#1115722.  Kumerow and Shogren were
supported by NSF DMS \#1115722 and \#1418869.  Brine-Doyle and Tibor
were supported by NSF DMS \#1418869. Cantarella and Shonkwiler were
supported by the Simons Foundation (\#524120 to JC, and \#354225 and
\#709150 to CS). Tibor, Annoni, Brine-Doyle, Kumerow, and Shogren
completed their work while students at the University of St.~Thomas.

\bibliographystyle{plain}
\bibliography{bibextra.bib,wonderbibv2_1.bib}

\begin{thebibliography}{10}

\bibitem{geodesic100}
Bob Allanson.
\newblock Martin's polyhedra.
\newblock http://members.ozemail.com.au/$\sim$llan/mpol.html, 2013.

\bibitem{plcurve}
Ted Ashton, Jason Cantarella, and Harrison Chapman.
\newblock plcurve.
\newblock http://www.jasoncantarella.com/wordpress/software/plcurve/, 2012.

\bibitem{expmath}
Ted Ashton, Jason Cantarella, Michael Piatek, and Eric~J. Rawdon.
\newblock Knot tightening by constrained gradient descent.
\newblock {\em Experiment. Math.}, 20(1):57--90, 2011.

\bibitem{stevedore}
Daniel B{\"o}linger, Joanna~I. Su{\l}kowska, Hsiao-Ping Hsu, Leonid~A. Mirny,
  Mehran Kardar, Jose~N. Onuchic, and Peter Virnau.
\newblock A stevedore's protein knot.
\newblock {\em PLoS Comput. Biol}, 6(4):e1000731, 2010.

\bibitem{openknotting}
Erin Brine-Doyle, Madeline Shogren, Emily Vecchia, and Eric~J. Rawdon.
\newblock Open knotting.
\newblock In Philipp Reiter, Simon Blatt, and Armin Schikorra, editors, {\em
  New Directions in Geometric and Applied Knot Theory}. De Gruyter, 2018.

\bibitem{jasoneqknots}
Jason Cantarella, Bertrand Duplantier, Clayton Shonkwiler, and Erica Uehara.
\newblock A fast direct sampling algorithm for equilateral closed polygons.
\newblock {\em J. Phys. A}, 49(27):275202, 2016.

\bibitem{fertility}
Jason Cantarella, Allison Henrich, Elsa Magness, Oliver O'Keefe, Kayla Perez,
  Eric~J. Rawdon, and Briana Zimmer.
\newblock Knot fertility and lineage.
\newblock {\em J. Knot Theory Ramifications}, 26(13):1750093, 2017.

\bibitem{Cantarella:2016iy}
Jason Cantarella and Clayton Shonkwiler.
\newblock {The symplectic geometry of closed equilateral random walks in
  $3$-space}.
\newblock {\em Ann. Appl. Probab.}, 26(1):549--596, 2016.

\bibitem{Chapman:2018ua}
Harrison Chapman.
\newblock On the structure and scarcity of alternating knots.
\newblock Preprint, {\tt arXiv:1804.09780 [math.GT]}, 2018.

\bibitem{knotprot2}
Pawel Dabrowski-Tumanski, Pawel Rubach, Dimos Goundaroulis, Julien Dorier,
  Piotr Su{\l}kowski, Kenneth~C. Millett, Eric~J. Rawdon, Andrzej Stasiak, and
  Joanna~I. Su{\l}kowska.
\newblock {KnotProt 2.0: a database of proteins with knots and other entangled
  structures}.
\newblock {\em Nucleic Acids Research}, 47(D1):D367--D375, 12 2018.

\bibitem{functionproteins}
Pawel Dabrowski-Tumanski, Andrzej Stasiak, and Joanna~I. Su\l{}kowska.
\newblock In search of functional advantages of knots in proteins.
\newblock {\em PLOS ONE}, 11(11):1--14, 11 2016.

\bibitem{joannapnaslinks}
Pawel Dabrowski-Tumanski and Joanna~I. Su\l{}kowska.
\newblock Topological knots and links in proteins.
\newblock {\em Proc. Natl. Acad. Sci. USA}, 114(13):3415--3420, 2017.

\bibitem{Devroye:1996wq}
Luc Devroye, L{\'a}szl{\'o} Gy{\"o}rfi, and G{\'a}bor Lugosi.
\newblock {\em {A Probabilistic Theory of Pattern Recognition}}, volume~31 of
  {\em Applications of Mathematics}.
\newblock Springer, New York, 1996.

\bibitem{knotencyclopedia}
Julien Dorier, Dimos Goundaroulis, Eric~J. Rawdon, and Andrzej Stasiak.
\newblock Open knots.
\newblock In Colin Adams, Erica Flapan, Allison Henrich, Louis~H. Kauffman,
  Lewis~D. Ludwig, and Sam Nelson, editors, {\em Encyclopedia of Knot Theory}.
  Chapman and Hall/CRC Press, 2020.
\newblock To appear.

\bibitem{modelsofrandomknots}
Chaim Even-Zohar.
\newblock Models of random knots.
\newblock {\em J. Appl. Comput. Topol.}, 1(2):263--296, 2017.

\bibitem{millettewing}
Bruce Ewing and Kenneth~C. Millett.
\newblock Computational algorithms and the complexity of link polynomials.
\newblock In Michel Boileau, Michel Domergue, Yves Mathieu, and Kenneth~C.
  Millett, editors, {\em Progress in Knot Theory and Related Topics}, volume~56
  of {\em Travaux en Cours}, pages 51--68. Hermann, Paris, 1997.

\bibitem{andrzejclassification}
Alessandro Flammini and Andrzej Stasiak.
\newblock Natural classification of knots.
\newblock {\em Proc. R. Soc. A}, 463(2078):569--582, 2007.

\bibitem{HOMFLY}
Peter~J. Freyd, David~N. Yetter, Jim Hoste, William B.~R. Lickorish, Kenneth~C.
  Millett, and Adrian Ocneanu.
\newblock A new polynomial invariant of knots and links.
\newblock {\em Bull. Amer. Math. Soc. (N.S.)}, 12(2):239--246, 1985.

\bibitem{andrzejchromatin}
Dimos Goundaroulis, Eriz Lieberman~Aiden, and Andrzej Stasiak.
\newblock Chromatin is frequently unknotted at the megabase scale.
\newblock {\em Biophys. J.}, 118(9):2268--2279, 2020.

\bibitem{Gukov:2020vr}
Sergei Gukov, James Halverson, and Fabian Ruehle.
\newblock Learning to unknot.
\newblock Preprint, {\tt arXiv:2010.16263 [math.GT]}, 2020.

\bibitem{knotscape}
Jim Hoste and Morwen Thistlethwaite.
\newblock Knotscape.
\newblock http://www.math.utk.edu/$\sim$morwen/knotscape.html, 2009.
\newblock Program for computing topological information about knots.

\bibitem{Hughes:2020hh}
Mark~C. Hughes.
\newblock {A neural network approach to predicting and computing knot
  invariants}.
\newblock {\em J. Knot Theory Ramifications}, 29(3):2050005, 2020.

\bibitem{fingerprints}
David A.~B. Hyde, Joshua Henrich, Eric~J. Rawdon, and Kenneth~C. Millett.
\newblock Knotting fingerprints resolve knot complexity and knotting pathways
  in ideal knots.
\newblock {\em J. Phys.: Condens. Matter}, 27:354112, 2015.

\bibitem{knotprot}
Michal Jamroz, Wanda Niemyska, Eric~J. Rawdon, Andrzej Stasiak, Kenneth~C.
  Millett, Piotr Su\l{}kowski, and Joanna~I. Su\l{}kowska.
\newblock Knotprot: a database of proteins with knots and slipknots.
\newblock {\em Nucleic Acids Research}, 43(D1):D306--D314, 2015.

\bibitem{Jejjala:2019gm}
Vishnu Jejjala, Arjun Kar, and Onkar Parrikar.
\newblock {Deep learning the hyperbolic volume of a knot}.
\newblock {\em Phys. Lett. B}, 799:135033, 2019.

\bibitem{torusenergy}
Denise Kim and Rob Kusner.
\newblock Torus knots extremizing the {M}\"obius energy.
\newblock {\em Experiment. Math.}, 2(1):1--9, 1993.

\bibitem{Levitt:2019ut}
Jesse S.~F. Levitt, Mustafa Hajij, and Radmila Sazdanovic.
\newblock Big data approaches to knot theory: Understanding the structure of
  the {J}ones polynomial.
\newblock Preprint, {\tt arXiv:1912.10086 [math.GT]}, 2019.

\bibitem{mansfield1}
Marc~L. Mansfield.
\newblock Are there knots in proteins?
\newblock {\em Nat. Struct. Biol.}, 1:213--214, 1994.

\bibitem{mansfield2}
Marc~L. Mansfield.
\newblock Fit to be tied.
\newblock {\em Nat. Struct. Biol.}, 4:166--167, 1997.

\bibitem{Millett:2016jz}
Kenneth~C. Millett.
\newblock {Knots in knots: A study of classical knot diagrams}.
\newblock {\em J. Knot Theory Ramifications}, 25(9):1641013, 2016.

\bibitem{mdsmacro}
Kenneth~C. Millett, Akos Dobay, and Andrzej Stasiak.
\newblock Linear random knots and their scaling behavior.
\newblock {\em Macromolecules}, 38(2):601--606, 2005.

\bibitem{spatial}
Kenneth~C. Millett and Eric~J. Rawdon.
\newblock Energy, ropelength, and other physical aspects of equilateral knots.
\newblock {\em J. Comput. Phys.}, 186(2):426--456, 2003.

\bibitem{Millett:2017iy}
Kenneth~C. Millett and Alex Rich.
\newblock {More knots in knots: A study of classical knot diagrams}.
\newblock {\em J. Knot Theory Ramifications}, 26(8):1750046, 2017.

\bibitem{PT}
J{\'o}zef~H. Przytycki and Pawe{\l} Traczyk.
\newblock Invariants of links of {C}onway type.
\newblock {\em Kobe J. Math.}, 4(2):115--139, 1987.

\bibitem{subknots}
Eric~J. Rawdon, Kenneth~C. Millett, and Andrzej Stasiak.
\newblock Subknots in ideal knots, random knots, and knotted proteins.
\newblock {\em Sci. Rep.}, 5:8928, 2015.

\bibitem{Rayleigh:1919do}
Lord Rayleigh.
\newblock {On the problem of random vibrations, and of random flights in one,
  two, or three dimensions}.
\newblock {\em Philosophical Magazine}, 37(220):321--347, 1919.

\bibitem{knotplot}
Robert~G. Scharein.
\newblock {KnotPlot}.
\newblock http://\-www.\-knotplot.\-com, 1998.
\newblock Program for drawing, visualizing, manipulating, and energy minimizing
  knots.

\bibitem{ourpnas}
Joanna~I. Su\l{}kowska, Eric~J. Rawdon, Kenneth~C. Millett, Jose~N. Onuchic,
  and Andrzej Stasiak.
\newblock Conservation of complex knotting and slipknotting patterns in
  proteins.
\newblock {\em Proc. Natl. Acad. Sci. USA}, 109(26):E1715--E1723, 2012.

\bibitem{unraveller}
Morwen Thistlethwaite.
\newblock unraveller, 2009.
\newblock Software to simplify knot diagrams.

\bibitem{turaev2012}
Vladimir Turaev.
\newblock Knotoids.
\newblock {\em Osaka J. Math.}, 49(1):195--223, 2012.

\bibitem{Vandans:2020cp}
Olafs Vandans, Kaiyuan Yang, Zhongtao Wu, and Liang Dai.
\newblock {Identifying knot types of polymer conformations by machine
  learning}.
\newblock {\em Phys. Rev. E}, 101:022502, 2020.

\bibitem{virnauweb}
Peter Virnau.
\newblock Protein knot server.
\newblock http://knots.mit.edu/.
\newblock Server that allows users to check for knots in {PDB} entries or other
  uploaded structures.

\end{thebibliography}

\pagebreak

\appendix
\section{Appendix}


The accuracies of the $\operatorname{PR}$, $\operatorname{PU}$, $\operatorname{SR}$, and $\operatorname{SU}$ classifiers for each of the knot types in our data set, shown as a function of subchain length.

\vfill

{\ }\hfill
\includegraphics[width=0.40\textwidth]{a01.fourmethods.pdf}
\hfill{\ }

\vfill

{\ }\hfill
\includegraphics[width=0.40\textwidth]{p31.fourmethods.pdf}
\hfill
\includegraphics[width=0.40\textwidth]{a41.fourmethods.pdf}
\hfill{\ }

\vfill

{\ }\hfill
\includegraphics[width=0.40\textwidth]{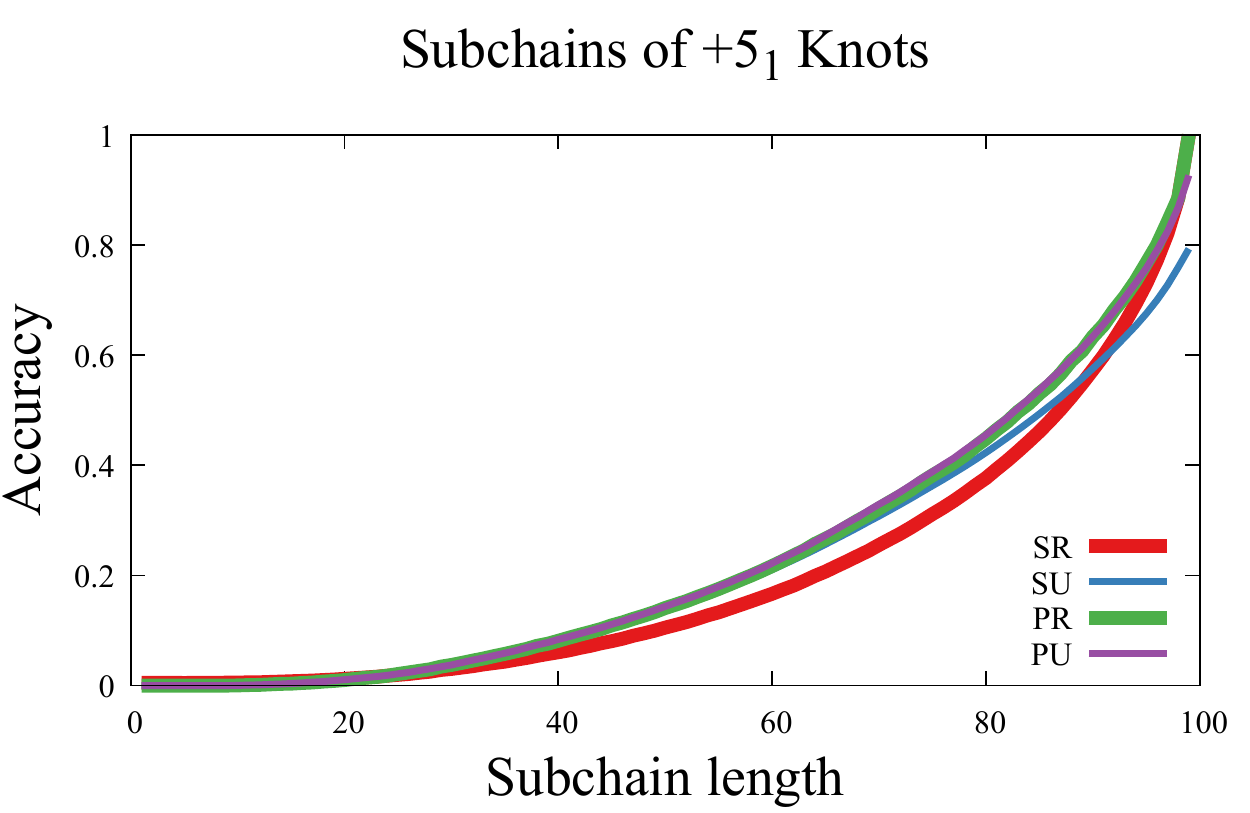}
\hfill
\includegraphics[width=0.40\textwidth]{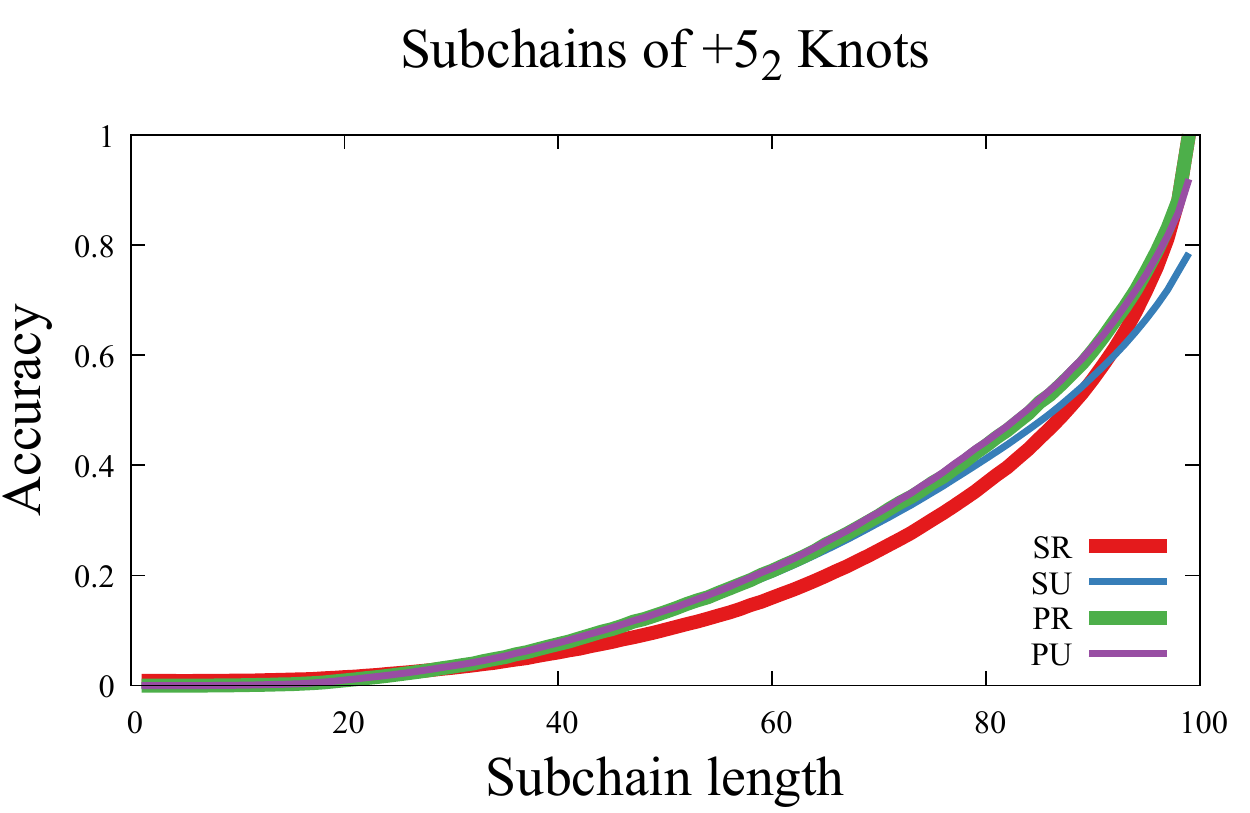}
\hfill{\ }

\vfill

{\ }\hfill
\includegraphics[width=0.40\textwidth]{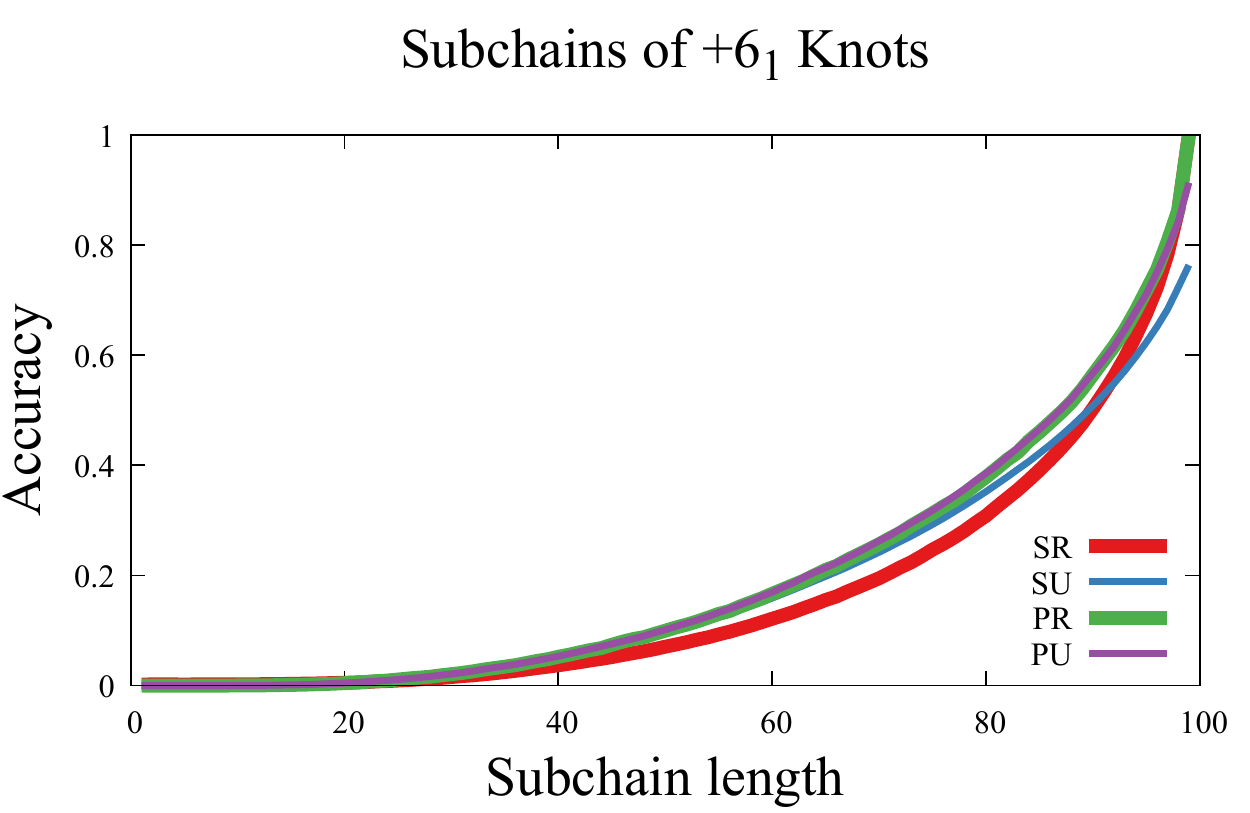}
\hfill
\includegraphics[width=0.40\textwidth]{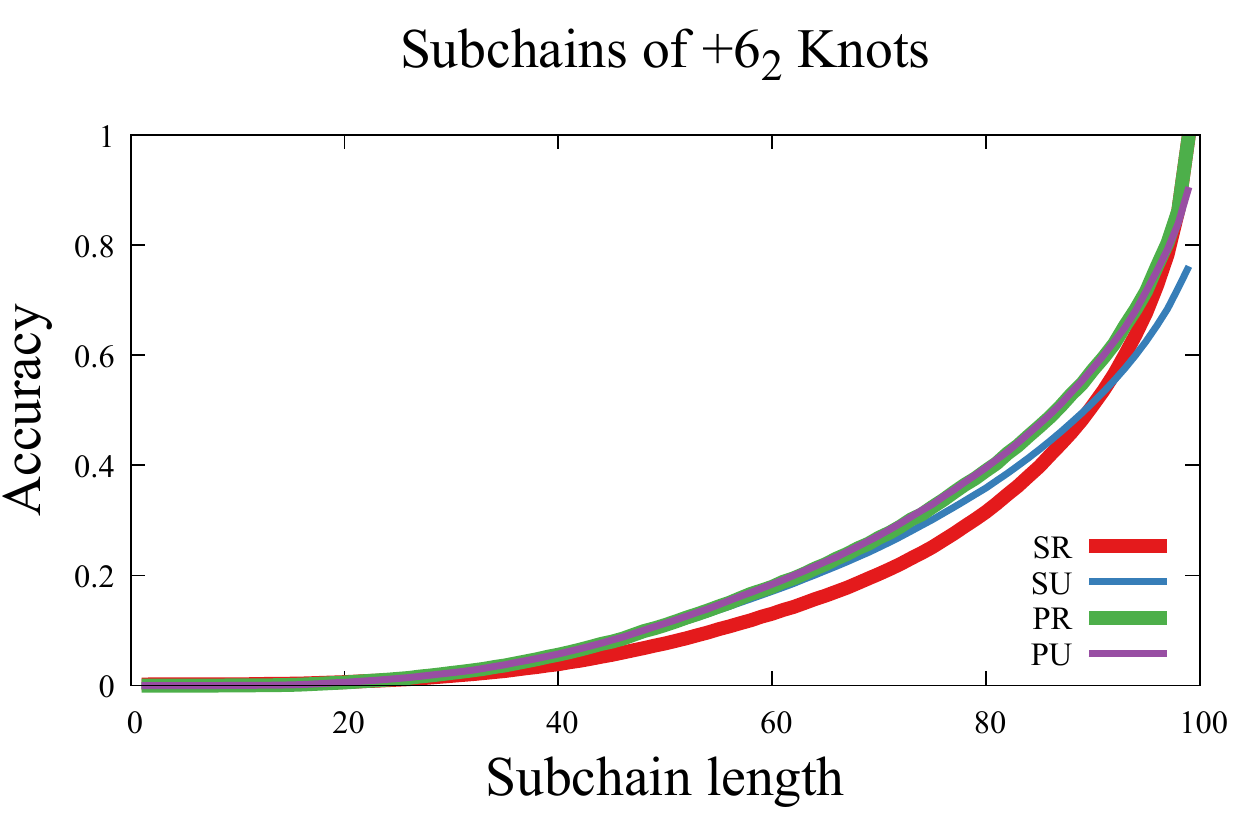}
\hfill{\ }

\vfill

\pagebreak

{\ }\hfill
\includegraphics[width=0.40\textwidth]{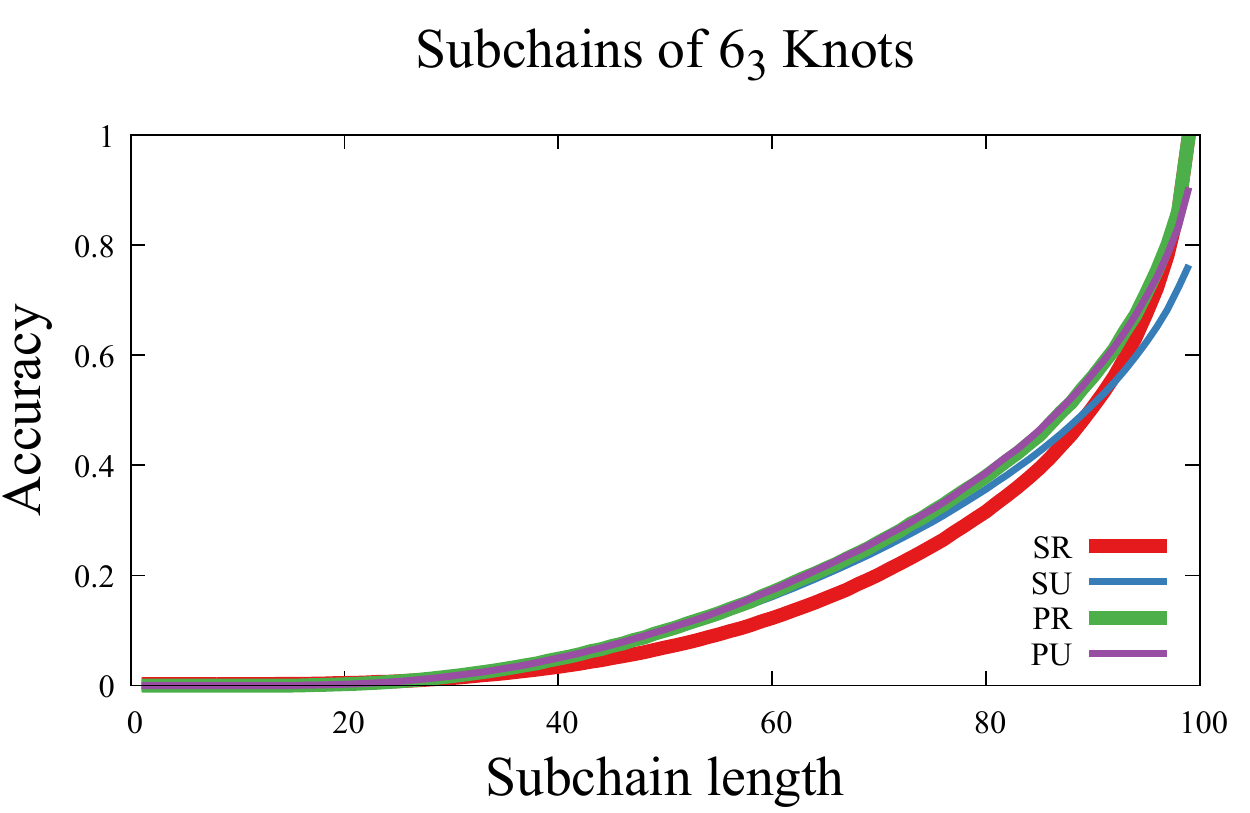}
\hfill
\includegraphics[width=0.40\textwidth]{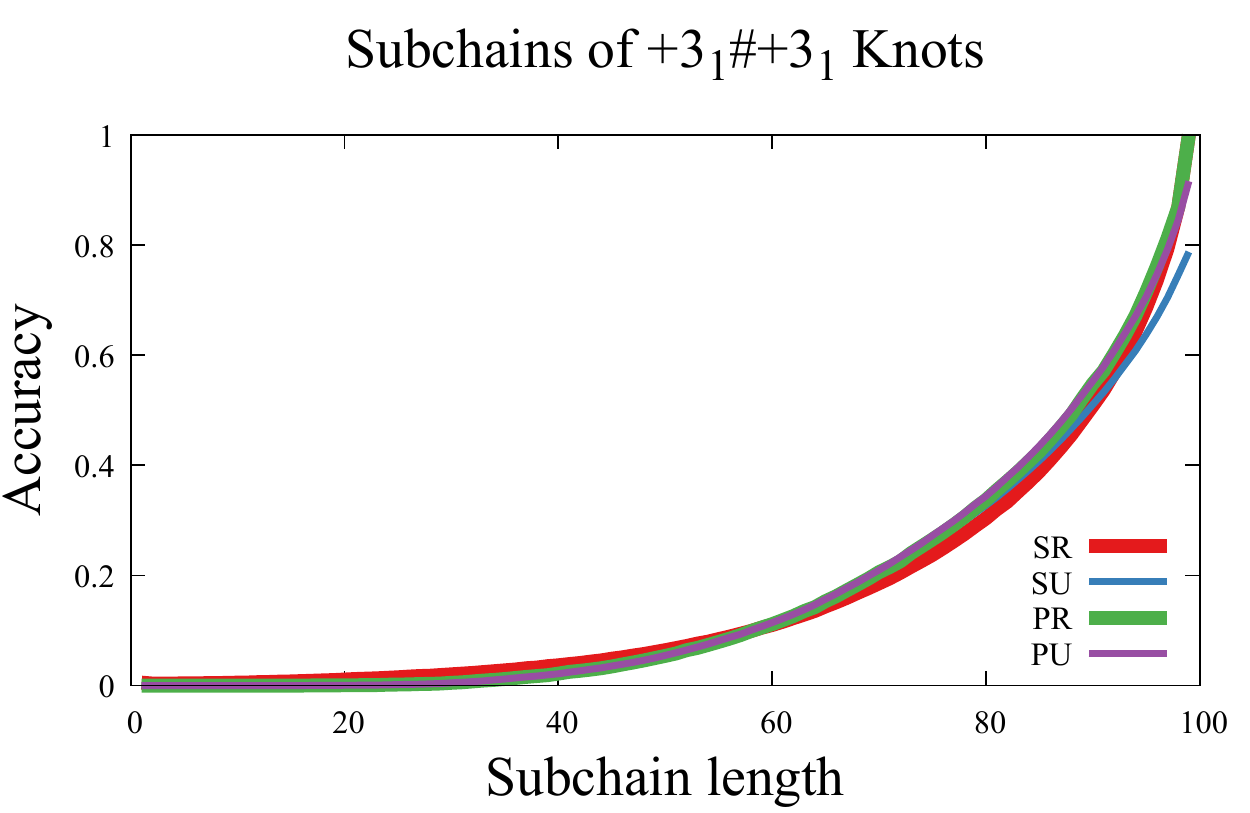}
\hfill{\ }

{\ }\hfill
\includegraphics[width=0.40\textwidth]{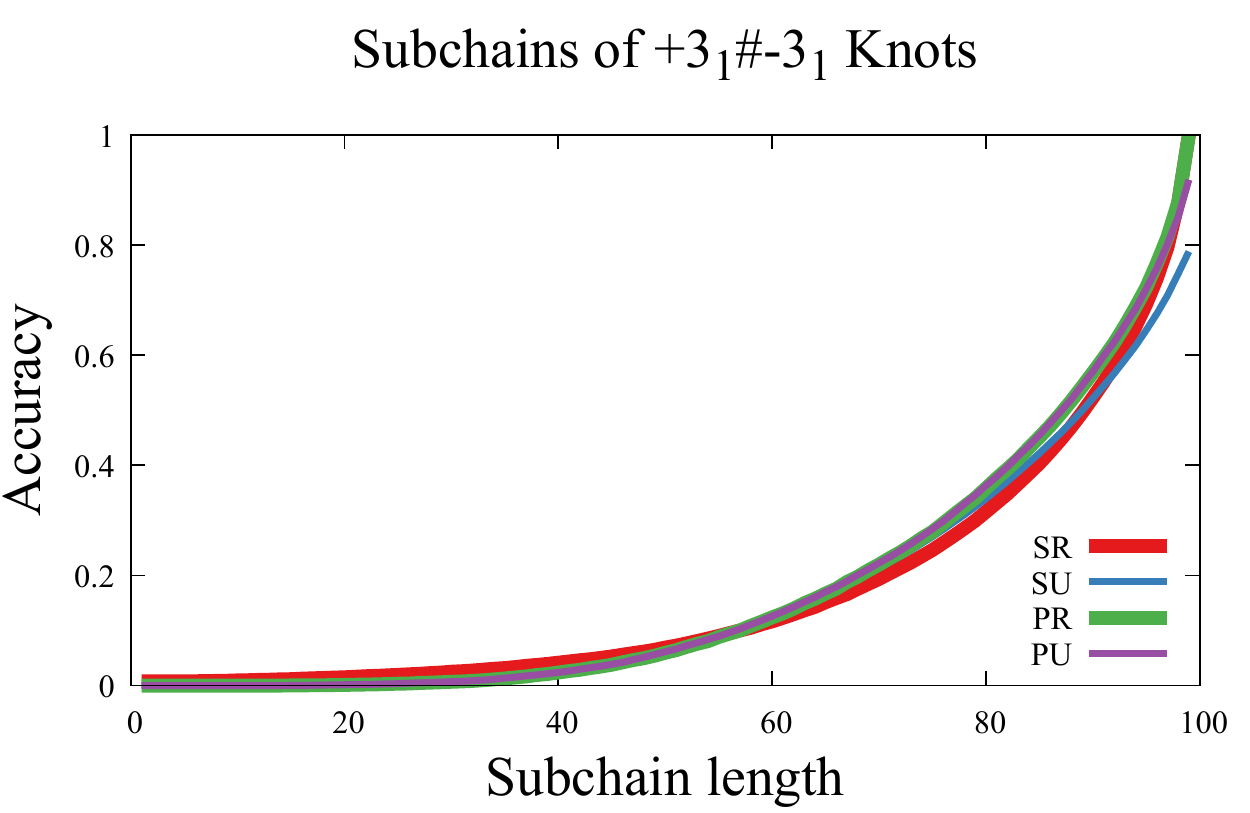}
\hfill
\includegraphics[width=0.40\textwidth]{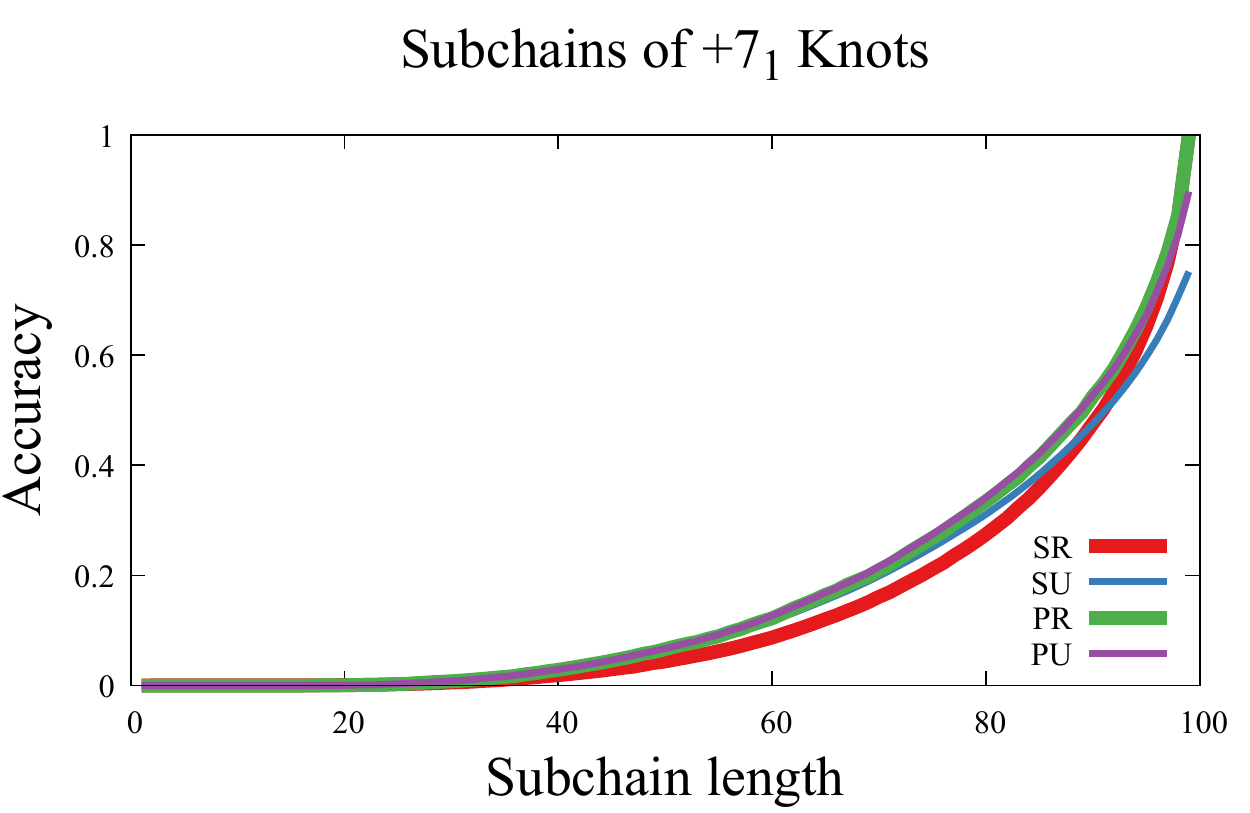}
\hfill{\ }

{\ }\hfill
\includegraphics[width=0.40\textwidth]{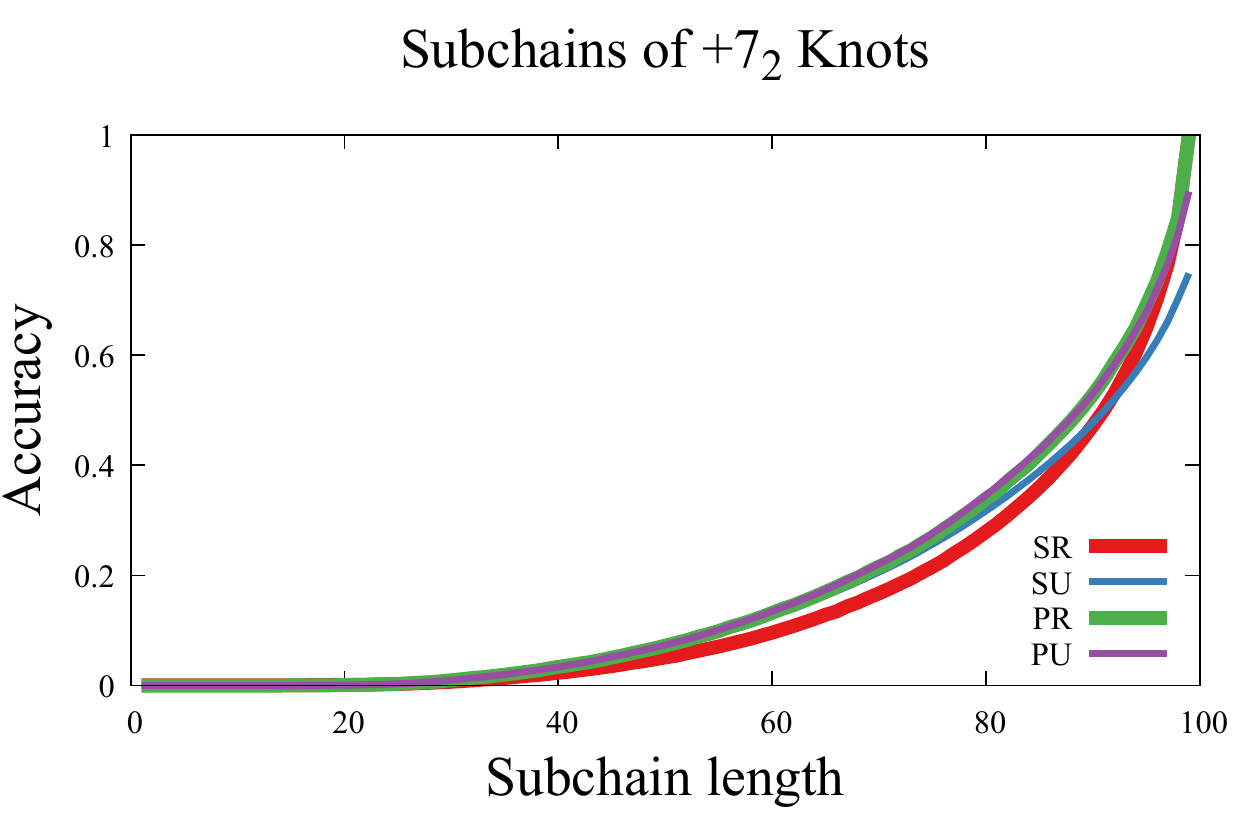}
\hfill
\includegraphics[width=0.40\textwidth]{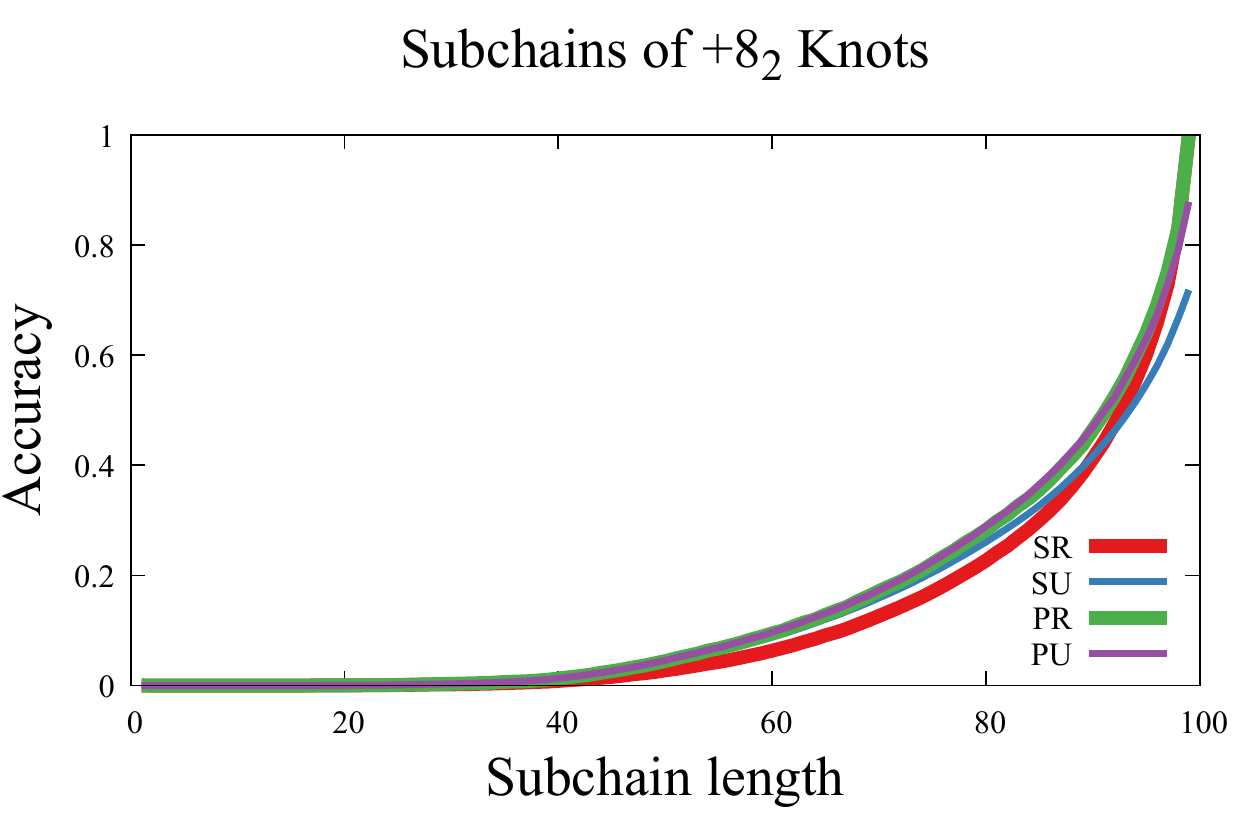}
\hfill{\ }

{\ }\hfill
\includegraphics[width=0.40\textwidth]{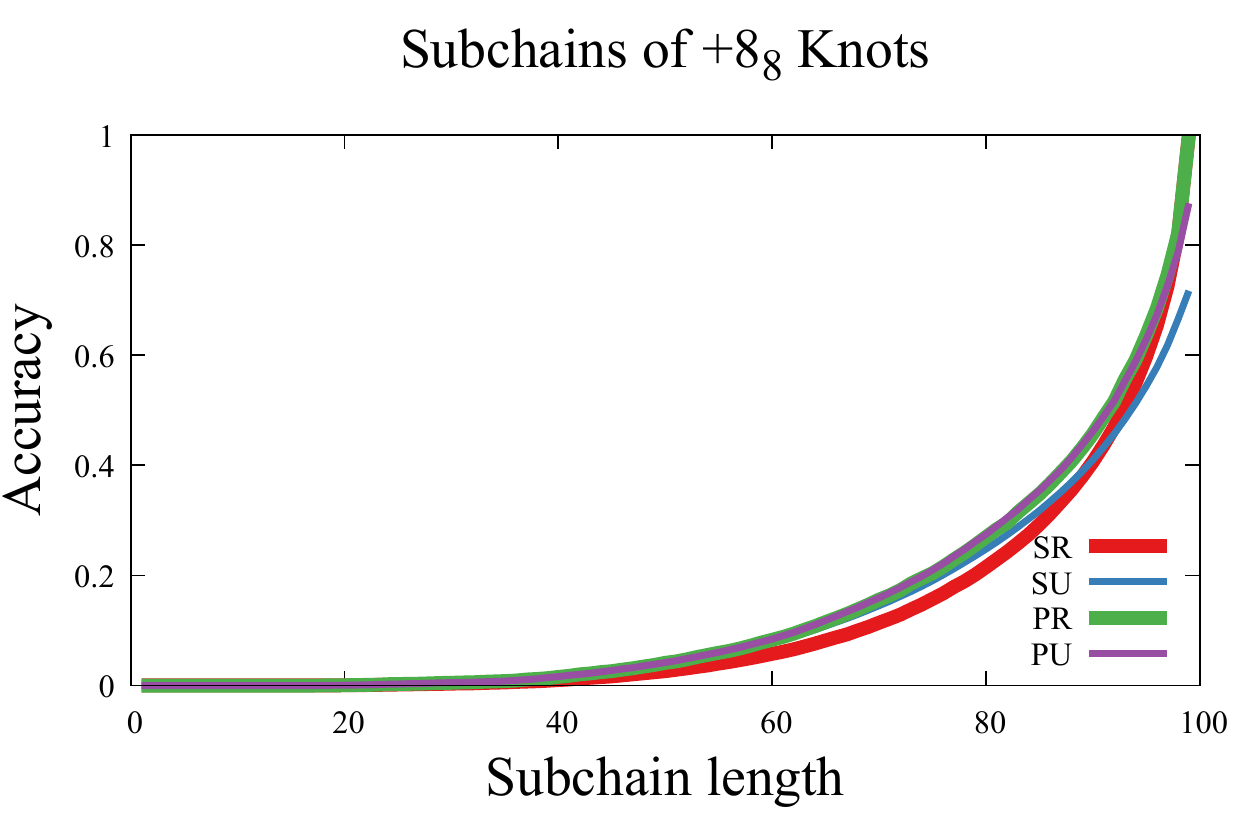}
\hfill
\includegraphics[width=0.40\textwidth]{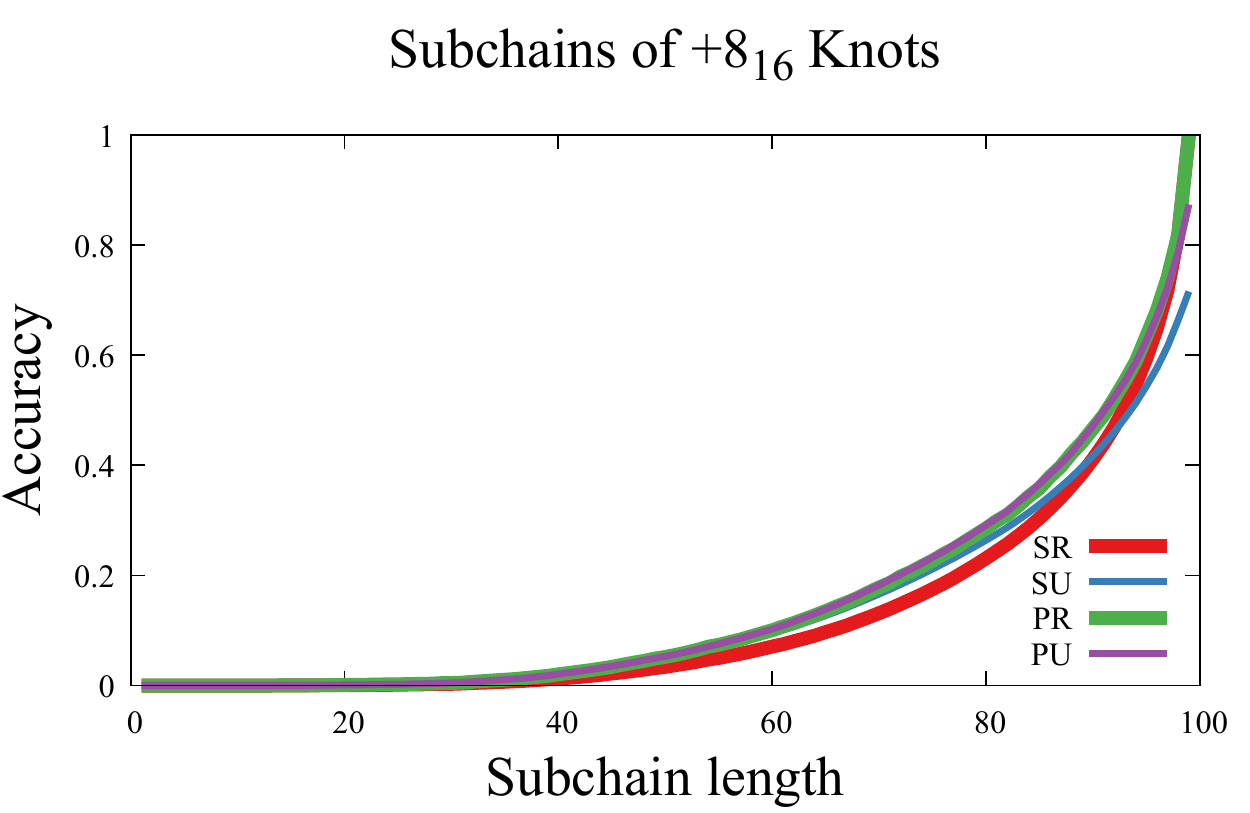}
\hfill{\ }

{\ }\hfill
\includegraphics[width=0.40\textwidth]{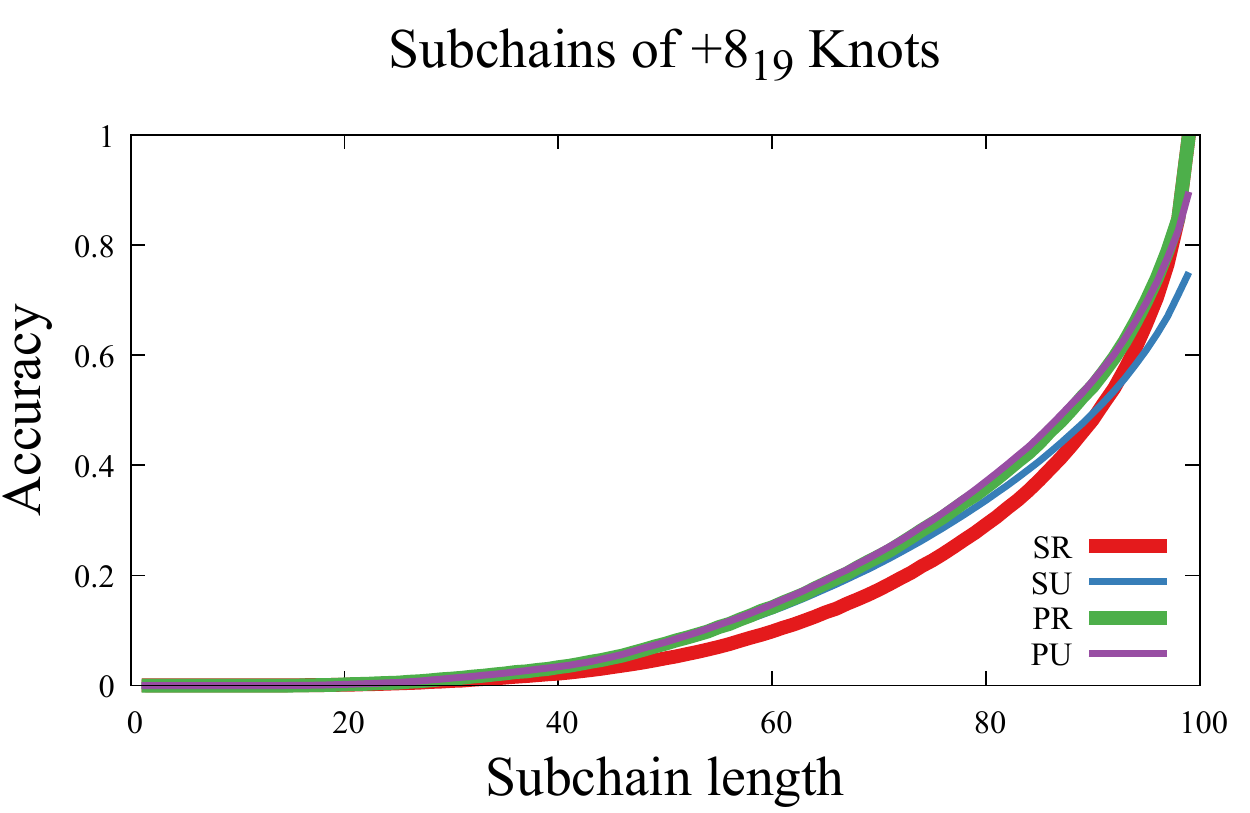}
\hfill
\includegraphics[width=0.40\textwidth]{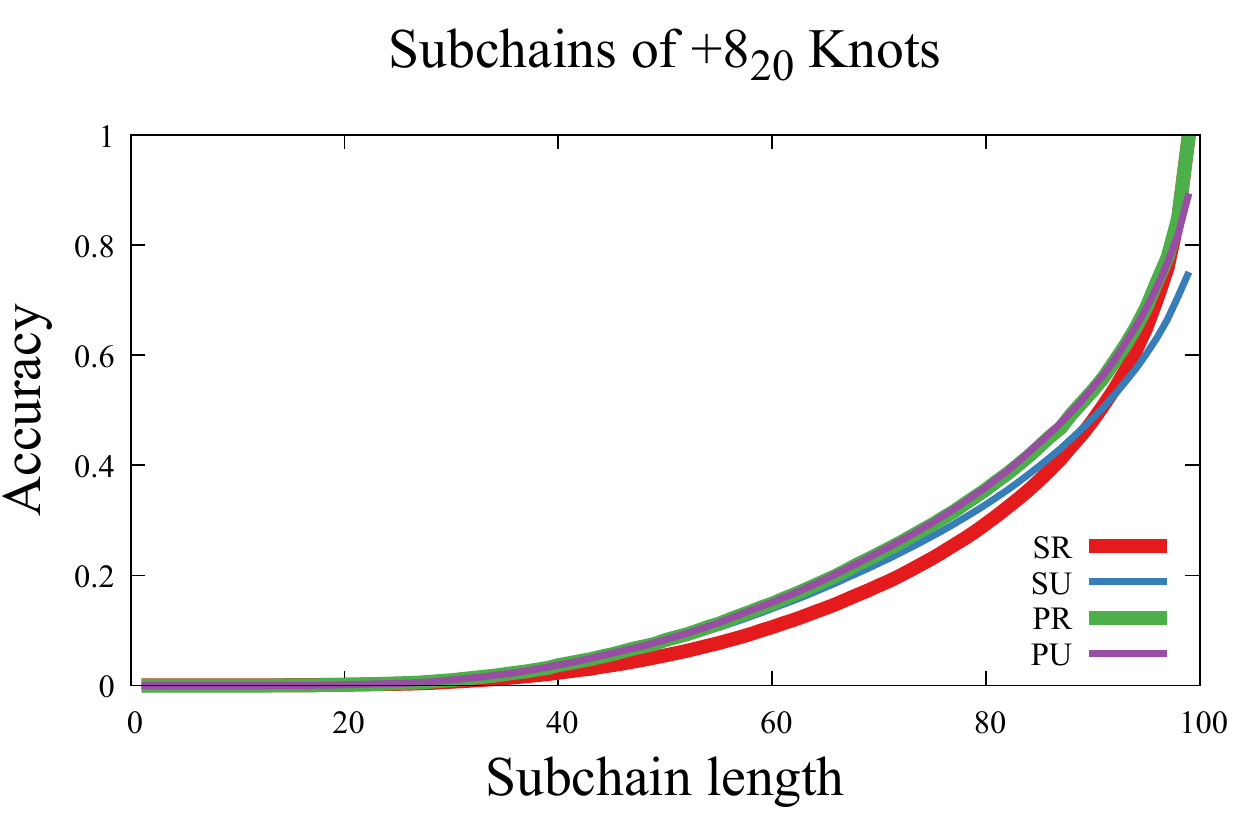}
\hfill{\ }

\pagebreak

\vfill

{\ }\hfill
\includegraphics[width=0.40\textwidth]{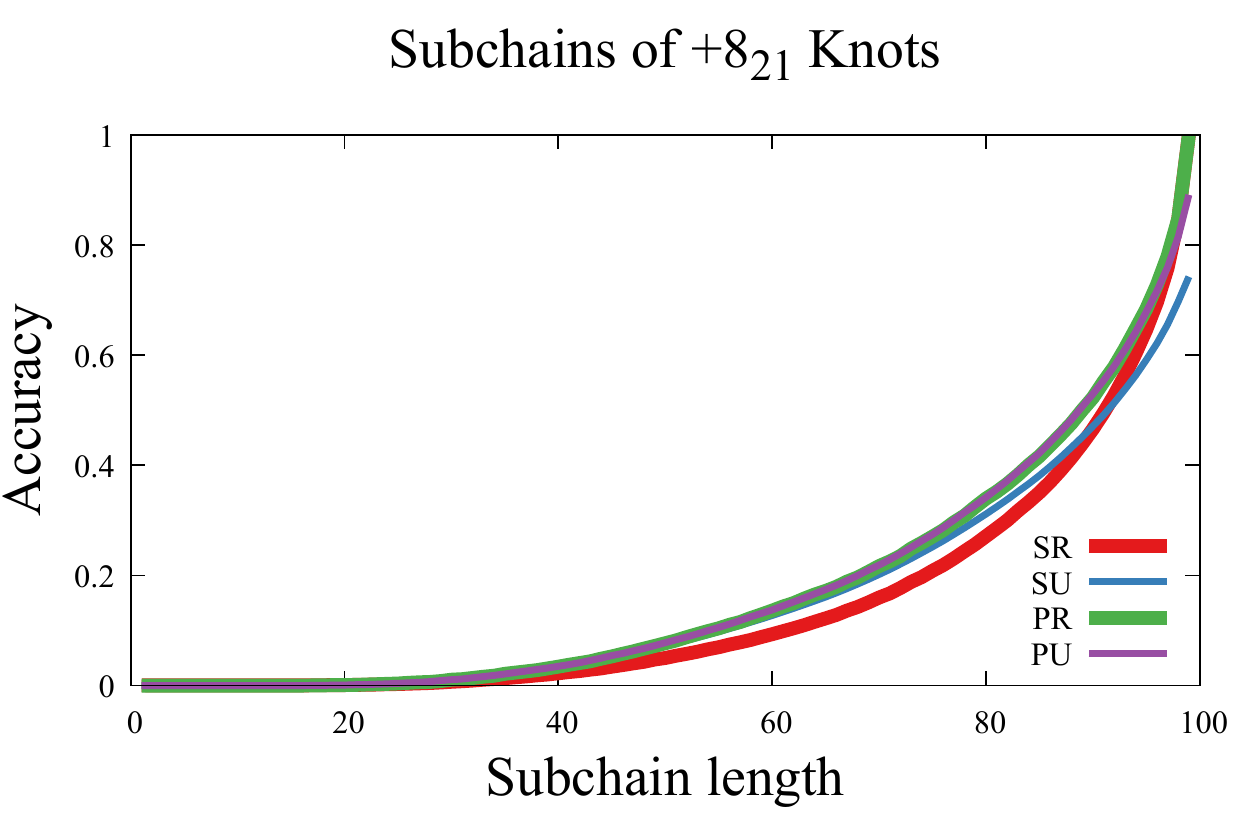}
\hfill
\includegraphics[width=0.40\textwidth]{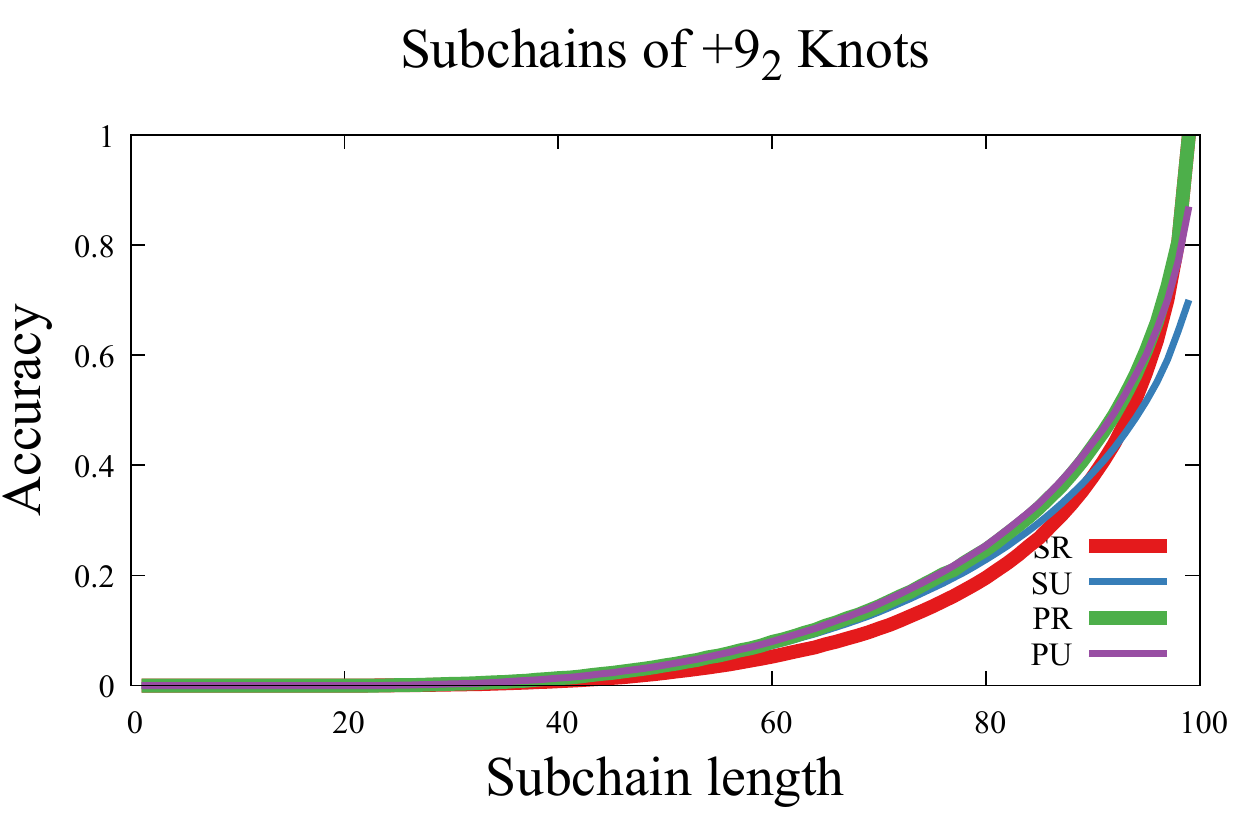}
\hfill{\ }

\vfill

{\ }\hfill
\includegraphics[width=0.40\textwidth]{p915.fourmethods.pdf}
\hfill
\includegraphics[width=0.40\textwidth]{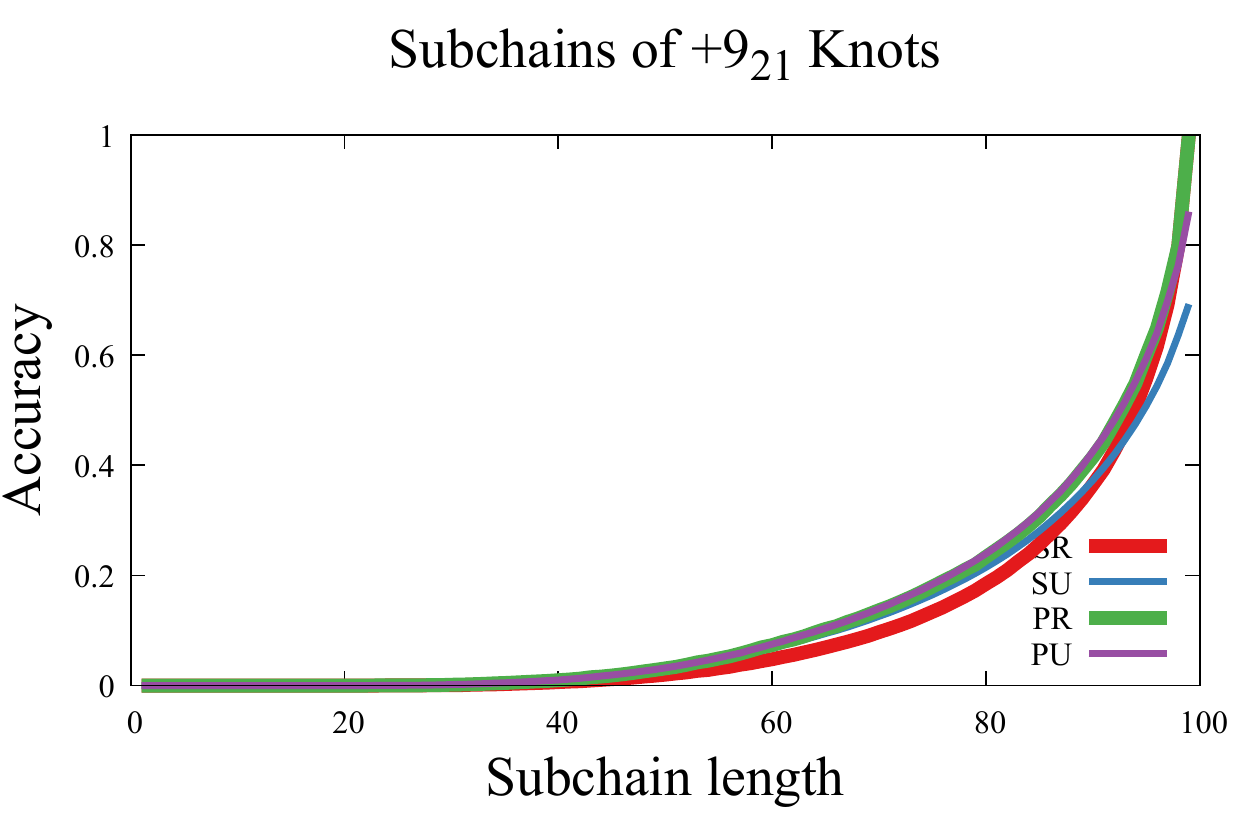}
\hfill{\ }

\vfill

{\ }\hfill
\includegraphics[width=0.40\textwidth]{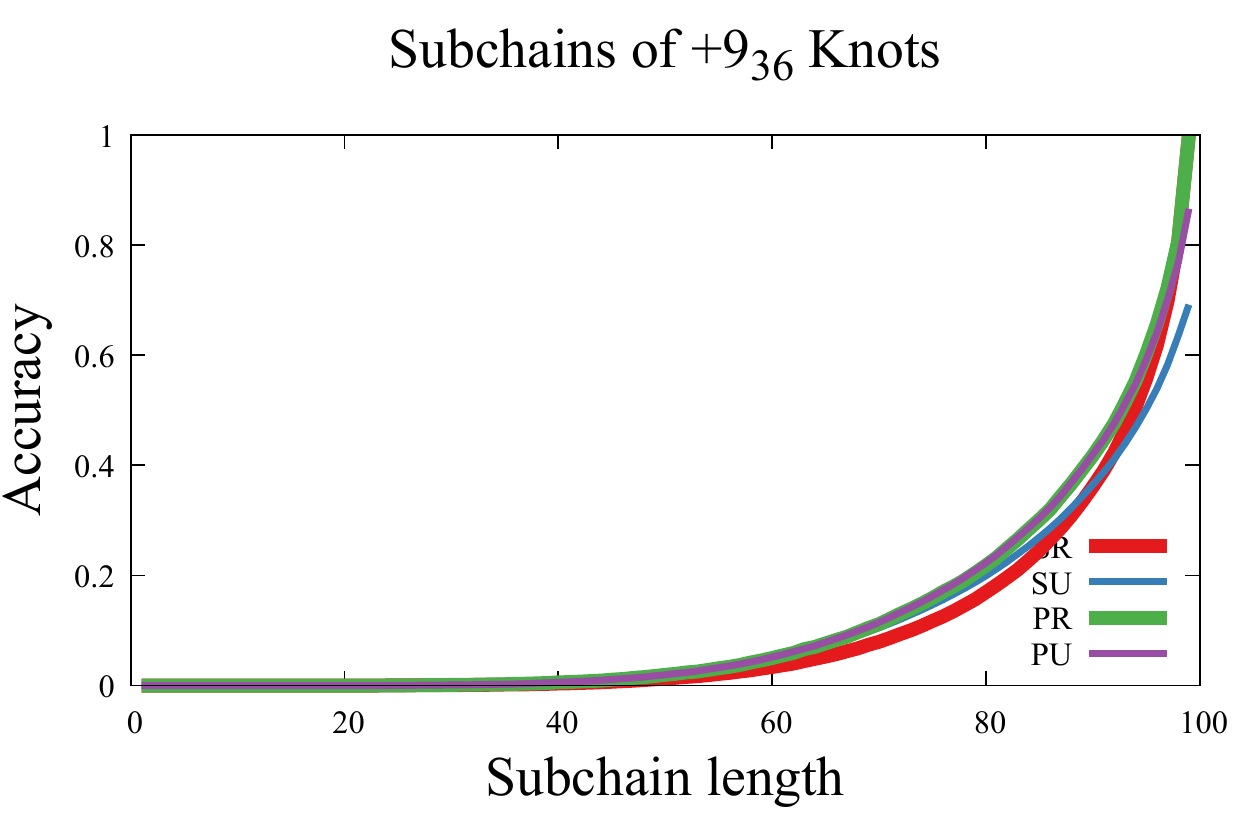}
\hfill
\includegraphics[width=0.40\textwidth]{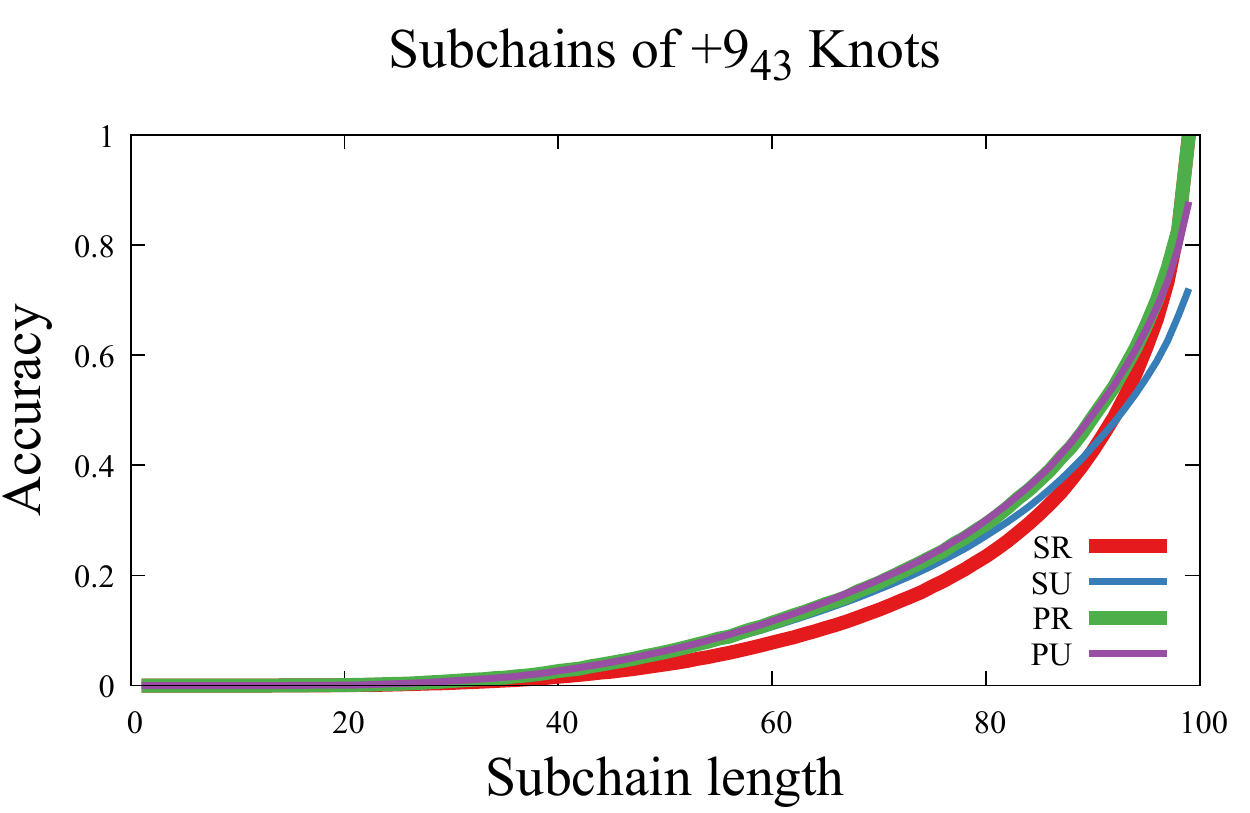}
\hfill{\ }

\vfill

{\ }\hfill
\includegraphics[width=0.40\textwidth]{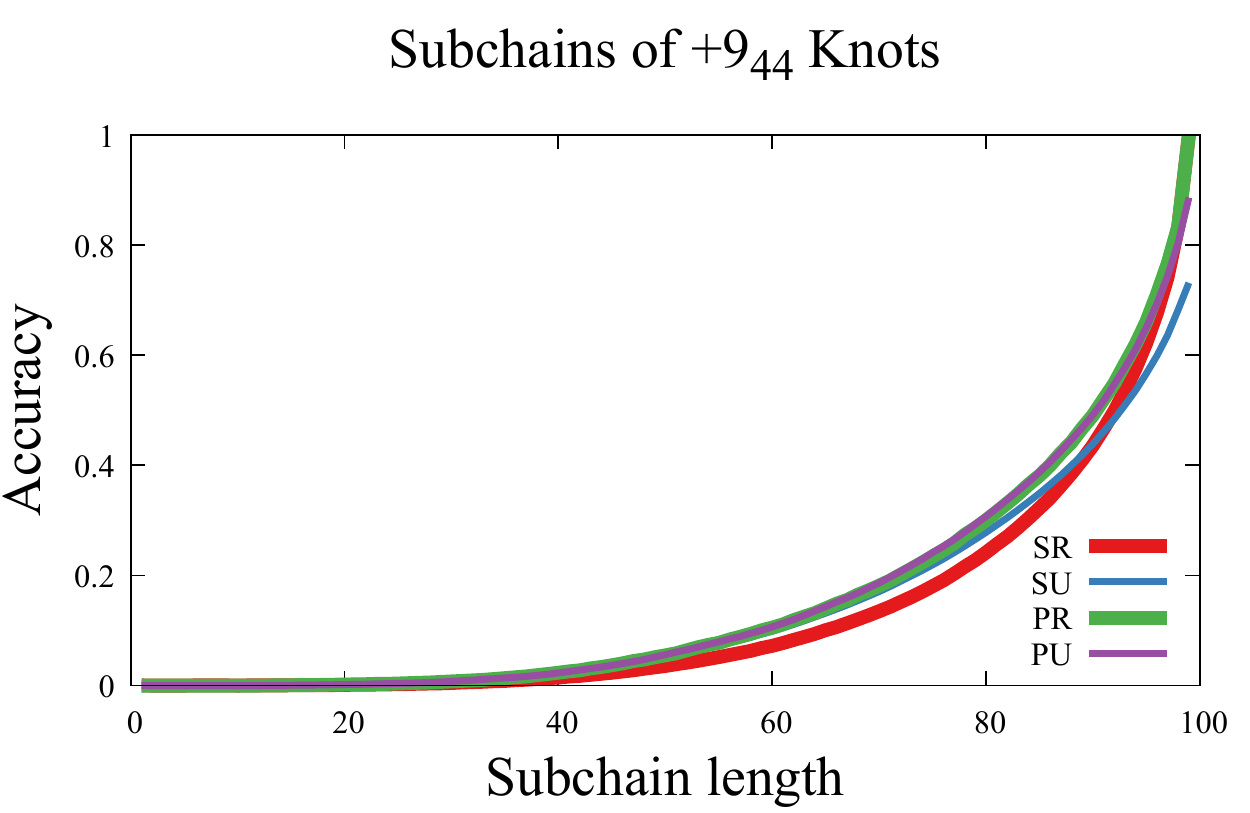}
\hfill
\includegraphics[width=0.40\textwidth]{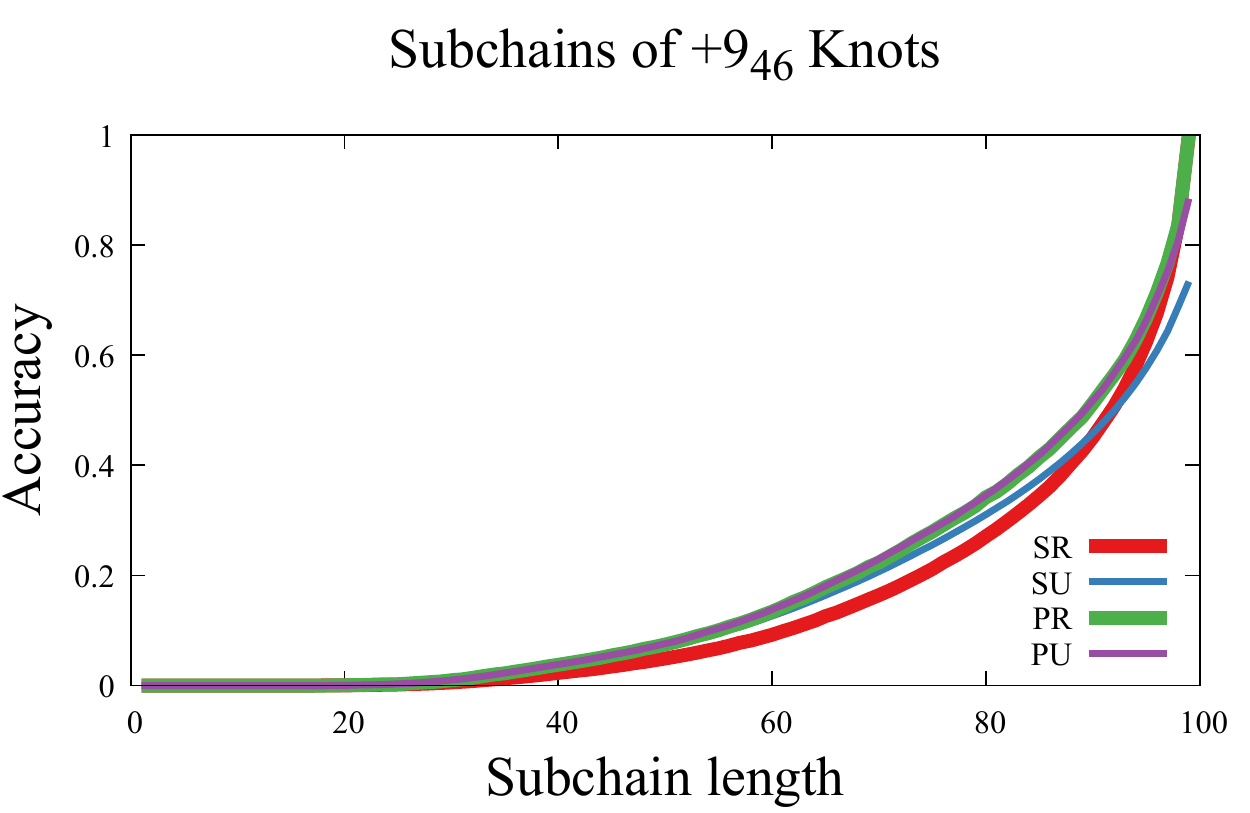}
\hfill{\ }

\vfill

\pagebreak

The PPV of the $\operatorname{PR}$, $\operatorname{PU}$, $\operatorname{SR}$, and $\operatorname{SU}$ classifiers for each of the knot types in our data set, shown as a function of subchain length.

\vfill

{\ }\hfill
\includegraphics[width=0.40\textwidth]{a01.bayes.pdf}
\hfill{\ }

\vfill

{\ }\hfill
\includegraphics[width=0.40\textwidth]{p31.bayes.pdf}
\hfill
\includegraphics[width=0.40\textwidth]{a41.bayes.pdf}
\hfill{\ }

\vfill

{\ }\hfill
\includegraphics[width=0.40\textwidth]{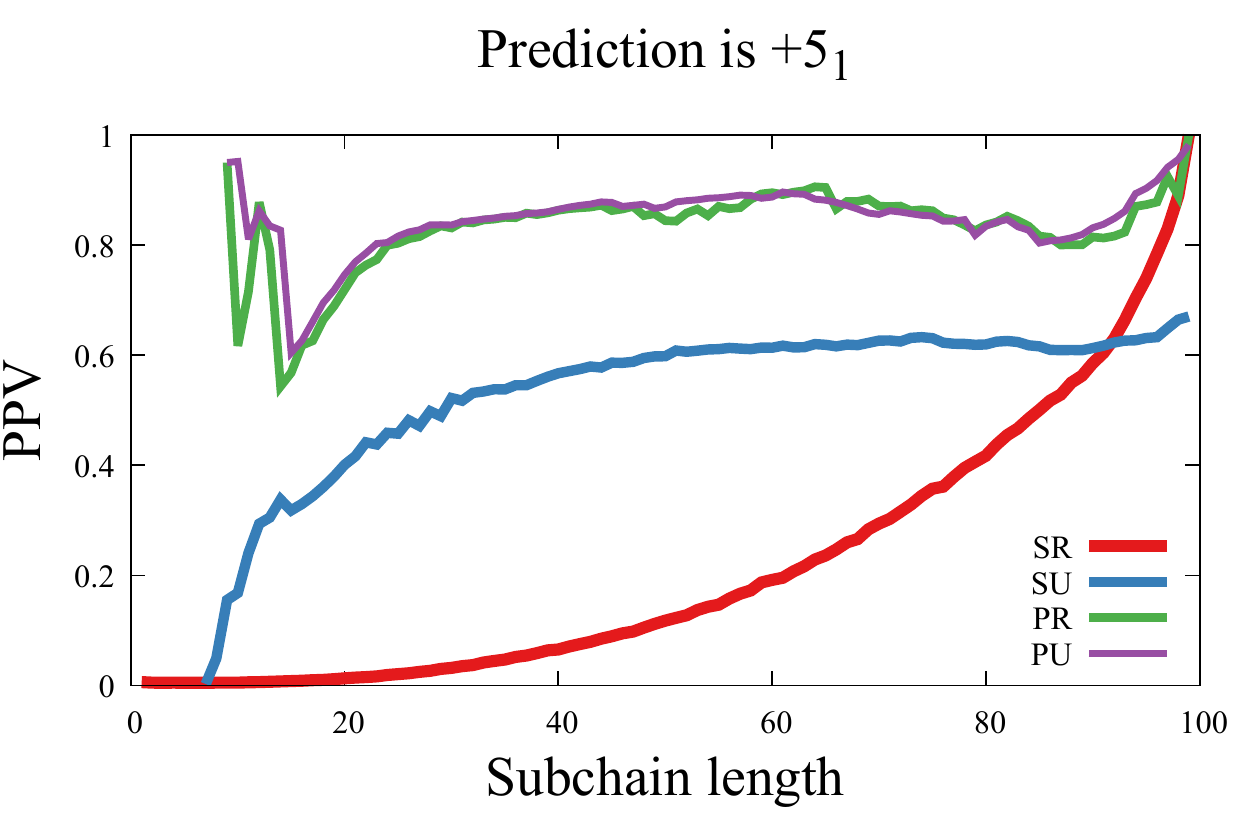}
\hfill
\includegraphics[width=0.40\textwidth]{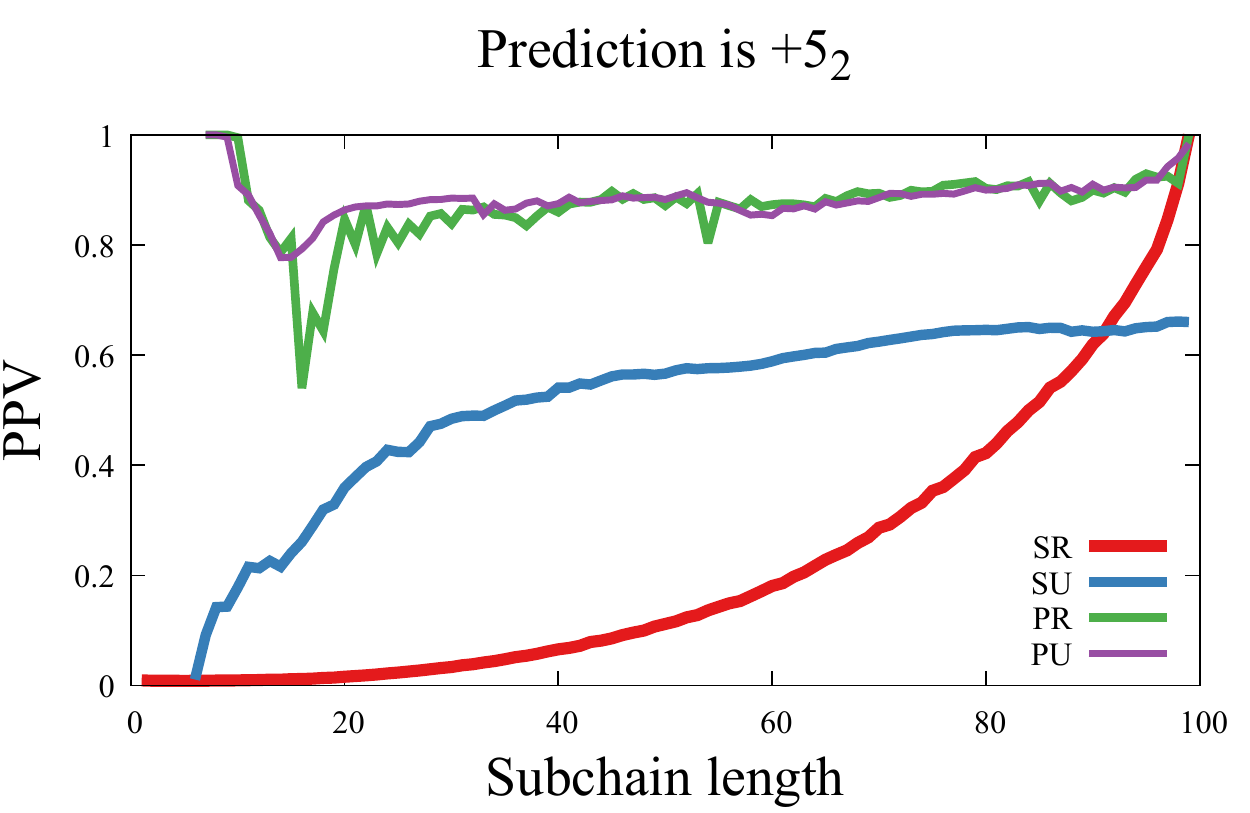}
\hfill{\ }

\vfill

{\ }\hfill
\includegraphics[width=0.40\textwidth]{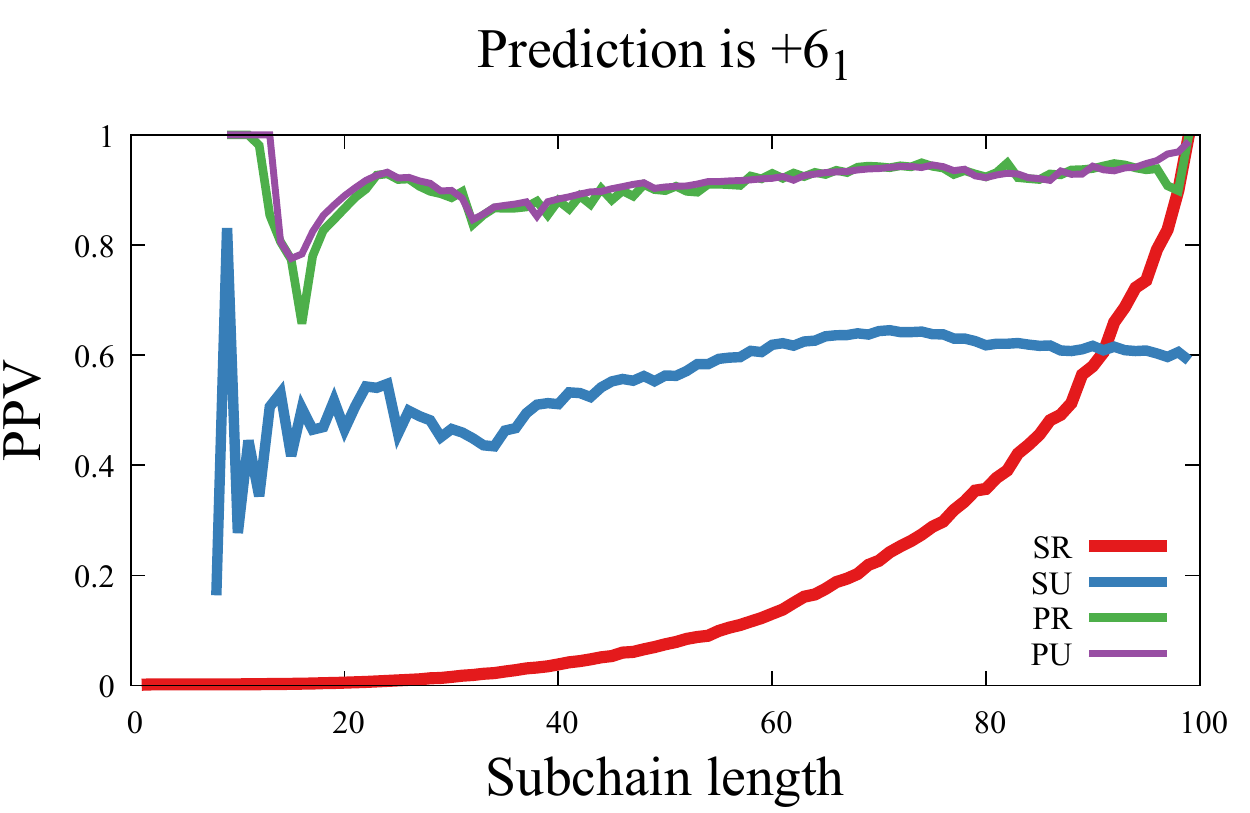}
\hfill
\includegraphics[width=0.40\textwidth]{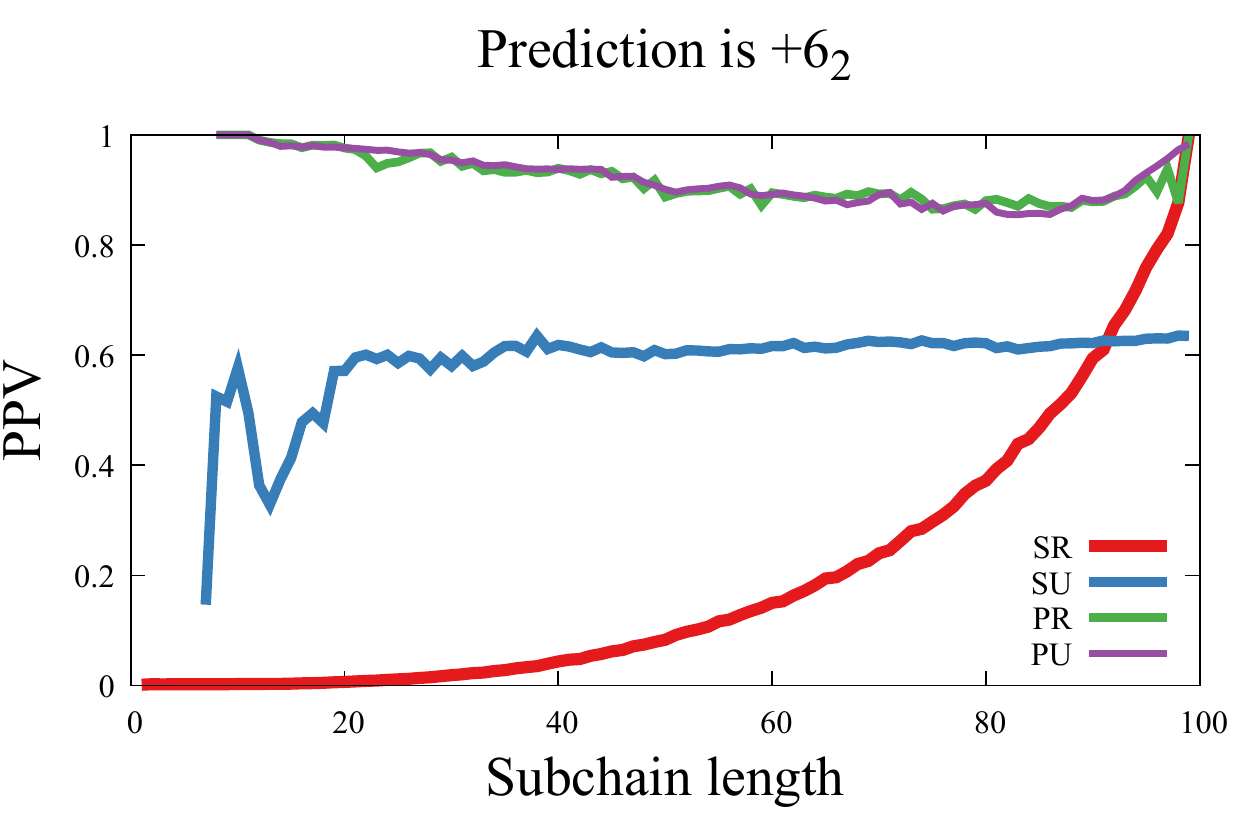}
\hfill{\ }

\vfill

\pagebreak

{\ }\hfill
\includegraphics[width=0.40\textwidth]{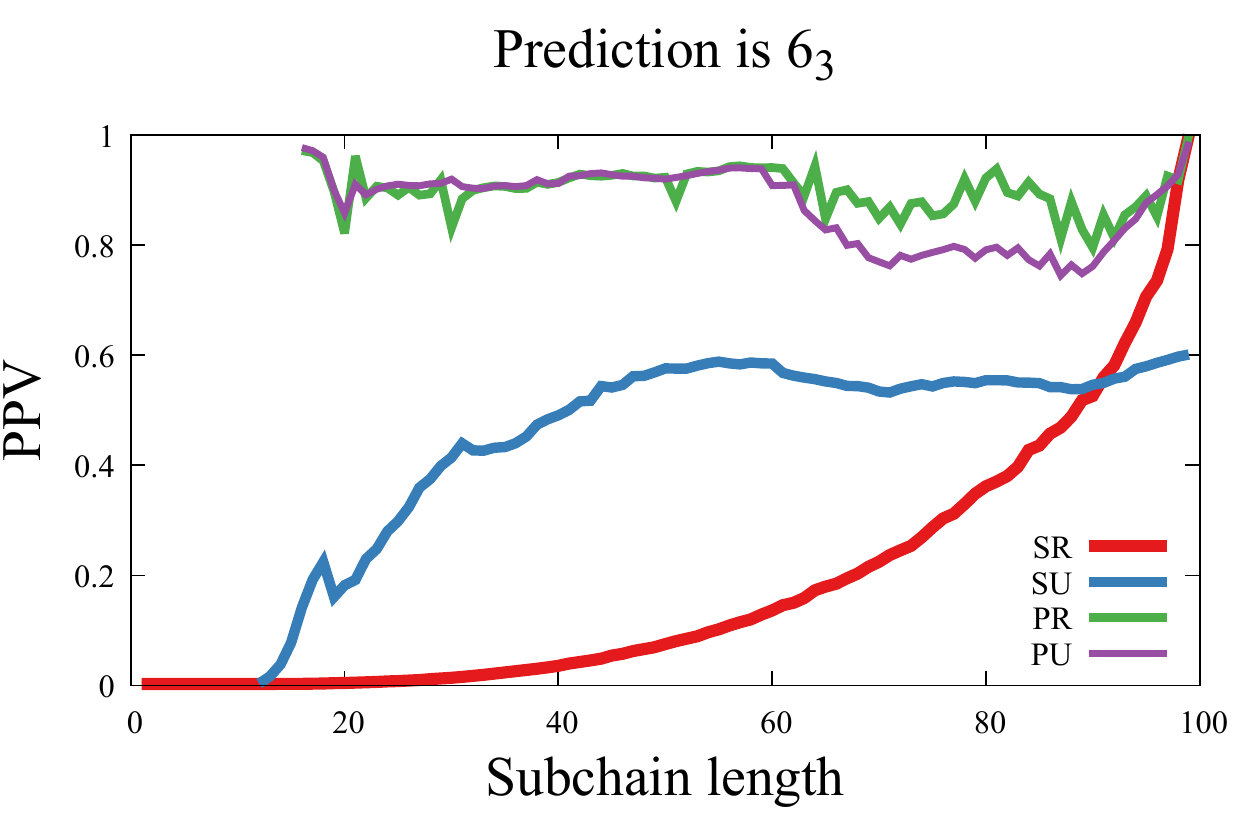}
\hfill
\includegraphics[width=0.40\textwidth]{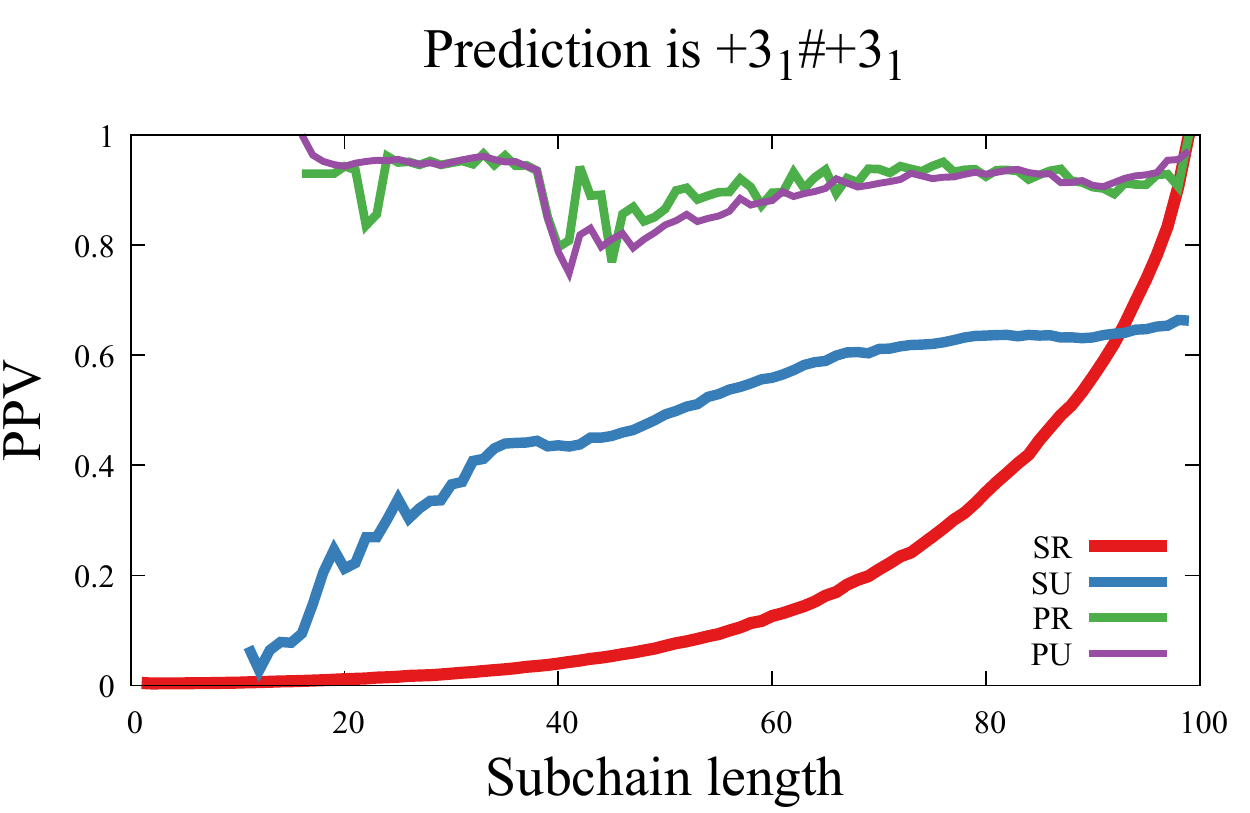}
\hfill{\ }

{\ }\hfill
\includegraphics[width=0.40\textwidth]{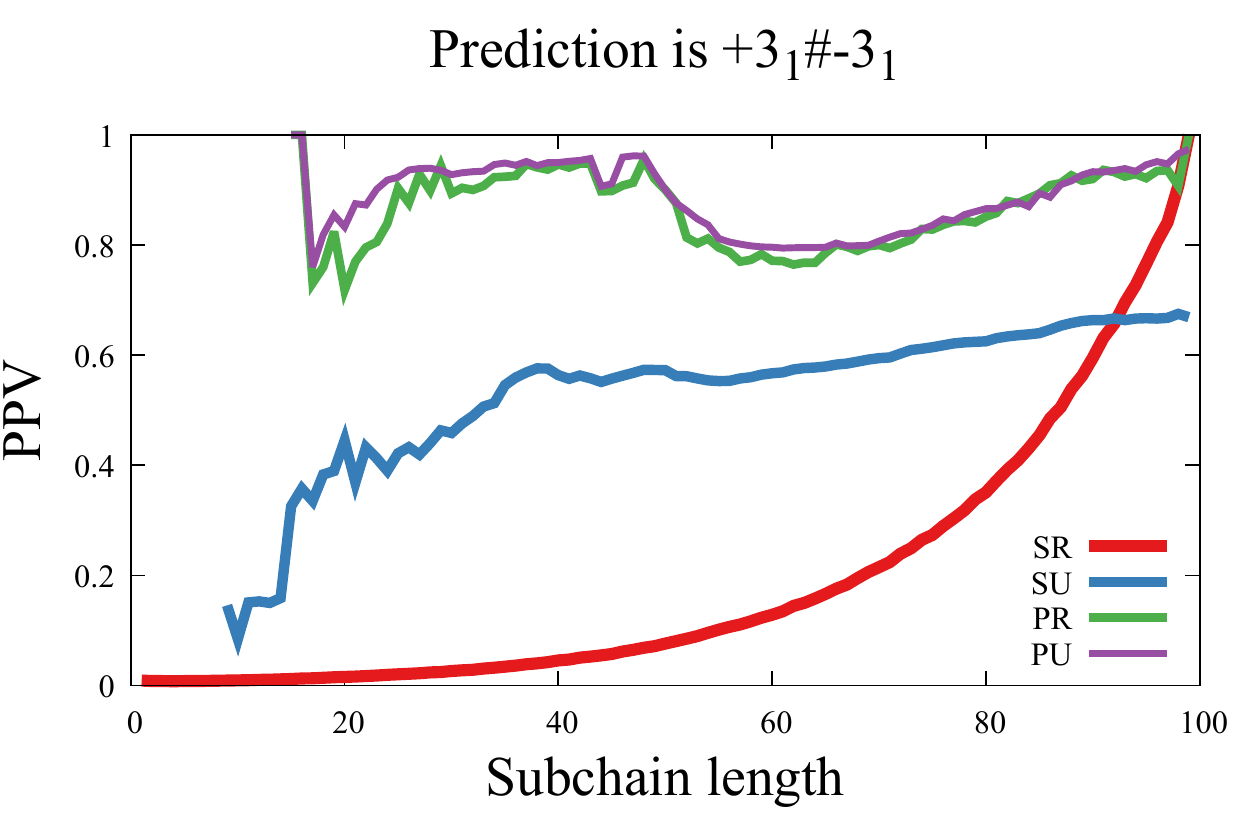}
\hfill
\includegraphics[width=0.40\textwidth]{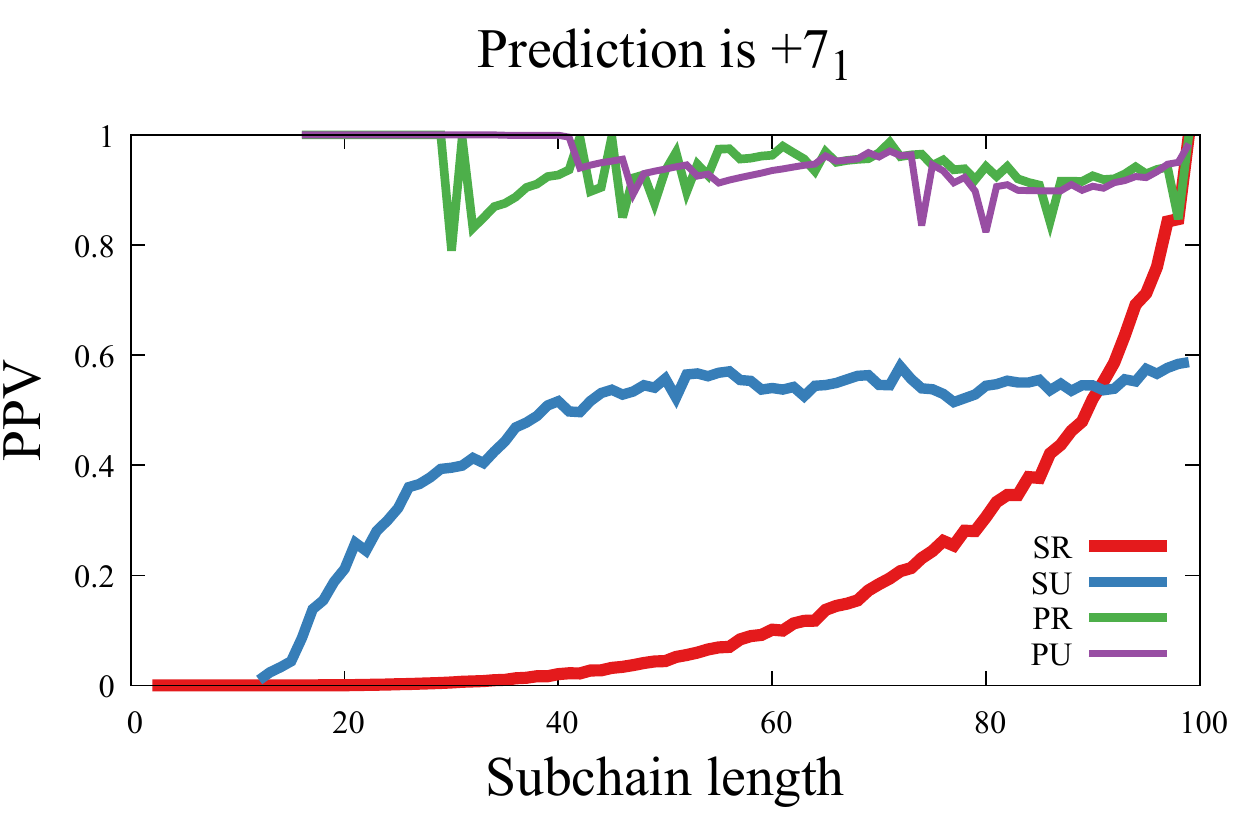}
\hfill{\ }

{\ }\hfill
\includegraphics[width=0.40\textwidth]{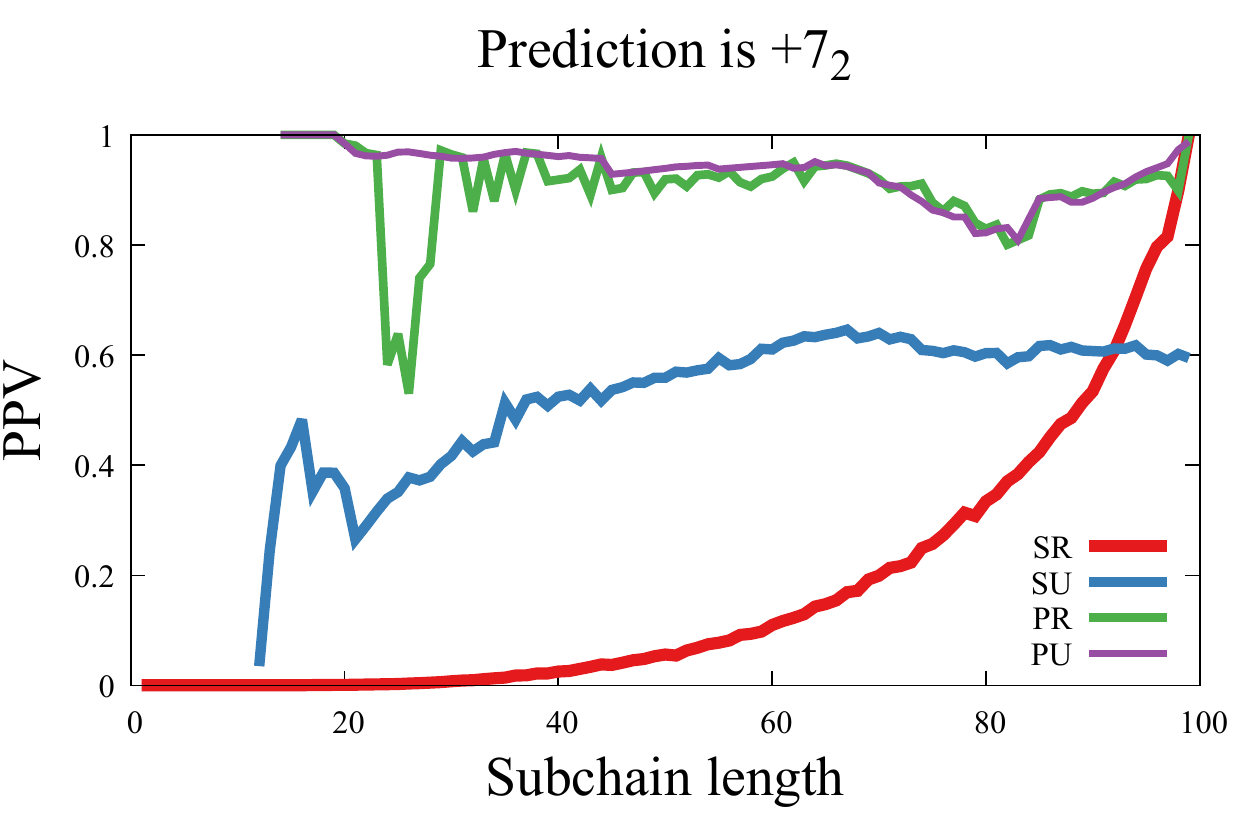}
\hfill
\includegraphics[width=0.40\textwidth]{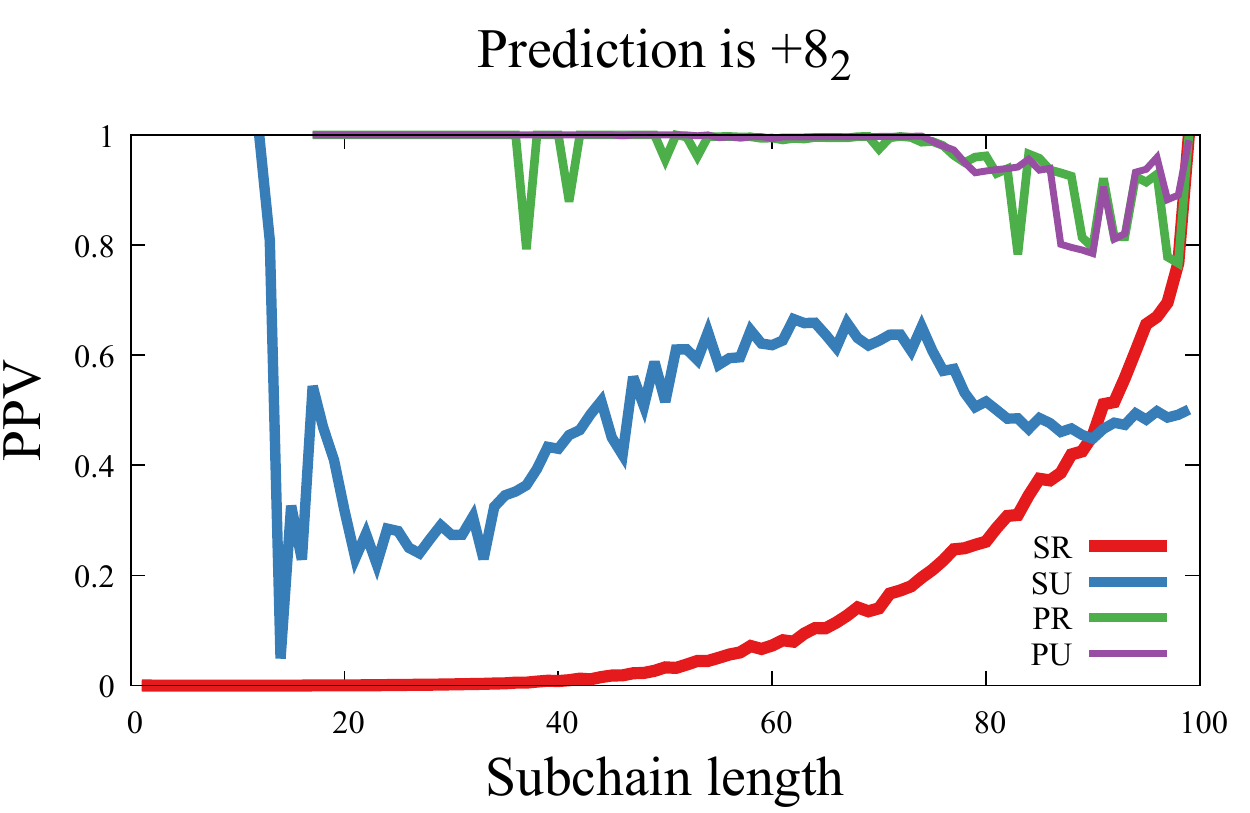}
\hfill{\ }

{\ }\hfill
\includegraphics[width=0.40\textwidth]{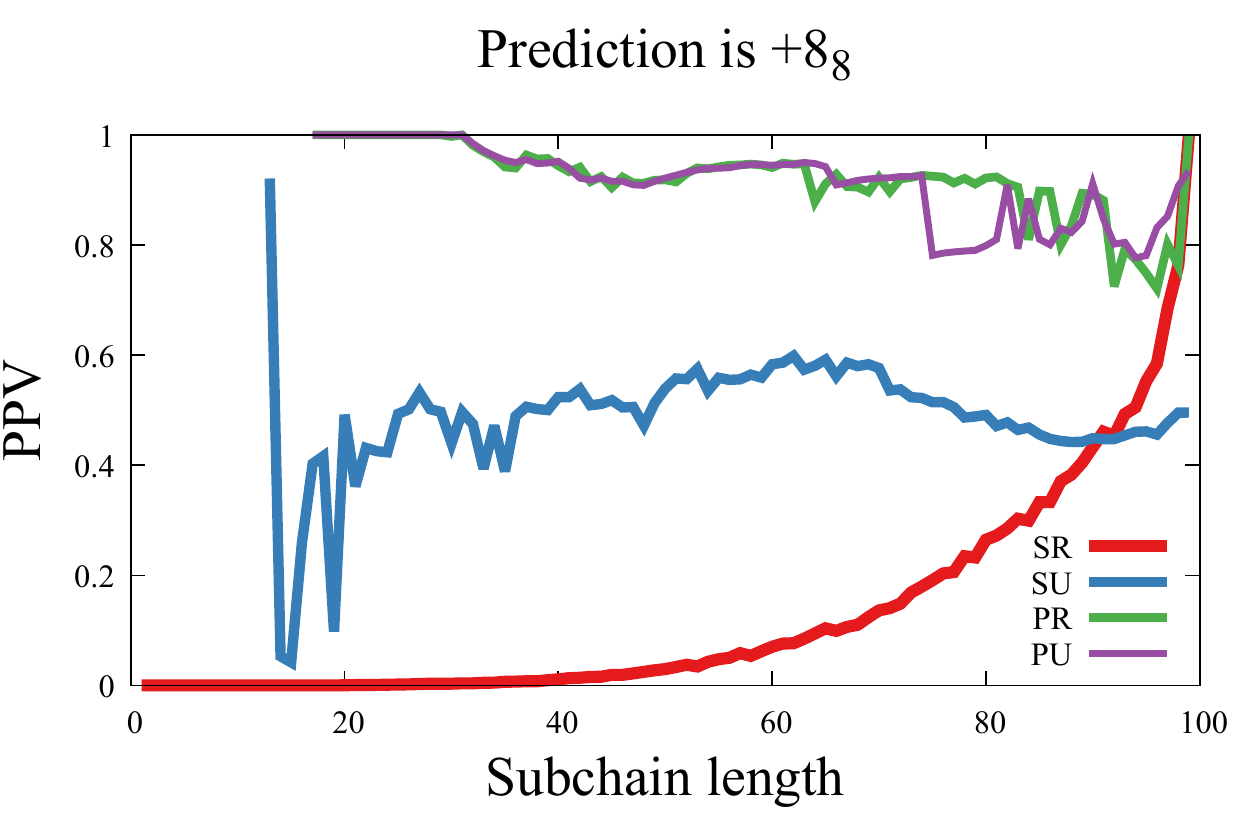}
\hfill
\includegraphics[width=0.40\textwidth]{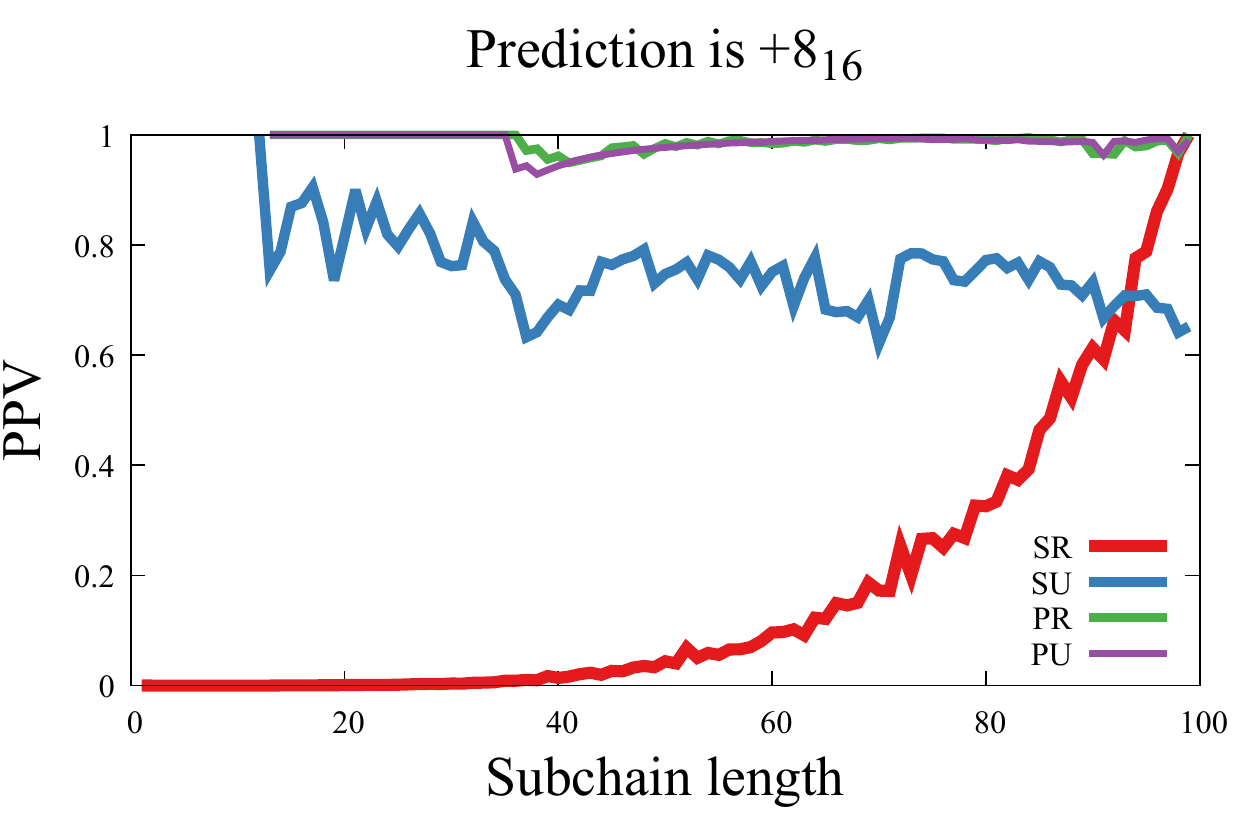}
\hfill{\ }

{\ }\hfill
\includegraphics[width=0.40\textwidth]{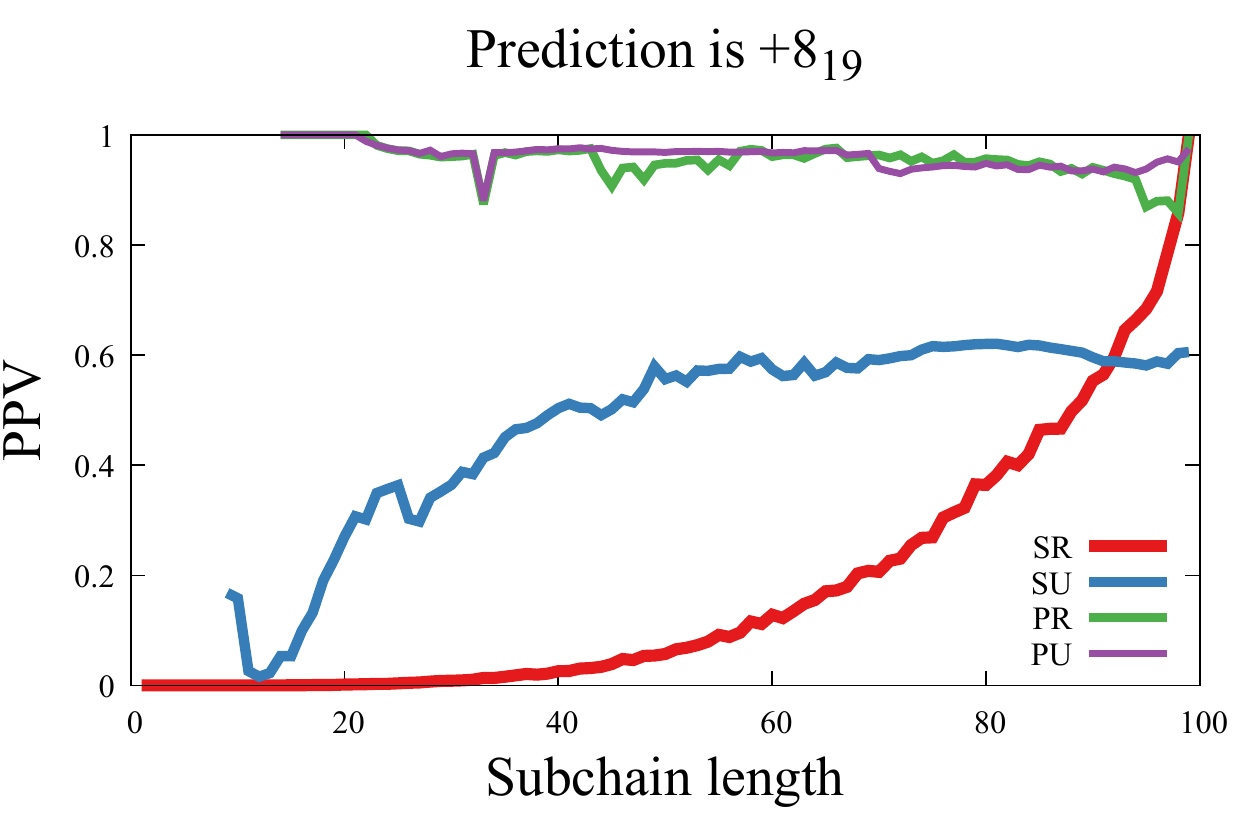}
\hfill
\includegraphics[width=0.40\textwidth]{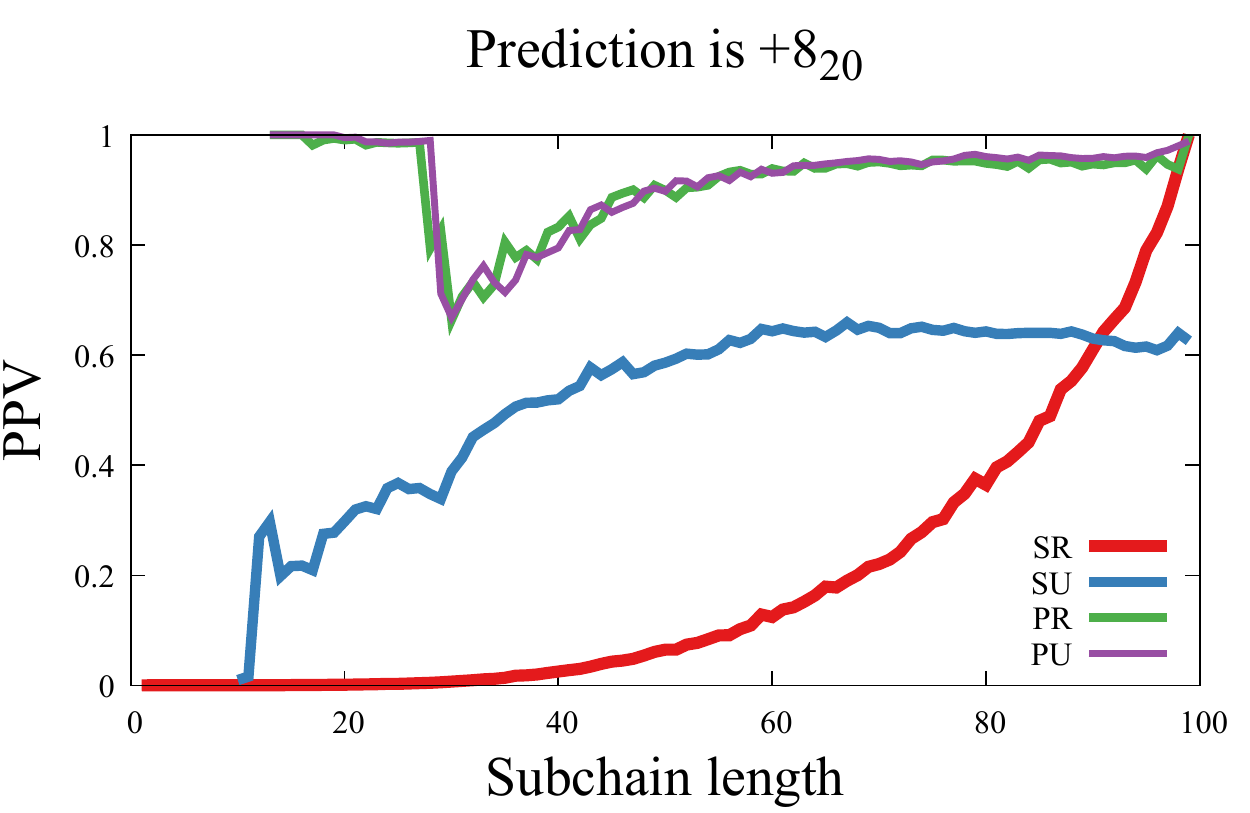}
\hfill{\ }

\pagebreak

\vfill

{\ }\hfill
\includegraphics[width=0.40\textwidth]{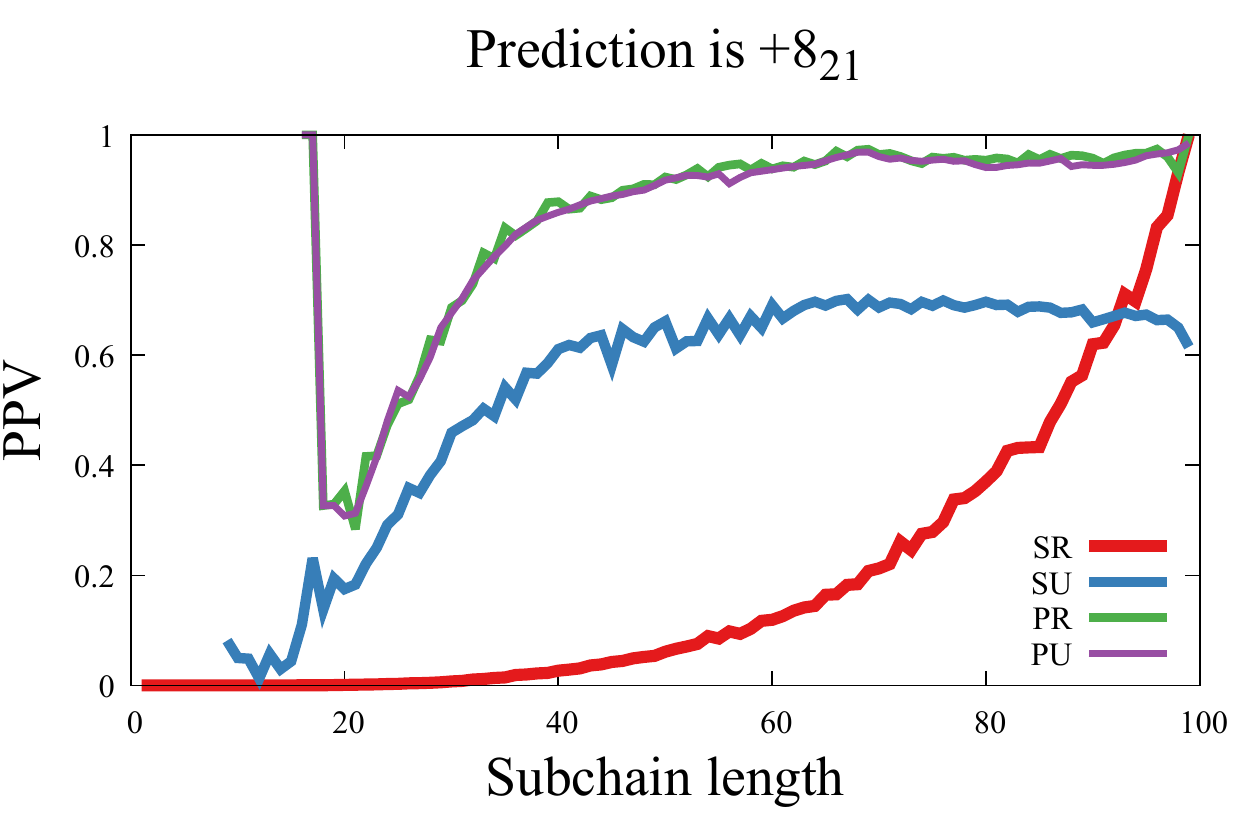}
\hfill
\includegraphics[width=0.40\textwidth]{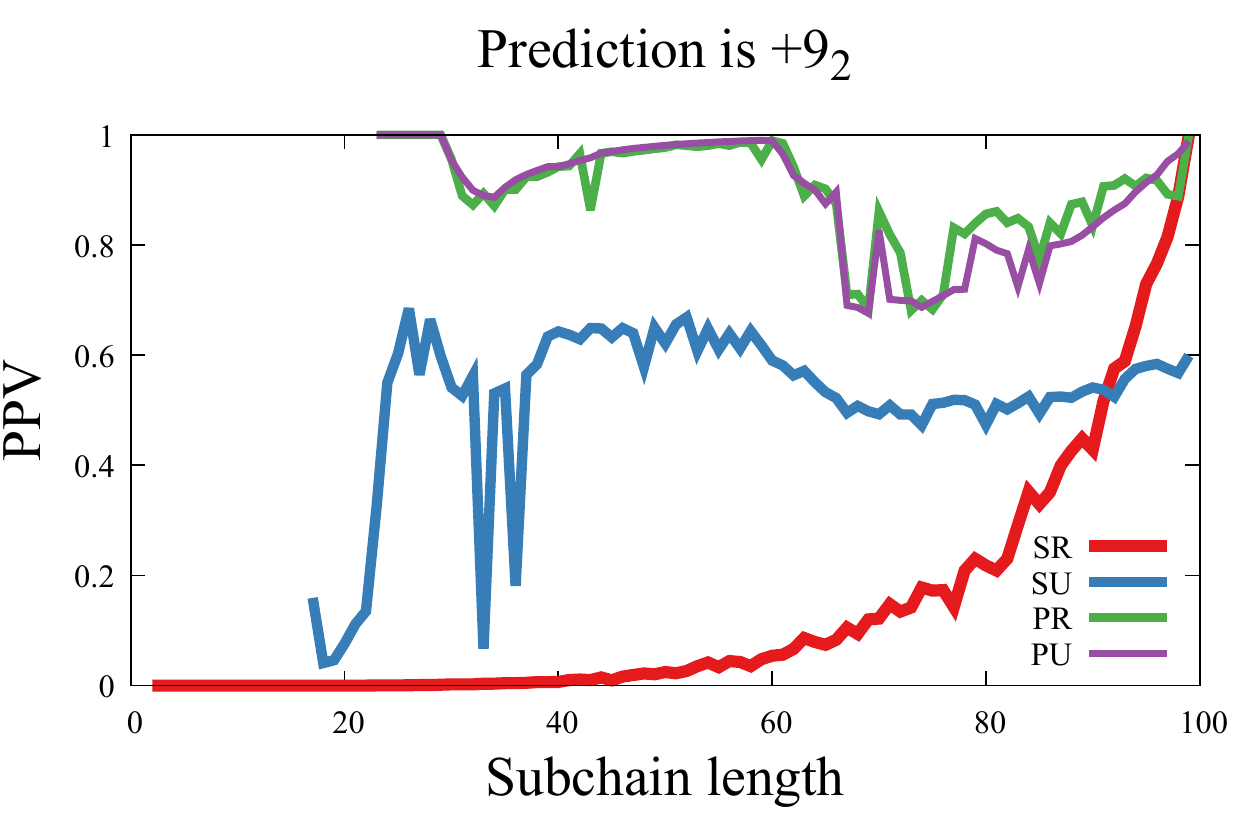}
\hfill{\ }

\vfill

{\ }\hfill
\includegraphics[width=0.40\textwidth]{p915.bayes.pdf}
\hfill
\includegraphics[width=0.40\textwidth]{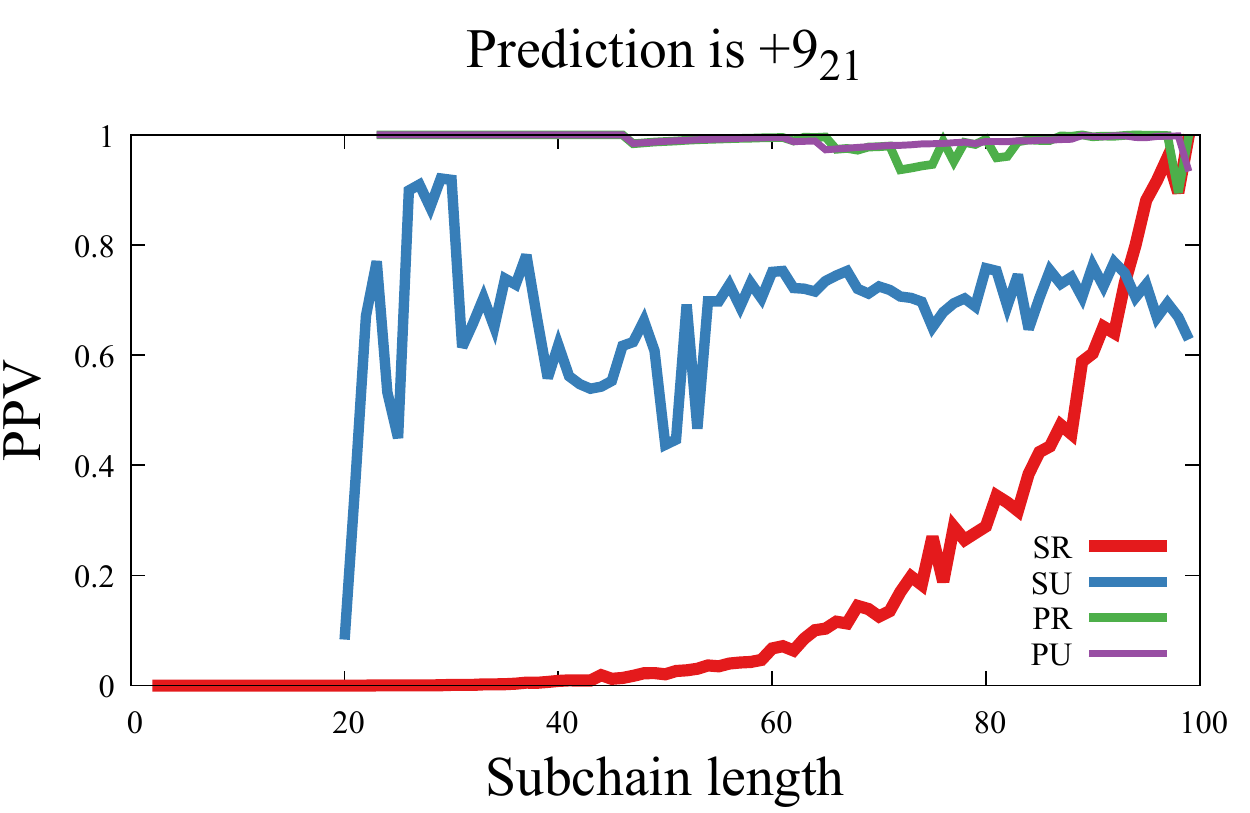}
\hfill{\ }

\vfill

{\ }\hfill
\includegraphics[width=0.40\textwidth]{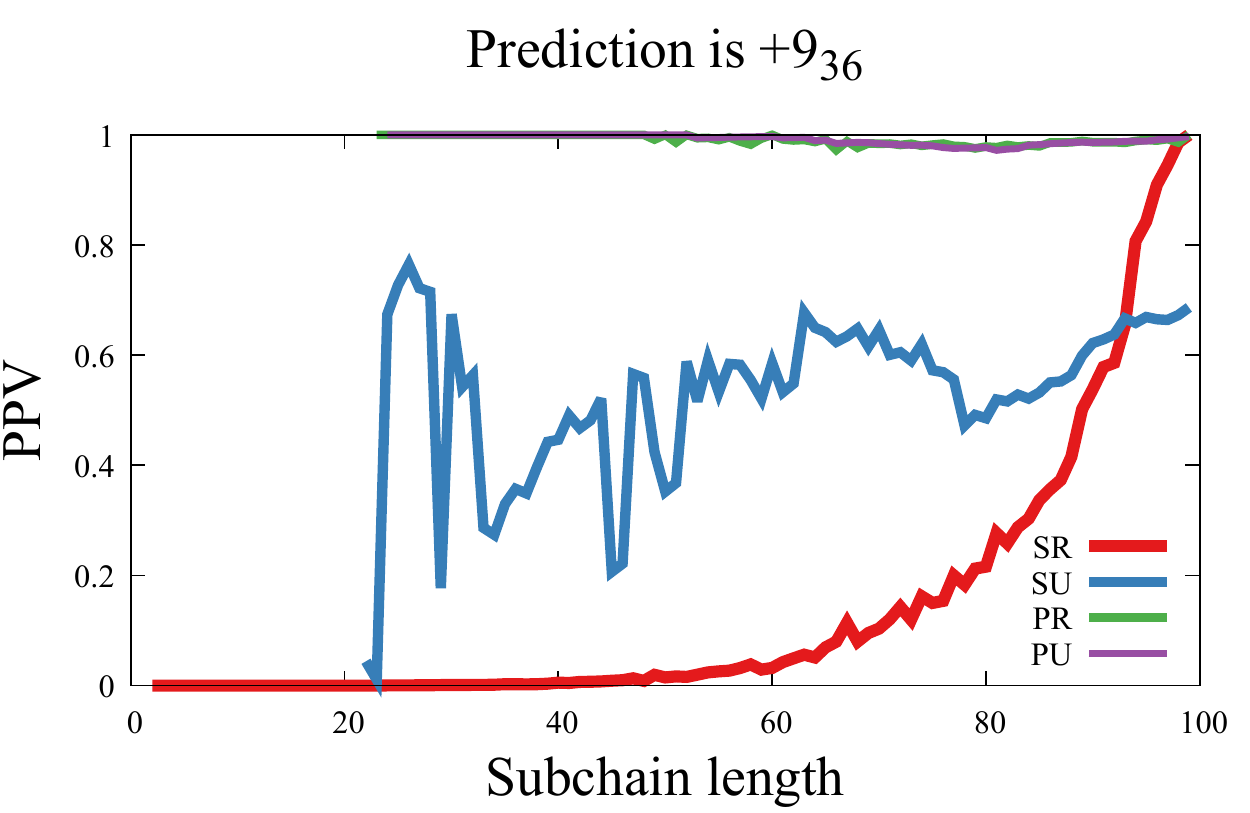}
\hfill
\includegraphics[width=0.40\textwidth]{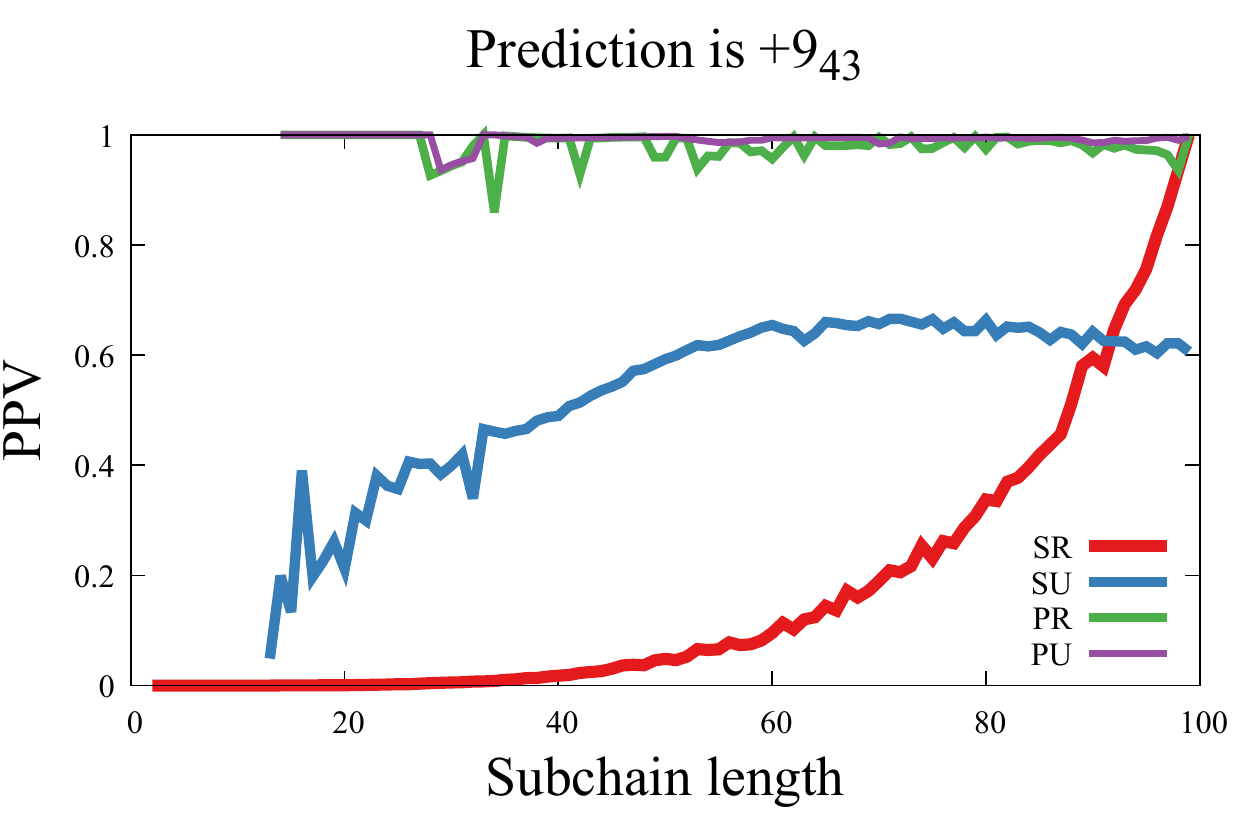}
\hfill{\ }

\vfill

{\ }\hfill
\includegraphics[width=0.40\textwidth]{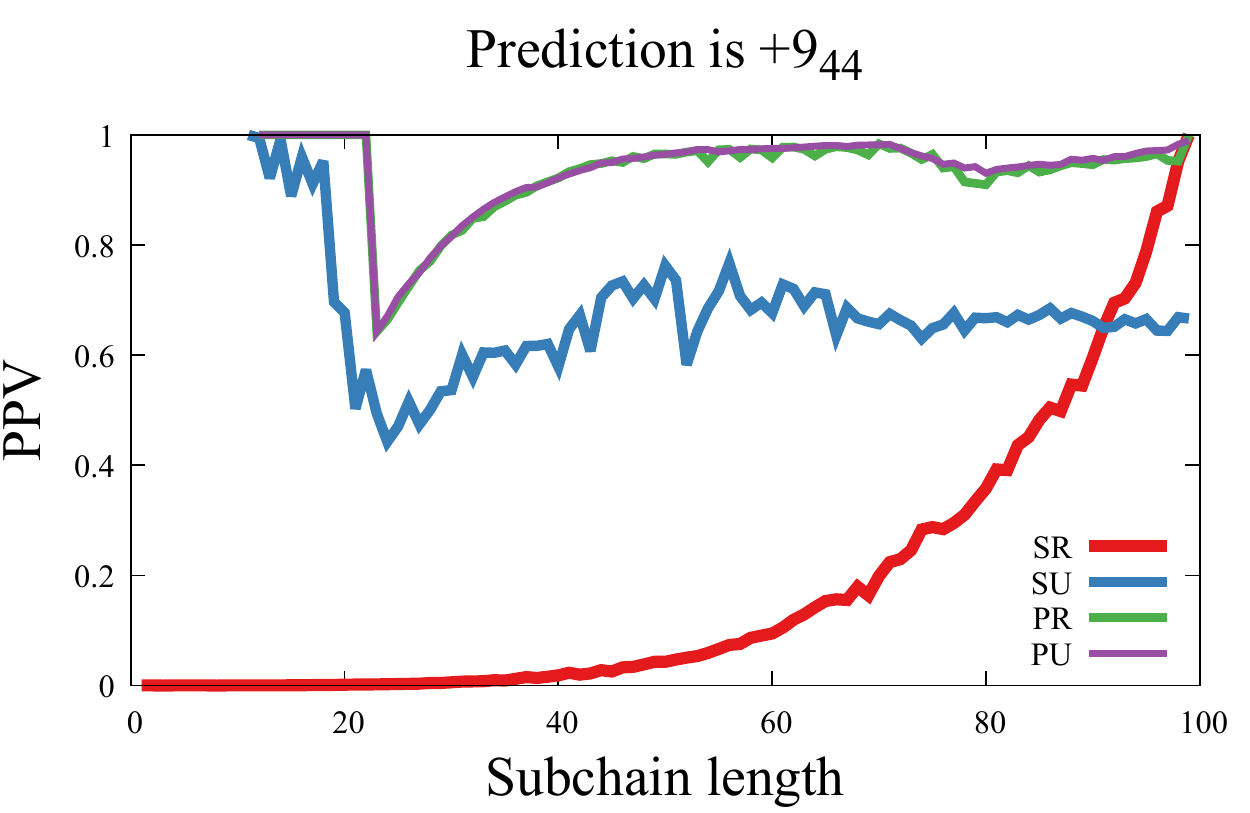}
\hfill
\includegraphics[width=0.40\textwidth]{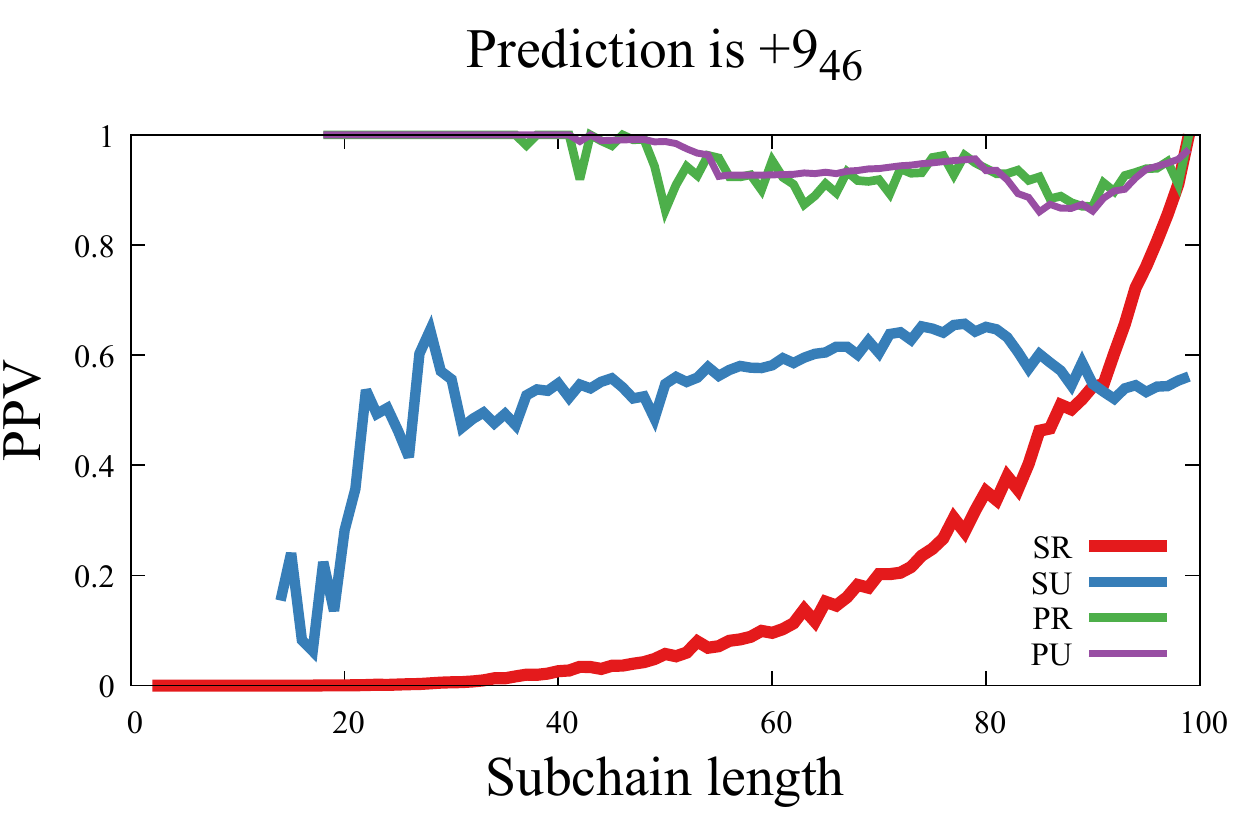}
\hfill{\ }

\vfill

\end{document}